\documentclass{article}
\usepackage[english]{babel}
\usepackage{amsfonts}
\usepackage{amscd}

\title
{\bf Skand theory and its applications. 
 \\ (A new look at non-well-founded sets)}

\author{Ju. T. Lisica\\
Mathematical Analysis and Function Theory Department\\
Russian University of Peoples' Friendship\\
Miklukho-Maklay str. 6\\
117198 Moscow, Russia\\
e-mail: jutlisica@yandex.ru}

\date{}
\begin{document}
\maketitle

\begin{abstract}
A new mathematical object called a {\it skand} is introduced, which turns out in general to be a non-well-founded set. Skands of finite lengths are ordinary well-founded sets, and skands of very long length (like the hyper-skand  of all ordinals) are hyper-classes. 

{\it Self-similar} skands are also considered, and they clarify the reflexivity of sets, i.e., the meaning of the relation $X\in X$; in particular, self-similar skands considered as non-well-founded sets are always reflexive, but not vice versa. The existence of self-similar skands shows at once that  Russell's well-known   paradox {\it is not a paradox} at all. The inconsistency of Russell's \grqq set" $R=\{X|\,\,X\not \in X\}$ is proved here not with the help of Russell's paradox (as it is traditionally given, which is incorrect), but via a simple method of the {\it maximality} ({\it universality})  of $R$ which goes back to Cantor and can be also applied to other set-theoretical paradoxes.

{\it Generalized skands} are also defined and a new look at the generalized skand-class of all ordinals is demonstrated. In particular, the last (class) ordinal called the {\it eschaton} is defined. 

The next application of skand theory is a description of all epsilon-numbers in the sense of Cantor. Another application is a generalized theory of one-dimensional continua of arbitrary powers and the construction of generalized real numbers as a non-Archimedean straight line of arbitrary power, and the introduction of the {\it absolute continuum} and the {\it absolute straight line} as the hyper-classes nearest to the class of sets.

\end{abstract}

\bigskip

\vspace{2cm}
{\bf 0. Epigraphs}

 \parindent=5,5cm A parson had a hound-dog,

One he loved a lot.

It ate a piece of mutton,

For which he had it shot.

He buried the hound,

Then wrote on its mound,

That

A parson had a hound-dog, etc. 

({\it ad infinitum} and {\it ad imum}).

\smallskip

[{\it Free translation of a Russian children's ditty}.]

\bigskip

\parindent=2,8cm Once a four-year-old son returned home from kindergarten, 

\parindent=2,8cm where he had been told that his father was a mathematician.  

\parindent=2,8cm When the son saw his father again, he asked him:

\grqq Is it true, Daddy, that you are a mathematician?"

\grqq Yes, sonny, it is," was the answer.

\grqq Well," responded the son, \grqq can you count to the last 

\parindent=2,8cm number?"

\grqq Ummm...ummmm," mumbled the father, stumped.

\smallskip

[Dialogue with a child which is in fact a problem of 

mathematical eschatology: 

\grqq What  is the 
'$\acute{\varepsilon}\sigma\chi\alpha\tau o\nu$ or {\it ad imum}?"]

\bigskip

\parindent=2,4cm \grqq A  soul is only a {\it skand}, i.e., an accidental aggregation of being".

\smallskip

[{\it From Buddhist doctrine.}]

\bigskip

\parindent=3,0cm \grqq The content of a concept diminishes as its extension 

increases; if its extension becomes all-embracing, its 

content must vanish altogether".

\smallskip

[{\it Gottlob Frege, \grqq The Foundations of Arithmetic".}]

\bigskip

\parindent=0 cm
{\bf 1. Introduction}

\parindent=0,5cm
We are going to clarify the notion of {\it reflexivity} in Set Theory, i.e., the meaning of a binary relation $X\in Y$ in the case when $Y=X$, that is $X\in X$. As to Russell himself, the relation $X\in X$ \grqq must be always meaningless" \cite{l55}, p. 81, since he was seriously frightened by this paradox, which was later on named after him: \grqq Thus $X\in X$ was held to be meaningless, because $\in$ requires that the relatum should be a class composed of objects which are of the type of the referent" (\cite{l47}, Chap. X, p. 107). Moreover, he went further and concluded that $X\not\in X$ \grqq must be always meaningless", too. (\cite{l55}, p. 81). He wrote: \grqq If $\alpha$ is a class, the statement \grq $\alpha$ is not a member of $\alpha$' is always meaningless, and there is therefore no sense in the phrase \grq the class of those classes which are not members of themselves'$\,\,$" (\cite{l55}, p. 66). In particular,  Russell wrote: \grqq A class consisting of only one member must not be identical with that one member". And he added immediately \grqq $X=\{X\}$ must be absolutely meaningless, not simply false" (\cite{l55}, p. 81). In this paper we shall see when Russell was right and when  was he not.

$$
$$

At the beginning we start out  within a von Neumann-Bernays-G\"{o}del-type set theory ($NBG$ for short) which includes  the axiom of choice ${\bf C}$  and the axiom of  foundation ${\bf FA}$. The basic set theory can be with individuals (called sometimes atoms or ur-elements, i.e., mathematical objects which are not sets or classes and which have no members), or it can be without them; it does not matter. In the latter case the only individual is just the empty set $\{\}=\emptyset$ and Set Theory in this case is often called the theory of \grqq pure sets". 

We denote  
   the class of all sets by  ${\bf V}[{\cal U}]$  in the set theory with the class ${\cal U}$ of  individuals (atoms, ur-elements) and ${\bf V}$ in the theory of \grqq pure" sets, i.e., when ${\cal U}=\emptyset$, respectively. In $NBG$ ${\bf V}[{\cal U}]$ there turns out to be the class ${\bf WF}$ of all well-founded sets. The {\it universal} class ${\bf U}$ is the union ${\bf V}[{\cal U}]\cup{\cal U}$ of ${\bf V}[{\cal U}]$ and ${\cal U}$. Note that ${\bf V}[{\cal U}]\cap{\cal U}=\emptyset$. 
   
   We also denote the class of all ordinals by ${\bf On}$  and the class of all cardinals by ${\bf Card}$. 
 
Then we shall consider $NBG^-$, i.e., $NBG$ without the axiom of foundation, and instead of the axiom of  choice we use the axiom ${\bf N}$ of von Neumann, ${\bf V}[{\cal U}]^-\approx{\bf On}$; i.e., these classes are bijective, which in the absence of foundation is stronger than choice.

 There is an important distinction in $NBG$ as well as in $NBG^-$ Set Theories: the distinction between {\it sets} and {\it proper classes}, or \grqq small" classes and \grqq large" classes. There are two ways to distinguish sets and proper classes in addition to axioms for sets. The first of them is the following: a subclass $X\subset{\bf U}$ of ${\bf U}$ is a set if and only if there exists a one-element object (singleton) $\{X\}\in{\bf U}$; the second one: a subclass $X\subset{\bf U}$ of ${\bf U}$ is a set if and only if there is no bijection $X$ on ${\bf U}$. Otherwise, $X\subseteq{\bf U}$ is not a set but a proper subclass of ${\bf U}$ (further, in short, a proper class). Moreover, all proper subclasses of ${\bf U}$ are bijective to each other.  Notice also that when we deal with proper classes we speak in the language of their elements but not of them as wholes or as units. In other words, sets are arguments (elements of ${\bf V}[{\cal U}]$), and classes are extensions of some predicates.

Recall that a set $X$  is {\it well-founded} (or ordinary) if every $\in$-descending chain in $X$ is finite, i.e., for each $x_0\in X$, every $\in$-descending chain starting with $x_0$ can be at most the following one: $x_0\owns x_1\owns x_2\owns...\owns x_n$, where  $x_n$ is some individual or the empty set; otherwise it is  {\it non-well-founded} or \grqq extraordinary", i.e., there exist infinite $\in$-descending chains.  The distinction between well-founded and non-well-founded sets was first articulated by Mirimanoff \cite{l2} and in his terminology the distinction was between \grqq ordinary" and \grqq extraordinary" sets. Later on, von Neumann \cite{l90} proposed an  axiom of regularity (\grqq Restrictive Axiom" in \cite{l321} and \grqq Axiom der Fundierung" in  \cite{l221}, i.e. the Foundation Axiom ${\bf FA}$) which excluded Mirimanoff's extraordinary sets, because according to the axiom any \grqq descending" sequence terminates, i.e., reaches its bottom or \grqq foundation".

The restriction axiom as ${\bf FA}$ was very important, first of all, for completeness of the extensionality axiom ${\bf Ext}$, because in $NBG[{\cal U}]^-$ it is impossible to prove in general that for different sets $X$ and $Y$ the two-element set $Z=\{X,Y\}$, which always exists by the axiom of pairing, $Z=X$, or $Z=Y$, or $Z$ is different from $X$ and $Y$. Only by means of ${\bf FA}$ can one prove that $\{X,Y\}$ is different from $X$ and $Y$. Secondly, ${\bf FA}$ avoids {\it vicious circle} phenomena, i.e., there is no set $X$ such that $X\owns X_1\owns...\owns X_n\owns X$; in particular, the reflexive sets $X\in X$, which appeared to be a source of paradoxes (which was actually not true), e.g.,  Russell, with reference to H. Poincar\'{e} wrote: \grqq An analysis of the paradoxes to be avoided shows that they all result from a certain kind of vicious circle"  \cite{l55}, p. 39. And Russell had formulated his famous \grqq vicious-circle principle" as follows: \grqq Whatever involves {\it all} of a collection must not be one of the collection" \cite{l55}, p. 40. Thus, with the help of  the {\bf Foundation Axiom},  Set Theory has been succesfully developed and the \grqq vicious-circle principle" has been satisfied. 

Nevertheless, many real problems concern circular phenomena in logic (e.g., the treatment of Liar-like paradoxes); linguistics, computer science, graph theory, game theory, streams, etc., all need  circular models which lie beyond the universe of well-founded sets.

From the early 20th century, many authors actually proposed their own {\it anti-foundation axioms} (${\bf AFA}$ for short) enriching and extending the Well-Founded Universe by the ${\bf AFA}$-Universe, e.g., \cite{l91}, \cite{l31}, \cite{l90}, \cite{l92}, \cite{l94} and others. All of them proposed theories of possibly non-well-founded sets which are consistent, assuming that $NBG$ is consistent. Nevertheless, the four axiom systems mentioned are non-comparable, and each one differs from the others in the strengthening of the extensionality criterion for set equality.

In the present paper we introduce a  {\it new object} called  a  {\it skand} ({\it Sanskrit}: jump, skip). Skands are \grqq definite and separate" objects (\grqq bestimmten wohlunterschiedenen Objekten" according to Cantor \cite{l1}, p. 481) or elements of the class ${\bf V}[{\cal U}]^-$ of all sets in $NBG^-$ and they also enrich ${\bf V}[{\cal U}]$; i.e., they can be well-founded or non-well-founded sets. Moreover, skands are essential extensions of some (not all) of Mirimanoff's extraordinary sets. 
 
 The class of all skands generates  a proper class ${\bf V}[{\cal U}]^{(1)}$ which is a subclass of the class  ${\bf V}[{\cal U}]^-$ by forming objects $X$ whose elements are ordinary sets or skands such that $\{X\}$ exists.  Clearly, ${\bf V}[{\cal U}]\stackrel{def}{=}{\bf V}[{\cal U}]^{(0)}\subset{\bf V}[{\cal U}]^{(1)}\subset{\bf V}[{\cal U}]^-$. Moreover, we successively continue such a process of enrichment for each ordinal number $\alpha$ and obtain the following embeddings: ${\bf V}[{\cal U}]^{(0)}\subset{\bf V}[{\cal U}]^{(1)}\subset...\subset{\bf V}[{\cal U}]^{(\alpha)}\subset...\subset {\bf V}[{\cal U}]^{\Omega}\subset{\bf V}[{\cal U}]^-$, where ${\bf V}[{\cal U}]^{\Omega}=\bigcup\limits_{\alpha\in{\bf On}}{\bf V}[{\cal U}]^{(\alpha)}$. We do not say that this cumulative hierarchy exhausts ${\bf V}[{\cal U}]^-$, i.e., ${\bf V}[{\cal U}]^{\Omega}={\bf V}[{\cal U}]^-$ but it makes ${\bf V}[{\cal U}]^-$ more structural. $\Box$
 
 \bigskip

\parindent=0 cm
{\bf 2. The notion of a skand}

\parindent=0,5cm
Let $(\alpha_0,\alpha)\stackrel{def}{=}\{\alpha'\in {\bf  On}|\,\,\alpha_0\leq\alpha'<\alpha\}$ be a segment of ${\bf On}$; if $\alpha_0=0$ and $\alpha\geq 1$, then we call $(0,\alpha)$ an initial segment of ${\bf On}$.

{\bf Definition 1}. Consider  a system of embedded curly braces or well-ordered set of embedding pairs of curly braces $\{_{\alpha_0}\{_{\alpha_0+1}...\{_{\alpha'}..._{\alpha'}\}..._{\alpha_0+1}\}_{\alpha_0}\}$, indexed by ordinals $\alpha'\in(\alpha_0,\alpha)$. We call it a {\it trivial skand} of length $l=\alpha-\alpha_0$ and denote it by ${\bf e}_{(\alpha_0,\alpha)}$.  It can be also called an {\it empty skand} because for each $\alpha'\in(\alpha_0,\alpha)$ the set of all elements between two neighboring opening braces $\{_{\alpha'}\,\,\{_{\alpha'+1}$ or possibly between $\{_{\alpha'}\,\,_{\alpha'}\}$, if $\alpha'=\alpha-1$,  is empty. In other words, each $\alpha'$-component ${\bf e}_{\alpha'}$ of ${\bf e}_{(\alpha_0,\alpha)}$ is empty. By a non-trivial {\it skand} $X_{(\alpha_0,\alpha)}$ of length $l=\alpha-\alpha_0$ we understand   a non-trivial system of embedded curly braces, indexed by ordinals $\alpha'\in(\alpha_0,\alpha)$; i.e., there is at least one index $\alpha'\in(\alpha_0,\alpha)$, and elements $x_0,x_1,x_2,...,x_\lambda,...$ of ${\bf U}$  such that the $\alpha'$-component $X_{\alpha'}\stackrel{def}{=}\{_{\alpha'}\,x_0,x_1,x_2,...,x_\lambda,...,\{_{\alpha'+1}$ (or $\{_{\alpha'}\,x_0,x_1,x_2,...,x_\lambda,...\,_{\alpha'}\}$, if $\alpha'=\alpha-1$) of $X_{(\alpha_0,\alpha)}$ is a set, i.e., $\{x_0,x_1,x_2,...,x_\lambda,...\}\in {\bf V}[{\cal U}]$. 

For simplicity we shall omit indexes $\alpha'$ of braces or elements in the cases when it is clear what they are, e.g., for $\{_{\alpha'}x^{\alpha'}_0,x^{\alpha'}_1,x^{\alpha'}_2,...,x^{\alpha'}_{\lambda_{\alpha'}},...,\{_{\alpha'+1}$ or possibly $\{_{\alpha'}x^{\alpha'}_0,x^{\alpha'}_1,x^{\alpha'}_2,...,x^{\alpha'}_{\lambda_{\alpha'}},..._{\alpha'}\}$ in the case $\alpha'=\alpha-1$, we write 

\parindent=0 cm 
$\{x^{\alpha'}_0,x^{\alpha'}_1,x^{\alpha'}_2,...,x^{\alpha'}_{\lambda_{\alpha'}},...,\{$ and $\{x^{\alpha'}_0,x^{\alpha'}_1,x^{\alpha'}_2,...,x^{\alpha'}_{\lambda_{\alpha'}},...\}$ in the case $\alpha'=\alpha-1$, respectively; or even simplify to $\{x_0,x_1,x_2,...,x_\lambda,...,\{$ and $\{x_0,x_1,x_2,...,x_\lambda,...\}$, respectively. All the more indexes $\alpha'$ of braces are conditional, e.g., if $X_{(0,\alpha)}$ is a skand, then $Y_{(0,\alpha)}=\{X_{(0,\alpha)}\}$ whose components are $Y_0=\emptyset$, $Y_1=X_0$,$\,\,$$Y_2=X_1$,$\,\,$...,$\,\,$$Y_{n+1}=X_n$,$\,\,$...$\,\,$, and $Y_{\alpha'}=X_{\alpha'}$, for all $\omega\leq\alpha'<\alpha$, is a skand, too.

\parindent=0,5cm 
 For the purpose of interpreting skands as sets, we  write a comma before the second open brace of non-trivial component $X_{\alpha'}$ of $X_{(\alpha_0,\alpha)}$,   i.e.,
 
  \parindent=0 cm
  $\{_{\alpha'}\,x_0,x_1,x_2,...,x_\lambda,...,\,_{\alpha'+1}\{$ or simply $\{x_0,x_1,x_2,...,x_\lambda,...,\{$. It says that braces are not only syntactical but also semantic in the definition of skand.

\parindent=0,5 cm
Thus, a general form of an arbitrary skand is the following:
\begin{equation}
\label{f023}
X_{(\alpha_0,\alpha)}=\{x^{\alpha_0}_0,x^{\alpha_0}_1,...,x^{\alpha_0}_{\lambda_{\alpha_0}},...,\{...\,\{x^{\alpha'}_0,x^{\alpha'}_1,...,x^{\alpha'}_{\lambda_{\alpha'}},...,\{...\}\}...\},
\end{equation}
where components $X_{\alpha'}=\{x^{\alpha'}_0,x^{\alpha'}_1,...,x^{\alpha'}_{\lambda_{\alpha'}},...,\{$ are sets, i.e., $\{x^{\alpha'}_0,x^{\alpha'}_1,...,x^{\alpha'}_{\lambda_{\alpha'}},...\}\in {\bf V}[{\cal U}]$ or empty, i.e., $X_{\alpha'}=\{\{$; if $\alpha$ is not a limit ordinal, then the last component $X_{\alpha-1}$ is an ordinary set, i.e., $\{x^{\alpha'}_0,x^{\alpha'}_1,...,x^{\alpha'}_{\lambda_{\alpha'}},...\}\in {\bf V}[{\cal U}]$ or empty $\{\}$.

If  all components $X_{\alpha'}$, $\alpha_0\leq\alpha'<\alpha$, of a skand $X_{(\alpha_0,\alpha)}$ are equal to the same set, e.g., $X=\{x_0,x_1,...,x_\lambda,...\}$, we shall denote this skand in a shorter way by $X_{(\alpha_0,\alpha)}(X)$; in particular, when $X=\{\gamma\}$, i.e., a one-element set $X$, we simplify the notation  to $X_{(\alpha_0,\alpha)}(\gamma)$. If $\alpha_0<\alpha_1<\alpha$, then together with $X_{(\alpha_0,\alpha)}$ we shall denote by $X_{(\alpha_1,\alpha)}$ the skand whose $\alpha'$-components $X_{\alpha'}$, $\alpha_1\leq\alpha'<\alpha$, are the same as those of $X_{(\alpha_0,\alpha)}$.

{\bf Remark 1.} Actually, between each pair of braces $\{\{$ or $\{\}$ of a skand $X_{(\alpha_0,\alpha)}$ there are different elements (sets or individuals), or their lack, which form a set; notice that the order in which the members of sets are written does not matter. In our notation we use $X_{\alpha'}=\{x_0,x_1,...,x_\lambda,...,\{$ or $X_{\alpha'}=\{x_0,x_1,...,x_\lambda,...\}$ only for simplicity of writing, since with the axiom of choice we can well-order any set $X$, $\lambda<\kappa$, for some ordinal $\kappa$, and obtain such notation.  
 It is clear that the definition of a skand does not depend on orderings  of elements of $X_{\alpha'}$, $\alpha_0\leq\alpha'<\alpha$, which may be different, or notations, which may also be alternative. Thus, an arbitrary skand is \grqq only an accidental aggregation" of sets or individuals each of which figuratively \grqq skips"  into its own place, $\{\{$ or $\{\}$, in the system of well-ordered embedded curly braces; meanwhile there can be empty places as well as a set, or an individual, and can \grqq skip" into different places and even be at all places at once; i.e., elements of each component are always different and at the same time it may happen that some or even all elements of different components may be equal. Finally, skands look like elements of an arbitrary direct product of sets with two differences: 1) components (coordinates) of elements of a direct product are elements of sets and components of a skand can be not only elements of sets but sets of elements or empty sets; 2) the notions of equality of elements of direct products and skands are absolutely different.

\parindent=0,5cm {\bf Definition 2}.
 Two skands $X_{(\alpha_0,\alpha)}$ and $Y_{(\beta_0,\beta)}$ are called {\it equal} if the segments $({\alpha_0,\alpha})$ and $({\beta_0,\beta})$ are isomorphic as well-ordered sets, i.e., their \grqq similarity" is given by $\varphi:(\alpha_0,\alpha)\rightarrow(\beta_0,\beta)$, and the corresponding  $\alpha'$- and $\beta'$-components  $X_{\alpha'}=\{_{\alpha'}x^{\alpha'}_1,x^{\alpha'}_2,...,x^{\alpha'}_\lambda,...,\{_{\alpha'+1}$ and $Y_{\beta'}=\{_{\beta'}y^{\beta'}_1,y^{\beta'}_2,...,y^{\beta'}_\lambda,...,\{_{\beta'+1}$ as well as the last components $X_{\alpha'}=\{_{\alpha'}x^{\alpha'}_1,x^{\alpha'}_2,...,x^{\alpha'}_\lambda,..._{\alpha'}\}$ and 

\parindent=0 cm 
$Y_{\beta'}=\{_{\beta'}y^{\beta'}_1,y^{\beta'}_2,...,y^{\beta'}_\lambda,..._{\beta'}\}$, if $\alpha'=\alpha-1$ and $\beta'=\beta-1$,   are equal as sets in $NBG$; i.e., $\{x^{\alpha'}_1,x^{\alpha'}_2,...,x^{\alpha'}_\lambda,...\}=\{y^{\beta'}_1,y^{\beta'}_2,...,y^{\beta'}_\lambda,...\}$, $\alpha_0\leq\alpha'<\alpha$ and $\beta_0\leq\beta'<\beta$, respectively, where $\beta'=\varphi(\alpha')$. 

\parindent=0,5cm 

{\bf Example 1.} $X_{(\alpha_0,\alpha)}=\{\{1,\{2,\{3,\{...\}\}\}\}\}$, $Y_{(\alpha_0,\alpha)}=\{1,\{\{2,\{3,\{...\}\}\}\}\}$.

Skands $X_{(\alpha_0,\alpha)}$ and $Y_{(\alpha_0,\alpha)}$ are not equal because  the $\alpha_0$-component $X_{\alpha_0}$ of $X_{(\alpha_0,\alpha)}$  is empty and  the $\alpha_0$-component $Y_{\alpha_0}$ of $Y_{(\alpha_0,\alpha)}$ consists of one element, which is equal to $1$. These skands have the same components, but in different order.

{\bf Remark 2}.
If $X_{(\alpha_0,\alpha)}=\{x^0_0,x_1^0,...,x^0_{\lambda_0},...,\{x^1_0,x^1_1,...,x^1_{\lambda_1},...,\{...\}\}\}$ is  an arbitrary skand, then actually it can be considered as a set in ${\bf V}[{\cal U}]^-$ whose elements are sets or individuals $x^0_0,x_1^0,...,x^0_\lambda,...$ and the skand $X_{(\alpha_0+1,\alpha)}=\{x^1_0,x^1_1,...,x^1_{\lambda_1},...,\{...\}\}$; i.e., $X_{(\alpha_0,\alpha)}=\{x^0_0,x_1^0,...,x^0_{\lambda_0},...,X_{(\alpha_0+1,\alpha)}\}$. In particular, for $l\geq\omega$  a trivial skand ${\bf e}_{(\alpha_0,\alpha)}$ can be considered as a set whose  only element is this set itself, i.e., ${\bf e}_{(\alpha_0,\alpha)}=\{{\bf e}_{(\alpha_0,\alpha)}\}$ because, by Definition 2, ${\bf e}_{(\alpha_0,\alpha)}={\bf e}_{(\alpha_0+1,\alpha)}$. So, Russell was not right saying that \grqq$X=\{X\}$ must be absolutely meaningless, not simply false".

Clearly, if $\alpha=\alpha_0+n$, where $n=1,2,...$ are natural numbers, then $X_{(\alpha_0,\alpha)}=\{x^0_0,x_1^0,...,x^0_{\lambda_0},...,\{x^1_0,x^1_1,...,x^1_{\lambda_1},...,\{...,\{x^{n-1}_0, x^{n-1}_1,...,x^{n-1}_{\lambda_{n-1}},...\}\}\}\}$ is an ordinary well-founded set whose elements are sets $x^0_0,x_1^0,...,x^0_{\lambda_0},...,X_{(\alpha_0+1,\alpha)}$ of ${\bf V}[{\cal U}]$.

\parindent=0,5cm {\bf Example 2}. $X_{(0,\omega)}=\{a_0,\{a_1,\{a_2,\{...\}\}\}\}$, where $a_i$, $i=0,1,2,...$, are individuals or sets in ${\bf U}$, where $\omega=\omega_0$ is the first infinite ordinal.

This is  a skand of length $\omega$ or a two-element set in ${\bf V}^-$; i.e., $X=\{a_0,X_{(1,\omega)}\}$, where $X_{(1,\omega)}=\{a_1,\{a_2,\{...\}\}\}$. It is a typical example of a non-well-founded set, or an extraordinary set in the sense of Mirimanoff \cite{l2}, \cite{l3}, \cite{l4}, where $a_i$, $i=0,1,2,...$, are individuals; e.g., non-negative integers $0,1,2,...$; i.e. $X_{(0,\omega)}=\{0,\{1,\{2,\{...\}\}\}\}$, or examples of Mirimanoff's circular extraordinary sets of period $n$, where $a_i=a_{i+n}$, $i=0,1,2,...$, and $n=1,2,...$ is a fixed natural number; e.g., $X_{(0,\omega)}=\{0,\{1,\{0,\{1,\{...\}\}\}\}\}$, where $n=2$, $a_0=0$ and $a_1=1$.

The following  is a more general example of a circular set.

{\bf Example 3}. Consider
\begin{equation}
\label{f165}
X_{(0,\omega)}=\{a_0,a_1,...,a_\lambda,...,\{b_0,b_1,..b_\mu,...,\{a_0,a_1,..a_\lambda,...,\{b_0,b_1,...,b_\mu,...,\{...\}\}\}\}\},
\end{equation} 
 where $a_\lambda$ and $b_\mu$  are individuals or sets in ${\bf U}$ such that $\{a_0,a_1,...,a_\lambda,...\}$ and $\{b_0,b_1,..b_\mu,...\}$ are elements of ${\bf V}[{\cal U}]$, $0\leq\lambda\leq\kappa$ and $0\leq\mu\leq\nu$, respectively. In other words, an even component  $X_{2n}$ is the set $\{a_0,a_1,...,a_\lambda,...\}$ and an odd component $X_{2n+1}$ is the set $\{b_0,b_1,...,b_\mu,...\}$, $0\leq n<\omega$.

This skand can be considered as circular sets, i.e., $X_{(0,\omega)}\ni X_{(1,\omega)}\ni X_{(0,\omega)}$ and $X_{(1,\omega)}\ni X_{(0,\omega)}\ni X_{(1,\omega)}$ because $X_{(0,\omega)}=\{a_0,a_1,..a_\lambda,...,X_{(1,\omega)}\}$, $X_{(1,\omega)}=\{b_0,b_1,..b_\mu,...,X_{(2,\omega)}\}$, $X_{(2,\omega)}=\{a_0,a_1,..a_\lambda,...,X_{(3,\omega)}\}$ and, by Definition 2, $X_{(0,\omega)}= X_{(2,\omega)}$ and $X_{(1,\omega)}= X_{(3,\omega)}$.

{\bf Example 4}. Consider two skands 
\begin{equation}
\label{f166}
\begin{array}{cc}
X_{(0,\omega 2)}=\{a_0,a_1,...,a_\lambda,...,\{b_0,b_1,...,b_\mu,...,\{a_0,a_1,..a_\lambda,...\\
\{a_0,a_1,...,a_\lambda,...,\{b_0,b_1,..b_\mu,\{...\}\}\}...\}\}\},
\end{array}
\end{equation}
 where $a_\lambda$ and $b_\mu$  are individuals or sets in ${\bf U}$,  $0\leq\lambda\leq\kappa$ and $0\leq\mu\leq\nu$, respectively, where  even components $X_{2\tau}$ are equal to the fixed set $\{a_0,a_1,...,a_\lambda,...\}$ and  odd components $X_{2\tau+1}$ are equal to the other fixed set $\{b_0,b_1,...,b_\mu,...\}$, $0\leq\tau<\omega 2$, and
 \begin{equation}
\label{f167}
\begin{array}{cc}
Y_{(0,\omega 2)}=\{a_0,a_1,...,a_\lambda,...,\{b_0,b_1,...,b_\mu,...,\{a_0,a_1,..a_\lambda,...\\
\{b_0,b_1,...,b_\mu,\{a_0,a_1,..a_\lambda,\{...\}\}\}...\}\}\}
\end{array}
\end{equation}
 which can be obtained from $X_{(0,\omega 2)}$ by changing  $2\tau$-components $X_{2\tau}$ on the set $\{b_0,b_1,...,b_\mu,...\}$ and $2\tau+1$-components $X_{2\tau+1}$ on  the set $\{a_0,a_1,...,a_\lambda,...\}$,  for all $\omega\leq\tau<\omega 2$.

They can also be considered as circular sets; i.e., $X_{(0,\omega 2)}\ni X_{(1,\omega 2)}\ni X_{(0,\omega 2)}$ and $X_{(1,\omega 2)}\ni X_{(0,\omega)}\ni X_{(1,\omega)}$, because $X_{(0,\omega 2)}=\{a_0,a_1,..a_\lambda,...,X_{(1,\omega 2)}\}$, $X_{(1,\omega 2)}=\{b_0,b_1,..b_\mu,...,X_{(2,\omega 2)}\}$, $X_{(2,\omega 2)}=\{a_0,a_1,..a_\lambda,...,X_{(3,\omega 2)}\}$ and, by Definition 2,  $X_{(0,\omega 2)}= X_{(2,\omega 2)}$ and $X_{(1,\omega 2)}= X_{(3,\omega 2)}$ as well as  $Y_{(0,\omega 2)}\ni Y_{(1,\omega 2)}\ni Y_{(0,\omega 2)}$ and $Y_{(1,\omega 2)}\ni Y_{(0,\omega 2)}\ni Y_{(1,\omega 2)}$,  because 

\parindent=0 cm
$Y_{(0,\omega 2)}=\{a_0,a_1,..a_\lambda,...,Y_{(1,\omega 2)}\}$, $Y_{(1,\omega 2)}=\{b_0,b_1,..b_\mu,...,Y_{(2,\omega 2)}\}$, $Y_{(2,\omega 2)}=\{a_0,a_1,..a_\lambda,...,Y_{(3,\omega 2)}\}$ and, by Definition 2, $Y_{(0,\omega 2)}= Y_{(2,\omega 2)}$ and $Y_{(1,\omega 2)}= Y_{(3,\omega 2)}$. Nevertheless,  $X_{(0,\omega 2)}\not=Y_{(0,\omega 2)}$, as well as $X_{(1,\omega 2)}\not=Y_{(1,\omega 2)}$. 
 It is also clear that $X_{(0,\omega)}\not=X_{(0,\omega 2)}$ and $X_{(0,\omega)}\not=Y_{(0,\omega 2)}$. $\Box$

 $$
 $$

 \parindent=0,5cm
Note that every well-founded set $X=\{x_0,x_1,...,x_\lambda,...\}$ can be considered not only as a skand $X_{(0,1)}$ of length $1$ but also as many other different skands  of different finite lengths in general. Indeed, fix, for example, one element of $X$. Let it be $x_0$. Since $X$ is well-founded we choose a descending $\in$-chain $x_0\ni x_0^1\ni x_0^2\ni...\ni x_0^{n}$ such that $x_0^{n}$ is an individual or the empty set, and fix it. 
Then $X=X_{(0,n+1)}=\{x_1,x_2,...,x_\lambda,...,x_0\}$, where $x_0=\{a^1_1,a^1_2,...,a^1_{\lambda_1},...,a^1_0\}$ and $a^1_0=x^1_0=\{a_1^2,a_2^2,...,a^2_{\lambda_2},..., a^2_0\}$, $a^2_0=x^2_0=\{a_1^3,a_2^3,...,a^3_{\lambda_3},...,a^3_0\}$,..., $a_0^{n-1}=x_0^{n-1}=\{a^{n}_1,a^{n}_2,...,a_{\lambda_{n}}^{n},a_0^{n}\}$,   $a_0^{n}=x_0^{n}$; i.e., 

\parindent=0 cm
$X=X_{(0,n+1)}=\{x_1,x_2,...,x_\lambda,...,\{a^1_1,a^1_2,...,a^1_{\lambda_1},...,\{...,\{a^{n}_1,a^{n}_2,...,a_{\lambda_{n}}^{n},...,a_0^{n}\}\}\}\}$. 

\parindent=0,5cm
If for example, $x_0^n=\emptyset$, then there is another skand which gives the same set $X=X_{(0,n+2)}=\{x_1,x_2,...,x_\lambda,...,\{a^1_1,a^1_2,...,a^1_{\lambda_1},...,\{...,\{a^{n}_1,a^{n}_2,...,a_{\lambda_{n}}^{n},...,\{\}\}\}\}\}$.

 {\bf Example 5}. $X_{(0,\omega+1)}=\{a_0,\{a_1,\{a_2,\{...\{a_\omega\}...\}\}\}$, where $a_i$, $i=0,1,2,...,\omega$, are individuals or sets in ${\bf U}$.

This is the simplest example of non-well-founded set whose only two elements are $a_0$ and skand $X_{(1,\omega+1)}=\{a_1,\{a_2,\{...\{a_\omega\}...\}\}\}$. It was not considered by Mirimanoff because by remaining silent on the issue of an extraordinary set, i.e., $X\ni X_1\ni X_2\ni...$, Mirimanoff seems to suggest that it is in our terminology a skand of length  $\omega$ only. Moreover, all such skands and others of length $\omega+\beta$, where $\beta\in{\bf On}$, $\beta\geq 1$, are here {\it essential generalizations} and {\it extensions} of some of Mirimanoff's extraordinary sets. $\Box$

\bigskip

\parindent=0 cm
{\bf 3. The universe ${\bf U}^{(1)}$}
 
 \parindent=0,5cm
 We want to extend the universe ${\bf U}\stackrel{def}{=}{\bf U}^{(0)}$  as well as  the set theory $NBG\stackrel{def}{=}NBG^{(0)}$ to  the universe which we denote by ${\bf U}^{(1)}$ and the set theory $NBG^{(1)}$, respectively. For this purpose we denote by ${\cal U}^{(1)}$ the class of all skands $X_{(\alpha_0,\alpha)}$ such that $X_{(\alpha_0,\alpha)}\notin {\bf V}[{\cal U}]={\bf V}[{\cal U}]^{(0)}$, i.e., ${\cal U}^{(1)}$ is the class of all skands of length $l\geq\omega$. Note that we consider a skand $X_{(\alpha_0,\alpha)}$ as the set whose elements are elements of the $\alpha_0$-component $X_{\alpha_0}$ of $X_{(\alpha_0,\alpha)}$ and one more element $X_{(\alpha_0+1,\alpha)}$. Skands of finite length do not enrich ${\bf V}[{\cal U}]$; moreover, on the class of all skands of finite length we define the following equivalence relation: $X_{(\alpha_0,\alpha)}\sim Y_{(\beta_0,\beta)}$ iff sets $X_{\alpha_0}\cup\{X_{(\alpha_0,\alpha)}\}$ and $Y_{\beta_0}\cup\{Y_{(\beta_0,\beta)}\}$, $0<\alpha,\beta<\omega$, are equal in $NBG$.

Then denote by ${\bf V}[{\cal U}]^{(1)}$ the class of sets {\it generated} by ${\bf V}[{\cal U}]={\bf V}[{\cal U}]^{(0)}$ and ${\cal U}^{(1)}$. More precisely, 

$a$) all well-founded sets are members of ${\bf V}[{\cal U}]^{(1)}$, i.e., ${\bf V}[{\cal U}]^{(0)}\subset{\bf V}[{\cal U}]^{(1)}$; 

$b$) all skands $X_{(\alpha_0,\alpha)}\in{\cal U}^{(1)}$ considered as non-well-founded sets (elements $X_{(\alpha_0+1,\alpha)}$ are non-well-founded) are elements of ${\bf V}[{\cal U}]^{(1)}$; 

$c$) the results of all set-theoretic operations  on sets (i.e., unions, intersections, power-sets, images of sets under mappings, products, coproducts, function sets, inverse and direct limits, etc.) and all possible constructions or implicit sets, whose constituents are elements of ${\bf V}[{\cal U}]$ or of
 ${\cal U}^{(1)}$, are elements of ${\bf V}[{\cal U}]^{(1)}$.
 
One can see that for sets and subclasses of ${\bf V}[{\cal U}]^{(1)}$ all axioms of $NBG=NBG^{(0)}$ (except the Axiom of Extensionality for sets and classes) are satisfied. To obtain $NBG^{(1)}$, we add to $NBG$ the {\bf Axiom of Skand Existence},  which is accepted as Definition 1, and refine the axiom of extensionality for sets and classes  in the following way:
 
 {\bf Axiom of Strong Extensionality}. Two sets (resp., subclasses) $X^{(1)}$ and $Y^{(1)}$ in ${\bf V}[{\cal U}]^{(1)}$ (resp., of ${\bf V}[{\cal U}]^{(1)}$) are equal if for each element $x\in X^{(1)}$ there is an element $y\in Y^{(1)}$ such that $x=y$ and for each element $y\in Y^{(1)}$  there is an element $x\in X^{(1)}$  such that $y=x$, where \grqq$=$" means the following: 
 
 $1)$ the equality of individuals, when $x,y\in {\cal U}\stackrel{def}{=}{\cal U}^{(0)}$; 
 
 $2$) the equality of well-founded sets, when $x,y\in{\bf V}[{\cal U}]$;
 
 $3)$ the equality of skands, when $x,y\in {\cal U}^{(1)}$; 
 
 $4)$ the (iterative) equality  of sets in ${\bf V}[{\cal U}]^{(1)}$. 
 
 This axiom is correct and complete since ${\bf V}[{\cal U}]^{(1)}$ turns out to be {\it pseudo-well-founded}; i.e., for each element $x\in{\bf V}[{\cal U}]^{(1)}$  an arbitrary finite $\in$-chain starting with $x$ terminates after a finite number of steps in the following sense: $x\owns x_1\owns...\owns x_{n-1}\owns x_n$, where $x_n=\emptyset$, $n\geq 0$,  or $x_n\in{\cal U}$,  $n\geq 1$,  or $x_n\in {\cal U}^{(1)}$  and $x_{n-1}\notin {\cal U}^{(1)}$, $n\geq 0$. Here in each case $x_0=x$; in the former and the latter cases we put formally $x_{-1}={\bf V}[{\cal U}]^{(1)}$, i.e., $\emptyset\in{\bf V}[{\cal U}]^{(1)}$ and $x\in{\bf V}[{\cal U}]^{(1)}$, respectively, and an $\in$-sequences degenerate into one-element $\emptyset$ and $x$ or two-element  sequences ${\bf V}[{\cal U}]^{(1)}\ni x$ and ${\bf V}[{\cal U}]^{(1)}\ni \emptyset$, respectively. Thus, the iteration (a descending $\in$-chain) for each element $x$ is finite (modulo elements of ${\cal U}^{(1)}$) and we definitely know that $X^{(1)}=Y^{(1)}$ and vice versa.

 That is why we can call, for a while, elements of ${\cal U}^{(1)}$ \grqq pseudo-individuals" (more precisely, \grqq $1$-pseudo-individuals"; see below), because they along with the empty set $\emptyset$ and individuals ${\cal U}$ are constituents of sets in ${\bf V}[{\cal U}]^{(1)}$.
 At last, we put ${\bf U}^{(1)}={\bf V}[{\cal U}]^{(1)}\cup{\cal U}$. 
 
 This construction gives the following theorem.
 
 {\bf Theorem 1.} The set theory $NBG^{(1)}$ is consistent on the assumption that $NBG$ is consistent.

The {\bf Proof} is similar to the case when one extends the set theory $NBG$ of pure sets to the set theory $NBG[{\cal U}]$ with individuals ${\cal U}$ \cite{l56}. In our case we have an additional class of individuals ${\cal U}^{(1)}$ as \grqq pseudo-individuals" in the above sense. $\Box$
 
Clearly, ${\bf V}[{\cal U}]={\bf V}[{\cal U}]^{(0)}\subset{\bf V}[{\cal U}]^{(1)}\subset {\bf V}[{\cal U}]^-$ and Examples 3, 4, 5, show that the former inclusion is proper; we shall see that the latter inclusion is proper, too.

\bigskip

\parindent=0 cm
{\bf 4. Self-similar skands}

\parindent=0,5cm
The following skands are of great interest in the study of the relation $X\in X$ in the Set Theory which we are going to clarify in
this paper.

{\bf Definition 3}.
 A skand $X_{(\alpha_0,\alpha)}$  is called  {\it self-similar} if for each $\alpha_1$, $\alpha_0\leq\alpha_1<\alpha$, there is an equality $X_{(\alpha_0,\alpha)}=X_{(\alpha_1,\alpha)}$.

It happens iff $\alpha=\omega^\kappa$, where $\omega=\omega_0$ is the  initial countable ordinal, $\kappa\geq 1$, and  each $\alpha'$-component $X_{\alpha'}$ of $X_{(\alpha_0,\omega^\kappa)}$, $\alpha_0\leq\alpha'<\omega^\kappa$,  is the same set $\{x_1,x_2,...,x_\lambda,...,\{$. 

Indeed, it is well known that each remainder $\rho$ of an ordinal $\alpha\not=0$ is equal to $\alpha$ if and only if $\alpha=\omega^\kappa$, $\kappa\geq 0$, (\cite{l5}, Ch. VII, \S 7, Theorem 7). Recall that an ordinal $\rho$ is called a {\it remainder} of an ordinal $\alpha$ if $\rho\not=0$ and there exists an ordinal $\sigma$ such that $\alpha=\sigma+\rho$. 

Then our assertion when skands are self-similar  follows immediately from Definition 2. 

Note that in the case $\kappa=0$ the skands $X_{(0,1)}$ which are ordinary well-founded sets can be formally considered as \grqq self-similar" skands of length $1$.

Recall also that ordinal numbers $\alpha$ of the form $\omega^\kappa$, $\kappa\geq 0$, turn out to be the {\it prime components} or {\it principal numbers of addition}, i.e., ordinal numbers $\alpha$ such that there is no decomposition $\alpha=\beta +\gamma$ where $\beta<\alpha$ and $\gamma<\alpha$ (\cite{l66}, \S 19, Chap. XIV, Theorem 1, p. 323).

{\bf Example 6}. $X_{(0,\omega)}=\{a_0,\{a_0,\{a_0,\{...\}\}\}\}$, where $a_0$ is an individual or set in ${\bf U}$. 

This is a self-similar skand of length $\omega$ or the simplest example of an  extraordinary circular set of period 1. 

More generally,  $X_{(0,\omega)}=\{a_0,a_1,a_2,...,a_\lambda,...\{a_0,a_1,a_2,...,a_\lambda,...\{...\}\}\}$, where $a_\lambda$, $0\leq \lambda<\kappa$, are individuals or sets in ${\bf U}$, e.g., for $\lambda=i$, $a_\lambda=i$, $0\leq i<\omega=\kappa$, $X_{(0,\omega)}=\{0,1,2,3,...,\{0,1,2,3,...,\{...\}\}\}$.

{\bf Example 7}. $X_{(0,\omega^2)}=\{a_0,\{a_0,\{a_0,\{...\}\}\}\}$, where $a_0$, is an individual or a set in ${\bf U}$.

This is an example of a self-similar skand of length $\omega^2$. We picture it in the following figure:

$$
\begin{array}{ccccccccccccc}
\omega&\uparrow\\
&\bullet_{\{a_0,}&\bullet_{\{a_0,}&\bullet_{\{a_0,}&\bullet_{\{a_0,}&\bullet_{\{a_0,}&\bullet_{\{a_0,}&\bullet_{\{a_0,}&\bullet_{\{a_0,}&\bullet_{\{a_0,}&\bullet_{\{a_0,}\,\,...\\
&\omega n&\omega n+1&\omega n+2&\cdot&\cdot&\cdot&\cdot&\cdot&\cdot&\omega n+m\\
&\cdot&\cdot&\cdot&\cdot&\cdot&\cdot&\cdot&\cdot&\cdot&\cdot\\
&\bullet_{\{a_0,}&\bullet_{\{a_0,}&\bullet_{\{a_0,}&\bullet_{\{a_0,}&\bullet_{\{a_0,}&\bullet_{\{a_0,}&\bullet_{\{a_0,}&\bullet_{\{a_0,}&\bullet_{\{a_0,}&\bullet_{\{a_0,}\,\,...\\
&\cdot&\alpha'&\alpha'+1&\cdot&\cdot&\cdot&\cdot&\cdot&\cdot&\cdot\\
&\cdot&\cdot&\cdot&\cdot&\cdot&\cdot&\cdot&\cdot&\cdot&\cdot\\
&\bullet_{\{a_0,}&\bullet_{\{a_0,}&\bullet_{\{a_0,}&\bullet_{\{a_0,}&\bullet_{\{a_0,}&\bullet_{\{a_0,}&\bullet_{\{a_0,}&\bullet_{\{a_0,}&\bullet_{\{a_0,}&\bullet_{\{a_0,}\,\,...\\
&\omega 2&\omega 2+1&\omega 2+2&\cdot&\cdot&\cdot&\cdot&\cdot&\cdot&\omega 2+m\\
&\bullet_{\{a_0,}&\bullet_{\{a_0,}&\bullet_{\{a_0,}&\bullet_{\{a_0,}&\bullet_{\{a_0,}&\bullet_{\{a_0,}&\bullet_{\{a_0,}&\bullet_{\{a_0,}&\bullet_{\{a_0,}&\bullet_{\{a_0,}\,\,...\\
&\omega&\omega+1&\omega+2&\cdot&\cdot&\cdot&\cdot&\cdot&\cdot&\omega+m\\
&\bullet_{\{a_0,}&\bullet_{\{a_0,}&\bullet_{\{a_0,}&\bullet_{\{a_0,}&\bullet_{\{a_0,}&\bullet_{\{a_0,}&\bullet_{\{a_0,}&\bullet_{\{a_0,}&\bullet_{\{a_0,}&\bullet_{\{a_0,}\,\,...&\longrightarrow\\
&0&1&2&.&.&.&\alpha_0&\alpha_0+1&.&m&\omega
\end{array}
$$
\begin{center}
{\bf Fig. 1.}
\end{center}
and conclude that $X_{(\alpha_{0},\omega^2)}\approx X_{(\alpha',\omega^2)}$, for all $0\leq\alpha_0<\alpha'<\omega^2$; i.e., we see that $X_{(0,\omega^2)}$ is really a self-similar skand of length $\omega^2$. Notice also that, by Definition 2, $X_{(0,\omega)}=\{a_0,\{a_0,\{a_0,\{...\}\}\}\}\not=X_{(0,\omega^2)}=\{a_0,\{a_0,\{a_0,\{...\}\}\}\}$.

{\bf Definition 4.} A set $X$ is called {\it reflexive} if $X\in X$.

It is not clear at once that reflexive sets do exist, in spite of Mirimanoff's \grqq ensembles de deuxi\`{e}me sort"  \cite{l2} and Eklund's \grqq Mengen, die Elemente ihrer selbst sind" \cite{l25}.

\parindent=0,5cm Consider now an \grqq indeterminate" object $X$ in ${\bf V}[{\cal U}]^-$ of the following form:
\begin{equation}
\label{f1}
X=\{x_0,x_1,x_2,...,x_\lambda,..., X\}.
\end{equation}

It is clear that the equation $(\ref{f1})$ is a general form of reflexive sets if they exist.

{\bf Proposition 1.} Reflexive sets do exist. Moreover, there are a huge number  of different solutions (which form a proper class) of (\ref{f1}) in ${\bf V}[{\cal U}]^-$.  

{\bf Proof.} Indeed, the following self-similar skands:
\begin{equation}
\label{f2}
X_{(0,\omega^\kappa)}=\{x_0,x_1,x_2,...,x_\lambda,...,\{x_0,x_1,x_2,...,x_\lambda,...,\{...\}\}\}, 
\end{equation}
for each $\kappa\in{\bf On}$, $\kappa>0,$
are solutions of $(\ref{f1})$ because, by the  {\bf Axiom of Skand Existence}, objects in the form $(\ref{f2})$ do exist and by Definition 2, $X_{(0,\omega^\kappa)}=X_{(1,\omega^\kappa)}$, and we obtain 
\begin{equation}
\label{f22}
X=X_{(0,\omega^\kappa)}=\{x_0,x_1,x_2,...,x_\lambda,...,X_{(1,\omega^\kappa)}\}=\{x_0,x_1,x_2,...,x_\lambda,...,X\}. 
\end{equation}
$\Box$

{\bf Remark 3.} Proposition 1 shows that the relation $X\in X$ is extremely mulivalued and Russell was more or less right to call it \grqq meaningless" because without additional description it is undefined. One can say the same thing about the relations $X\in Y\in X$ which are also multivalued and, without additional description, are undefined, as Examples 3 and 4 tell us.

{\bf Remark 4.} Many years ago the author noticed (better to say perceived) \cite{l99} that there are self-similar skands of length geater than $\omega$ and for a long time has been thinking that only such skands of length $\omega^\kappa$, $\kappa\geq 2$, are solutions of $(\ref{f1})$ which differ from each other and from the solution of length $\omega$, i.e., Mirimanoff's extraordinary set solution. Now it is clear that not only, e.g.,
\begin{equation}
\label{f23}
X_{(0,\omega)}=\{x_0,x_1,x_2,...,x_\lambda,...,\{x_0,x_1,x_2,...,x_\lambda,...,\{...\}\}\} 
\end{equation}
and 
 \begin{equation}
\label{f24}
X_{(0,\omega^2)}=\{x_0,x_1,x_2,...,x_\lambda,...,\{x_0,x_1,x_2,...,x_\lambda,...,\{...\}\}\} 
\end{equation}
whose components are the same set $X_{\alpha'}=\{x_1,x_2,...,x_\lambda,...\}$, $0\leq\alpha'<\omega$ and $0\leq\alpha'<\omega^2$, respectively, are different solutions of (\ref{f1}),  but also, e.g.,
 \begin{equation}
\label{f25}
X_{(0,\omega+1)}=\{x_0,x_1,x_2,...,x_\lambda,...,\{x_0,x_1,x_2,...,x_\lambda,\{...\{1\}...\}\}\} 
\end{equation}
and 
 \begin{equation}
\label{f26}
Y_{(0,\omega+1)}=\{x_0,x_1,x_2,...,x_\lambda,...,\{x_0,x_1,x_2,...,x_\lambda,\{...\{2\}...\}\}\} 
\end{equation}
 are  also different solutions of (\ref{f1}) because $X_{(0,\omega+1)}=X_{(1,\omega+1)}$ and $Y_{(0,\omega+1)}=Y_{(1,\omega+1)}$ and $X_\omega=\{1\}\not=\{2\}=Y_\omega$, although these solutions {\it are not self-similar skands}. So, reflexive sets {\it need not be} self-similar skands. On the other hand, all these different solutions above are {\it isomorphic extraordinary sets} in the sense of Mirimanoff \cite{l1}, p. 40-41, and form a proper class. Moreover, they are {\it identically isomorphic} extraordinary sets in Mirimanoff's sense and thus must be equal. So, identically isomorphic extraordinary sets, or equal extraordinary sets {\it need not be equal} in $NBG^{(1)}$. In other words, reflexive extraordinary sets in \cite{l2} as well as in \cite{l25} are {\it not well-defined}, if we do not restrict them to skands of length $\omega$. Judging by his silence on the issue, it seems that Mirimanoff tacitly supposed that the length of the skands was equal to $\omega$. Consequently, in Mirimanoff's approach we see only $\omega$-phenomena and ignore trans-$\omega$-phenomena. This is indeed the origin of the following error in logic which we are now going  to clarify. $\Box$
 
$$
$$ 
 
 \bigskip

\parindent=0 cm
{\bf 5. Applications to Russell's paradox and its variants}

\parindent=0,5cm {\bf 5.1 Russell's paradox.} 

In 1903 Russell published the famous paradox  he had discovered two years previously and of which he had informed other mathematicians by correspondence. Here is the original quotation: \grqq We examined the contradiction resulting from the apparent fact that if $w$ be the class of all classes which as single terms are not members of themselves as many, then $w$ as one can be proved both to be and not to be a member of itself as many"(\cite{l47}, Chap. X, p. 107). Thus, in modern terminology, he defined the following set:

\begin{equation}
\label{f3}
R=\{X|\,\,X\notin X\},
\end{equation}
where $X$ are sets; i.e., $R$ is the set of {\it all sets} that are not members of themselves, or {\it the universal set} formed by the property $X\notin X$, which he and then all mathematicians considered to be a paradoxical  set, or Russell's antinomy. Therefore, the set $R$ was supposed to be inconsistent in Cantor's Na\"{i}ve Set Theory and the latter was called inconsistent, e.g., \cite{l131}, p. 488-489. 

The property or predicate $X\notin X$ was called Russell's condition.

Russell's argument is the following:
\begin{equation}
\label{f4}
R\in R\Longleftrightarrow R\notin R
\end{equation}
which is a contradiction, and that is why the set $R$ is inconsistent. Hence, $R$ does not exist from Poincar\'{e}'s point of view: \grqq En math\'{e}matiques le mot {\it exister} ne peut avoir qu'un sens: il signifie {\it exempt de contradiction}" (\cite{l6}, p. 162), as well as from Russell's: \grqq The contradiction proves that the class as one, if it ever exists, is certainly sometimes absent" (\cite{l47}, Chap. X, p. 107).

Later on, when sets and proper classes were being distinguished and any predicate (in particular, Russell's condition) formed a class which existed by one of the class-formation axioms, it was supposed that Russell's paradox said nothing other than that $R$ was a proper class, not a set.

\parindent=0,5cm  Now, by Proposition 1, one can immediately  obtain and see   that  Russell's well-known  paradox (\ref{f4})  {\it is not a paradox} at all,  because the assertion of implication 
\begin{equation}
\label{f5}
R\notin R\Longrightarrow R\in R
\end{equation}
 {\it is false}  since on the  {\it assumption} that $R$  {\it is a set}  
 
 $1)$ Russell's condition $X\notin X$ becomes  {\it impredicative}; 
 
 $2)$ by (\ref{f1}) and (\ref{f2}), the relation  $R\in R$ is  {\it multi-valued},
 
\parindent=0 cm   and therefore,  on such assumption the definition of $R$ becomes  {\it indeterminate} in $NBG^-$. It is so because $R$ {\it cannot be at once} a self-similar skand of different length $\omega^\kappa$, $\kappa\geq 1$, and this applies not only to self-similar skands but to many other reflexive sets which are solutions of $(\ref{f1})$ with indeterminate object $X=R$, where for each set $x\notin x$ there is a $\lambda\in {\bf On}$ such that $x=x_\lambda$. 

\parindent=0,5cm
Such  was the situation in the period of pre-axiomatic approaches to set theory. But it remained the same later after the creation of set-axiomatic systems and up to the present.
Although in $NBG$ a logical falsehood of Russell's paradox is a trivial fact of the calculation of the truth function,  because  if $R$  is a {\it set} (more precisely, $R=R'\cap{\bf V}[{\cal U}]={\bf WF}$, where $R'=\{X|\,\,X\notin X\}$ is the Russellian set in ${\bf V}[{\cal U}]^-$), then the relation $R\in R$ is  {\it always false}  and the relation $R\notin R$  is  {\it always true} since all sets in $NBG$ are well-founded and, consequently, by the truth-function-table values for implication, $R\in R\Longrightarrow R\notin R$ is  {\it true} and $R\notin R\Longrightarrow R\in R$ is  {\it false};  nevertheless, in $NBG$ as well as in $ZFA$ or other axiomatic set theories mathematicians use Russell's argument in proving that $R$ is a proper class, not a set, e.g., \cite{l7}, Chap. 4, \S 1, \cite{l5}, \cite{l88}, \cite{l14}, \cite{l98}, etc.

\parindent=0,5cm Thus, the set $R$ {\it is not  inconsistent because of Russell's famous paradox, or antinomy}, as  has been believed for  more than one hundred years \cite{l025} but {\it it is  inconsistent} because of a different argument (we call it the {\it Maximality Principle}).

The {\it paradox of Russell's paradox} is that the assertion of 
Proposition 2 below
is valid and true although it has usually  been proven incorrectly with the help of Russell's \grqq paradox", i.e., the \grqq classical inconsistency" of $R$. 
 
 \parindent=0,5cm Here is a correct proof of it without reference to Russell's \grqq argument".

In order to prove it, one needs only two axioms: 

({\bf I}) {\bf Axiom of the singleton}. If $X$ is a set, then there exists a set $\{X\}$ which contains $X$ and only $X$ as an element.

 ({\bf II}) {\bf Axiom of the union of two disjoint sets}. If $X$ and $Y$ are two sets which have no common elements, then there exists a set $X\cup Y$ whose elements are the elements of $X$ and of $Y$, taken together.

 These two axioms are assumed in all  set theories and even in  Cantor's Na\"{i}ve Set Theory. (We purposely formulate the axiom of union as a Cantorian version \cite{l1}, p. 481). Of course, we need something more as definitions of a set and of a class, and the distinctions between them. Nothing else. (Cantor's definition of an aggregate is quite enough:  \grqq By an \grq aggregate' we are to understand any collection into a whole $M$ of definite and separate objects $m$ of our intuition and our thought" (\cite{l132}, p. 85); Cantor called proper classes   the \grqq inconsistent systems"; as for the rest of the aggregates he called them the \grqq sets".)

{\bf Proposition 2.} {\it Russell's collection $R$ in $(\ref{f3})$ is a proper class, not a set.}

{\bf Proof}. Suppose the contrary assumption, i.e., $R$ is a set.  If $R\in R$ in any possible solution of the above equation $(\ref{f1})$, then $R$ is an element of $R$ and contains itself (in each case!); i.e., $R\in R\in R$. Then, by formula (\ref{f3}), $R\notin R$. 

(In Na\"{i}ve Set Theory this part of the proof may be the following:
if $R\in R$, \grqq unlikely, but not obviously impossible", then $R\in R\in R$, and hence, by formula (\ref{f3}), $R\notin R$.)

By ({\bf I}), $\{R\}$ is also a set, and by ({\bf II}), $R\cup\{R\}$ is a set. Since $R\notin R$ as we have already proved above, we obtain $R\subset R\cup\{R\}$; i.e., $R$ is a proper subset of $R\cup\{R\}$, which contradicts the maximality of $R$, since it is the set of {\it all} sets which do not contain themselves as elements. $\Box$

{\bf 5.2 The parametrization of Russell's set.}

Later in \cite{l48}, \cite{l8} and \cite{l98} a parametric version of a {\it real} Russellian set (not a proper class) was introduced. For any set $a\in {\bf V}[{\cal U}]$ (well-founded or non-well-founded, it does not matter) the following set 
\begin{equation}
\label{f11}
R_a=\{b\in a|\,\,b\notin b\}
\end{equation}
always exists as a set because it is an intersection of the proper class $R$ and the set $a$. As has been observed: \grqq There is nothing paradoxical about $R_a$. The reasoning that seemed to give rise to paradox only tells us that $R_a\notin a$". (See \cite{l98}, p. 60). That means that the assumption $R_a\in a$ gives the paradox $R_a\in R_a\Longleftrightarrow R_a\notin R_a$. This is the usual Russellian argument. Consider, e.g., the proof, given by Halmos in \cite{l48}, p. 6, changing only the notation, which we take from \cite{l8}. \grqq Can it be that $R_a\in a$? We proceed to prove that the answer is no. Indeed, if $R_a\in a$, $R_a\in R_a$ also (unlikely, but not obviously impossible), or else $R_a\notin R_a$. If $R_a\in R_a$, then, by $(\ref{f11})$, the assumption $R_a\in a$ yields $R_a\notin R_a$ $-$ a contradiction. If $R_a\notin R_a$, then, by $(\ref{f11})$ again, the assumption $R_a\in a$ yields $R_a\in R_a$ $-$ a contradiction again. This completes the proof that $R_a\in a$ is impossible, so that we must have $R_a\notin a$". We see that here the implication $R_a\notin R_a\Longrightarrow R_a\in R_a$ is false, too, as it was in Russell's argument.

Here is a correct proof of this statement.

{\bf Proposition 3.} $R_a\notin a.$

{\bf Proof.} It is clear, that $R_a\notin R_a$ because otherwise, (i.e., $R_a\in R_a$), it would be a member of $a$ and, therefore, by $(\ref{f11})$,  $R_a\notin R_a$. Thus, we have proved that $R_a\notin R_a$. Notice the implication that $R_a\in R_a\Longrightarrow R_a\notin R_a$ is always true in similar situations. 

Suppose now that $R_a\in a$. Then $R_a\cup\{R_a\}$ is a subset of $a$; moreover, it is a subset of $a$ whose elements do not contain themselves, in particular, $R_a$, as we have proved above. Therefore, $R_a$ is a proper subset of $R_a\cup\{R_a\}$; i.e., $R_a\subset R_a\cup\{R_a\}\subset a$ which, by  $(\ref{f11})$, is in contradiction with the maximality of $R_a$, since it is the set of {\it all} elements $b$ of $a$ such that $b\notin b$. Consequently, $R_a\notin a$. $\Box$

{\bf Remark 5}. In spite of the gap in Halmos's proof, the author of this paper cannot refrain from relating a witty observation made by Paul R. Halmos (the author once met him and appreciated his wit as well as his brilliant lectures and books): \grqq The most interesting part of this conclusion is that there exists something (namely $R_a$) that does not belong to $a$. The set $a$ in this argument was quite arbitrary. We have proved, in other words, that {\it nothing contains everything}, or, more spectacularly, {\it there is no universe}.  \grq Universe' here is used in the sense of \grq universe of discourse,' meaning, in any particular discussion, a set that contains all the objects that enter into that discussion. In older (pre-axiomatic) approaches to set theory, the existence of a universe was taken for granted, and the argument in the preceding paragraph was known as {\it the Russell paradox}. The moral is that it is impossible, especially in mathematics, to get something for nothing. To specify a set, it is not enough to pronounce some magic words (which may form a sentence such as \grq $x\notin x$'); it is necessary also to have at hand a set to whose elements the magic words apply". (Ibid., p. 6-7).

{\bf Remark 6}. Here we notice, contrary to the main idea of this paper, that {\it there are  axiomatic set theories} where, on the assumption that Russell's class $R$ is a set, {\it Russell's paradox really arises}; moreover, the classical proof that $R$ is not a set as well as $R_a\notin a$ because of Russell's paradox is {\it correct}. For example, such is Aczel's theory \cite{l93}, based on an extremely natural extension of the Zermelo conception of a set-axiomatic system. He rejected the foundation axiom in ${\bf ZFC}$ and proposed his own anti-foundation axiom ${\bf AFA}$. Then, after improving the  axiom of extensionality, changed it to the axiom of strong extensionality (bisimulation, in current terminology), and obtained an axiomatic theory ${\bf ZFC}^-+{\bf AFA}$ of so-called hypersets. In this theory the implications $R\notin R\Longrightarrow R\in R$ and $R_a\notin R_a\Longrightarrow R_a\in R_a$ are true, contrary to the general situation considered above (i.e., in systems without ${\bf AFA}$) because ${\bf AFA}$ implies a unique solution of an equation $(\ref{f1})$. But working in ${\bf ZFC}^-+{\bf AFA}$, the authors of \cite{l8} and \cite{l98},  very surpisingly, do not refer to ${\bf AFA}$ and employ the classical Russellian argument, which, in the absence of ${\bf AFA}$, is wrong as we saw above. Notice also, that the {\it method of maximality} offered above in this paper works well in ${\bf ZFC}^-+{\bf AFA}$ without reference to Russell's paradox or to ${\bf AFA}$.

{\bf 5.3  Zermelo's paradox.} 

Zermelo founded his own paradox independently of Russell and said in 1908 that he had mentioned it to Hilbert and other people already before 1903. Later, in 1936, in his letter to Scholz, he wrote that the set-theoretical paradoxes were often discussed in the Hilbert circle around 1900, and he himself had at that time given a precise formulation of the paradox which was later named after Russell (see \cite{l121}). We shall see below that Zermelo's opinion was mistaken, and that Russell's and his  paradoxes are similar but {\it not the same}. Zermelo's paradox is the following: a set $M$ that comprises as elements all of its subsets is inconsistent. Indeed, consider the set $M_0$ of all elements of $M$ which are not elements of themselves (e.g., the empty set is in $M_0$) This set is a subset of $M$ and hence by assumption on $M$, $M_0\in M$. If $M_0\in M_0$, then $M_0$ is not a member of itself. Hence $M_0\notin M_0$ and since $M_0\in M$, $M_0\in M_0$: contradiction. (See \cite{l12}, \S$\,$ 2.4; \cite{l025}, p. 507). We notice the same mistake: the implication $M_0\notin M_0$ and $M_0\in M \Longrightarrow M_0\in M_0$ is false. Here is the solution of Zermelo's paradox.

{\bf Proposition 4.} There is no set $M$ such that it comprises as elements all of its subsets.

{\bf Proof.} Suppose the contrary, and such set $M$ exists. Then $M_0=\{X\in M|\,\,X\notin X\}$ is well-defined and $M_0\in M$. If $M_0\in M_0$ in any possible way, then $M_0\in M_0\in M_0$ and, clearly, by definition of it, $M_0\notin M_0$. (The first part of Zermelo's proof is correct.) By axiom ({\bf I}), the singleton $\{M_0\}$ also exists  and it is a subset of $M$, i.e., $\{M_0\}\subset M$ (because $M_0\subset M$ and, by definition of $M$, $M_0\in M$). By axiom ({\bf II}), we obtain that $M_0\cup\{M_0\}$ is also a subset of $M$ and consists of elements which are not elements of themselves because in addition to all elements of $M_0$ with such a condition, $M_0$ itself satisfied the same condition, i.e., $M_0\notin M_0$, as was proved in the first part. Clearly,  $M_0$ is a proper subset of $M_0\cup\{M_0\}$; i.e., $M_0\subset M_0\cup\{M_0\}$ which is in contradiction with the maximality of $M_0$. Consequently, $M$ doest not exist. 

Notice also that contrary to Russell's \grqq set" $R$, which is not a set but  does exist as a proper class, Zermelo's set $M$ does not exist as a proper class at all; in other words, there is no mathematical object such as Zermelo's set $M$, but there is an object $R$ in $NBG^-$ and in $NBG$. Moreover, in the latter it is called the Universe ${\bf V}[{\cal U}]$ or the class of all well-founded sets ${\bf WF}$. That is why Zermelo's and Russell's paradoxes are meaningfully different. $M$ cannot be a proper class because in its rigid definition $M$ has to be an element of $M$, which is impossible for proper classes. $\Box$

\parindent=0,5cm

{\bf 5.4 The  Paradox of Propositions} (\cite{l025}, p. 260; \cite{l12}, 2.3). 

Let $m$ be a set of propositions and let $\Pi m$ be the proposition that \grqq every proposition of $m$ is true" (regarded as a possibly infinitary conjunction);  then, if $m$ and $n$ are different sets, the propositions $\Pi m$ and $\Pi n$ are different; i.e., the map associating to $m$ its product $\Pi m$ is injective. Therefore, if we consider the set
\begin{equation}
\label{f23}
R=\{p|\,\,\exists m(\Pi m=p\,\,\,\&\,\,\, p\notin m)\}
\end{equation}
we have, by injectivity, a contradiction. Notice that we do not identify here equivalent propositions; if we did, no contradiction could be arrived at. This is an example of a Russell-type paradox; its classical exposition is the following: $\Pi R\in R\Longleftrightarrow \Pi R\notin R$. As above  $\Pi R\in R\Longrightarrow \Pi R\notin R$ is true because of the injectivity of the map considered above,  and $\Pi R\notin R\Longrightarrow \Pi R\in R$ is false, because the latter relation $\Pi R\in R$ is multivalued. Hence, $R$ is not maximal, since $R$ is a proper subclass of $R\cup\{\Pi R\}$; i.e., $R\subset R\cup\{\Pi R\}$.
The contradiction with the maximality of $R$ tells us that $R$ is a proper class, not a set, and moreover, the proposition $\Pi R$ does not exist, since the infinitary conjunction is possible only for sets of conjunctions. $\Box$

{\bf Remark 7.} The state of affairs is that the \grqq set" $R$ {\it is inconsistent}  because of the {\it maximality} or {\it universality}  of it,  since the supposed set $R$ is a family of {\it all} and at the same time {\it not all} objects with a certain condition (here $X\notin X$), but not because of Russell's paradox (\ref{f4}), which is false at least in axiomatic systems without ${\bf AFA}$ or perhaps other similar axioms. In particular, its logical falsity in an $NBG$-type set theory with ${\bf FA}$ is an immediate consequence of the truth-function table values for implication.  $R$ is consistent and exists as a proper class, contrary to Mirimanoff \cite{l3}, p. 31, and others who accept Poincar\'{e}'s point of view, mentioned above; and Russell's condition excludes class quantifiers.  Hence {\it there is no contradiction} in such a universal family or collection of objects as all sets which are not members of themselves. The singleton  $\{R\}$ is not empty, contrary to Mendelson (\cite{l7}, Chap. 4, \S 1, Axiom N),  and it is also a mathematical object, but not a fact in the world ${\bf V}^-$, clearly, as well as in ${\bf V}$  (in the terminology of Barwise and Etchemendy \cite{l8});  it is a member of a hyper-class (Quine and Morse \cite{l9}). As to the set $R_a$ it is consistent and always exists, as opposite to the Zermelo's set $M$, which is inconsistent and hence does not exist as a set   or as a proper class, which was already stated above. $\Box$

\bigskip

\parindent=0 cm
{\bf 6. Other set-theoretical paradoxes}

\parindent=0,5cm The same argument is valid for the \grqq set" ${\bf On}$ (Burali-Forti's \grqq paradox" \cite{l10}), for the \grqq set" ${\bf Card}$ of all cardinal numbers (Cantor's \grqq paradox," first mentioned in the letter to Dedekind of August 31. 1899, and later in \cite{l11}),  for the \grqq set" ${\cal U}$ (Hilbert's \grqq paradox" \cite{l025}, p. 505), and for the \grqq set" ${\bf WF}$ (Mirimanoff's \grqq paradox" \cite{l2}). 

We shall also prove that the family ${\bf NWF}$ of all non-well-founded sets is not a set but a proper class and it does not concern Russell-like paradoxes at all. 
\parindent=0,5cm
All proofs will be based on  the method of the maximality (universality) of a collection of objects satisfying a given condition  presented above (our main pattern or patent: the Maximality Principle). Admittedly, it was implicitly used by Cantor in proving that the second number-class set is not countable \cite{l111}, \S 16, [Theorem] D. His proof takes up at least one page--more precisely, 32 lines). Now this proof requires only five  lines. 

{\bf Proposition 5.} For each initial ordinal $\omega_\alpha$, the set of all ordinals $\lambda$ such that $\lambda<\omega_\alpha$ cannot be of a power smaller than the power of $\omega_\alpha$.

{\bf Proof.} Indeed, if the power of the set $X=\{\lambda|\,\,\lambda<\omega_\alpha\}$ were smaller than the power of $\omega_\alpha$, then the power of the next ordinal number $\beta=X\cup\{X\}$ would be smaller than the power of $\omega_\alpha$, too. Since $\beta$ is greater than each ordinal $\lambda$ in $X$, it is not a member of $X$. Hence $X\subset X\cup\{\beta\}$; i.e., $X$ is a proper subset of $X\cup\{\beta\}$, which is in contradiction with the maximality of $X$. $\Box$

{\bf 6.1 Burali-Forti's  paradox, or the antinomy of the greatest ordinal.}

The earliest antinomy  in set theory was published in 1897 in \cite{l10}. It has, however, a non-trivial history, starting out from the fact that there was nothing paradoxical in it, since in \cite{l10} there was no contradiction. This is because Burali-Forti had misconstrued Cantor's definition of a well-odered set, and used his own notion of a different kind of ordered sets  which he called \grqq perfectly ordered classes", and proved that such classes are non-well-ordered. Russell \cite{l47}, p. 323, reformulated the argument of Burali-Forti as a contradiction and gave it its present name. (This observation was borrowed from \cite{l913}, p. 306-307, and  \cite{l025}, p. 350). And we will finish this story with an amusing remark given by Halmos: \grqq The contradiction, based on the assumption that there is a set of all ordinals, is called the Burali-Forti paradox. (Burali-Forti was one man, not two.)" (See \cite{l48}, p. 80.) We are going to show that Russell and others did not supply Burali-Forti's  paradox with a contradiction, although there are axiomatic systems where this contradiction arises.

Look at the modern explanation of Russell's correction of Burali-Forti's antinomy of the greatest ordinal, e.g., in \cite{l025}, p. 350. \grqq Let $\Omega$ be the ordinal number of the well-ordered set of all ordinals, ${\bf No}$. But, for every ordinal $\alpha$, $\alpha+1>\alpha$. Thus $\Omega+1>\Omega$. But, for every $\alpha\in{\bf No}$, $\alpha\leq \Omega$. Thus $\Omega+1\leq\Omega$".

Recall that, due to Cantor, \grqq ordinal numbers" are ordinal types of \grqq well-ordered aggregates", i.e., well-ordered sets. Proper classes of equivalent well-ordered sets can be eliminated due to von Numann's method \cite{l090}, by  choosing a canonical set-representative of each proper class, that is to define an ordinal as a set $\alpha$ of sets $X$ which is well-ordered by the relation $\in$ between its elements and transivity, i.e., if $Z\in Y\in X$, then $Z\in X$, starting, e.g., with the empty set, i.e., $\emptyset;$ $\{\emptyset,\{\emptyset\}\};$ $\{\emptyset,\{\emptyset\},\{\emptyset,\{\emptyset\}\}\};$\,\,...\,\,. (It would perhaps have been fairer to say that the idea of this method for the first time was given by Mirimanoff \cite{l2}, p. 46.)   The collection of {\it all ordinals} ${\bf On}$ is well ordered by the relation $\in$  (see \cite{l13}, p. 8) and class transitive; i.e., if $Z\in Y\in {\bf On}$, then $Z\in {\bf On}$. Note that the class of all ordinals ${\bf No}$ above is in bijection  with ${\bf On}$; moreover, they are order isomorphic. Note also that in $NBG$ all ordinals are well-founded sets. Thus, for each $ X\in\alpha$ we have $X\not\in X$, in particular, $\alpha\notin\alpha$, and moreover, we can even omit \grqq well" in \grqq well-ordered" (see, e.g., \cite{l7}, Chap. 4, \S 5). The situation in $NBG^-$ and in Na\"{i}ve Set Theory is more complicated: \grqq well-ordered" is sufficient and $\alpha\in\alpha$ can happen. Nevertheless, we remember that in general the relations $\alpha\in\alpha$ as well as $\alpha\in\beta\in\alpha$ are indeterminate.

{\bf Proposition 6.} {\it  ${\bf On}$ is a proper class, not a set.}

{\bf Proof}.
Suppose that ${\bf On}$  is a set. Thus, it is a transitive set well-ordered by the relation $\in$ and hence is an ordinal.  Then, by ({\bf I}), $\{{\bf On}\}$ is also a set, and, by ({\bf II}), ${\bf On}\cup\{{\bf On}\}$ is a set which is evidently well-ordered by the relation $\in$ and transitive, too. Consequently, ${\bf On}\cup\{{\bf On}\}$ is an ordinal. (Notice that in $NBG$ ${\bf On}\notin{\bf On}$ and $\{{\bf On}\cup\{{\bf On}\}\}\notin{\bf On}$ because all sets in $NBG$ are well-founded; as to $NBG^-$ and Na\"{i}ve Set Theory, ${\bf On}\in{\bf On}$ may occur, e.g., as self-similar skands of different lengths, but the relations ${\bf On}\in{\bf On}\cup\{{\bf On}\}\in{\bf On}$  are indeterminate.)  Moreover, by the same argument, $\{{\bf On}\cup\{{\bf On}\}\}$ and ${\bf On}\cup\{{\bf On}\cup\{{\bf On}\}\}$ are sets and thus ${\bf On}$ is a proper subset of ${\bf On}\cup\{{\bf On}\cup\{{\bf On}\}\}$ because of the relation noted above; i.e., ${\bf On}\cup\{{\bf On}\}\notin{\bf On}$. The relation ${\bf On}\subset{\bf On}\cup\{{\bf On}\cup\{{\bf On}\}\}$ is in contradiction with the maximality of ${\bf On}$. Thus ${\bf On}$ is not a set. It is a proper class since its existence is guaranteed by the predicate: {\it to be  an ordinal}. $\Box$

{\bf Remark 8.} Burali-Forti's antinomy of the greatest ordinal is seen now as two inconsistent inequalities $\Omega+1>\Omega$ and $\Omega+1\leq\Omega$. Consequently, we have the relations $\Omega<\Omega+1\leq\Omega$ in ${\bf No}$ and hence, by the isomorphism between ${\bf No}$ and ${\bf On}$ mentioned above, we obtain the relations $\Omega\in\Omega\cup\{\Omega\}\in\Omega$ in ${\bf On}$, which are indeterminate, as we have already seen many times. So in the classical argument, to conclude by definition that $\Omega+1\in \Omega$ is false. (The true state of affairs is that after the supposition that ${\bf On}=\Omega$ is a set, the definition of ${\bf On}$ has been radically changed; i.e., each possible relation ${\bf On}\in {\bf On}$ or $\Omega\in\Omega\cup\{\Omega\}\in\Omega$ changes the set ${\bf On}$ itself, and that is why ${\bf On}$ is indeterminate.) Consequently, when Russell supposedly {\it showed a contradiction} in the Burali-Forti example, {\it he was wrong}. If we suppose that ${\bf On}$ is a set, then naturally ${\bf On}\cup\{{\bf On}\}=\Omega+1$ is an ordinal and $\Omega+1>\Omega={\bf On}$. But $\Omega\cup\{\Omega\}\notin\Omega$, i.e., $\Omega+1\not\leq\Omega$. Because $\Omega={\bf On}$ is a well-founded set in $NBG$, and in $NBG^-$ as well as in Na\"{i}ve Set Theory, the relations $\Omega\cup\{\Omega\}\in\Omega$ are indeterminate. In our proof of the inconsistency of the set ${\bf On}$ we omit this erroneous step ${\bf On}\cup\{{\bf On}\}\in{\bf On}$ but add a new ordinal ${\bf On}\cup\{{\bf On}\}$ (different from each element of ${\bf On}$) to ${\bf On}$ and obtain a larger set than ${\bf On}$, and this is a real contradiction, because ${\bf On}$ is a collection of {\it all} ordinals. (We are aware that the relations $\Omega\in\Omega\cup\{\Omega\}\in\Omega $ here are much more complicated than in Examples 3 and 4; i.e., $X\in Y\in X$. Therefore, it is much more indeterminate. So it is a good idea to omit such relations, which can be done in the maximality method.) In ${\bf ZFC}+{\bf AFA}$, e.g., Russell's argument and the classical proof are correct, but only with reference to ${\bf AFA}$. In ${\bf ZFC}+{\bf AFA}$ our proof works, too, but there is no need for such a reference.

{\bf 6.2 Cantor's paradox.}

The historically second set-theory antinomy belongs to Cantor: the set ${\bf V}^*$ of all sets  (which was called later \grqq Cantor's paradise") is inconsistent or paradoxical because the power set ${\bf P}{\bf V}^*$   is always  of  a greater cardinality than   ${\bf V}^*$, but at the same time ${\bf P}{\bf V}^*\subseteq{\bf V}^*$ since ${\bf V}^*$ is the most inclusive set of sets. Notice that the same mistake obtains: the relation ${\bf P}{\bf V}^*\subseteq{\bf V}^*$ implies, in particular, ${\bf V}^*\in{\bf V}^*$ since ${\bf V}^*\in {\bf PV}^*$. But the relation ${\bf V}^*\in{\bf V}^*$ is multivalued  and indeterminate, as we saw above, and therefore, false. Thus there is no real paradox or contradiction in Cantor's argument; moreover, in this case Cantor's paradox is not a paradox at all. But what was correctly asserted via Cantor's \grqq paradox" was the following:

{\bf Proposition 7.} {\it  ${\bf V}^*$ is a proper class, not a set, or in Cantor's denomination ${\bf V}^*$ is an \grqq inconsistent system".}

 {\bf Proof}. Suppose that ${\bf V}^*$  is a set. Then, by the {\bf Power-set Axiom}, the collection (system) ${\bf P}{\bf V}^*$ of all subsets of ${\bf V}^*$ is a set. It is clear, by ({\bf I}), that there is an injection $i:{\bf V}^*\rightarrow {\bf P}{\bf V}^*$, given by $i(X)=\{X\}$. On the other hand, ${\bf P}{\bf V}^*$ is not a subset of ${\bf V}^*$. Otherwise, there would be an embedding $j:{\bf P}{\bf V}^*\rightarrow{\bf V}^*$ and, by the Cantor-Bernstein theorem, ${\bf P}{\bf V}^*$ and ${\bf V}^*$ would be equivalent, which is false because of Cantor's theorem $|{\bf P}{\bf V}^*|>|{\bf V}^*|$. Thus, there exists an element $X\in {\bf P}{\bf V}^*$ such that $X$ is a set and $X\notin{\bf V}^*$.  Then, by ({\bf I}), $\{X\}$ is a set and, by ({\bf II}), ${\bf V}^*\cup\{X\}$ is a set. Moreover, ${\bf V}^*\subset{\bf V}^*\cup\{X\}$, i.e., ${\bf V}^*$ is a proper subset of ${\bf V}^*\cup\{X\}$, which is in contradiction with the maximality of ${\bf V}^*$. Thus ${\bf V}^*$ is not a set but a proper class, given by the predicate $X=X$, i.e., ${\bf V}^*\stackrel{def}{=}\{X|\,\,X=X\,\,\&\,\,\exists \{X\}\}$ and called the universal class (see \cite{l14}, p. 124). (Of course, we could take $X={\bf V}^*$ at the beginning, and since ${\bf V}^*\in{\bf V}^*$ is indeterminate, conclude that ${\bf V}^*\notin{\bf V}^*$ and, by supposition that ${\bf V}^*$  is a set, obtain, by $(\bf I)$ and $(\bf II)$, that ${\bf V}^*$ is a proper subset of ${\bf V}^*\cup\{{\bf V}^*\}$. In other words, we could repeat our method of maximality; but we wished to find a mistake in Cantor's argument.) $\Box$
 
 Cantor himself reached the faulty conclusion that ${\bf P}{\bf V}^*\subseteq {\bf V}^*$, which is in contradiction with their cardinality, and obtained a supposed \grqq paradox". But it was only a presumption, which was proved to be false by Cantor's theorem $|{\bf P}{\bf V}^*|>|{\bf V}^*|$, nothing more. Therefore, there is no Cantor's paradox. It would be a paradox if Cantor could prove that ${\bf P}{\bf V}^*\subseteq {\bf V}^*$. Meanwhile his conclusion ${\bf P}{\bf V}^*\subseteq {\bf V}^*$ was made via the Definition of ${\bf V}^*$; but after the assumption that ${\bf V}^*$ was a set, the Definition of ${\bf V}^*$ was radically changed, and ${\bf V}^*\in {\bf V}^*$ became automatically multivalued. 
 
 Absolutely by a similar  argument Cantor proved that the collection ${\bf Card}$ of all cardinal numbers is not a set (August 31, 1899, letter to Dedekind) and Cantor made the following conclusion: \grqq... the system ${\bf Card}$ { is not a set. That is why there exist certain} pluralities which are not at the same time wholes (unities), i.e., pluralities for which the real \grq mutual being of all their elements' is {\it impossible}. They are those I call the \grq inconsistent systems'; as for the rest I call them the \grq sets'$\,\,$" \cite{l125}.

 Here is a correct proof of Cantor's statement.
 
{\bf Proposition 8.} {\it  ${\bf Card}$ is a proper class, not a set, or in Cantor's denomination ${\bf Card}$ is an \grqq inconsistent system".}

{\bf Proof}. Suppose that ${\bf Card}$  is a set. It is known that ${\bf Card}$ is in a bijection with the family ${\bf In}$ of all finite numbers together with the all initial ordinals $\omega_\kappa$ of ${\bf On}$ which we identify with the initial segments $(0,\omega_\kappa)$ of ${\bf On}$, $\kappa\in{\bf On}$ (if $\kappa=0$, then the corresponding segment is the empty set).  Since we supposed that ${\bf Card}$  is a set, therefore ${\bf In}$ is a set, by the {\bf Replacement Axiom} for sets, and we obtain that the discrete sum $X=\coprod\limits_{\omega_\kappa\in {\bf In}} [0,\omega_\kappa)$ is also a set, by the {\bf Union Axiom} for sets, and, by the {\bf Power-set Axiom}, ${\bf P}X$ is a set. Hence there exists an initial ordinal $\omega_\lambda$ such that ${\bf P}X$ and $[0,\omega_\lambda)$ are bijective. It is clear that $\omega_\lambda>\omega_\kappa$, for each $\omega_\kappa\in{\bf In}$ (because of Cantor's theorem: $|{\bf P}[0,\omega_\kappa)|>|[0,\omega_\kappa)|$, $\kappa\in{\bf On}$), which is in contradiction with the maximality of ${\bf In}$. Thus ${\bf Card}$  is not a set but a proper class, since ${\bf Card}\subset{\bf V}^*$.  $\Box$

 {\bf Remark 9.} Notice that from the fact that ${\bf On}$ is not  a set but a proper class, as has been proved above, it does not follow that its subclass ${\bf Card}$ is also a proper class. On the other hand, there is a bijection between ${\bf On}$ and ${\bf Card}$ because there is an injection $j:{\bf On}\rightarrow{\bf Card}$, given by $j(\kappa)=\omega_\kappa$,  $\kappa\in{\bf On}$, and we apply the Cantor-Bernstein theorem. Thus ${\bf Card}$ is not a set, but a proper class.

{\bf 6.3 Hilbert's paradox.}
 
 At least earlier than 1905 Hilbert formulated a paradox of his own which he had considered \grqq purely mathematical" in the sense that it did not make use of notions from Cantor's theory of cardinals and ordinals. Hilbert never published the paradox: \grqq I have never published this contradiction, but it is known to set theorists, especially to G. Cantor" \cite{l491}, p. 204.
 
 Let us quote it from \cite{l025}, p. 505-506.  \grqq The paradox is based on a special notion of set which Hilbert introduces by means of two set formation principles starting from the natural numbers. The first principle is the {\it addition principle} ({\it Additionsprinzip}). In analogy to the finite case, Hilbert argued that the principle can be used for uniting two sets together \grq into a new conceptual unit', a new set that contains the elements of both sets. This operation can be extended: \grq In the same way, we are able to unite several sets and even an infinitely many into a union.' The second principle is called  the {\it mapping principle} ({\it Belegungsprinzip}). Given a set ${\cal M}$, he introduces the set ${\cal M}^{{\cal M}}$ of {\it self-mappings} ({\it Selbstbelegungen}) of ${\cal M}$ to itself. A self-mapping is just a total function (\grq transformation') which maps the elements of ${\cal M}$ to elements of ${\cal M}$.
 
 Now, he considers all sets which result from the natural numbers \grq by applying the operations of addition and self-mapping an arbitrary number times.' By the addition priciple, which allows us to form the union of arbitrary sets, one can \grq unite them all into a sum set ${\cal U}$, which is well-defined.' In the next step the mapping principle is applied to ${\cal U}$, and we get ${\cal F}={\cal U}^{{\cal U}}$ as the set of all self-mappings of ${\cal U}$. Since ${\cal F}$ was built from the natural numbers, by using the two principles alone, Hilbert concludes that it has to be contained in ${\cal U}$. From this fact he derives a contradiction.
 
 Since there are \grq not more' elements in ${\cal F}$ than in ${\cal U}$ there is an assignment of the elements $u_i$ of ${\cal U}$ to elements $f_i$ of ${\cal F}$ such that all elements of $f_i$ are used. Now one can define a self-mapping $g$ of ${\cal U}$ which differs from all $f_i$. Thus, $g$ is not contained in ${\cal F}$. Since ${\cal F}$ contains all self-mappings, we have a contradiction. In order to define $g$, Hilbert used Cantor's diagonalization method. If $f_i$ is a mapping $u_i$ to $f_i(u_i)=u_{f_i^{i}}$ he chooses an element $u_{g^{i}}$ different from $u_{f_i^{i}}$ as the image of $u_i$ under $g$. Thus, we have $g(u_i)=u_{g^{i}}\not=u_{f_i^{i}}$ and $g$ \grq is distinct from any mapping $f_k$ of ${\cal F}$ in at least one asssignment.' Hilbert finishes his argument with the following obsorvation:
 
 \grqq We could also formulate this contradiction so that, according to the last consideration, the set ${\cal U}^{{\cal U}}$ is always bigger [of greater cardinality] than ${\cal F}={\cal U}$, but according to the former, is an element of ${\cal F}={\cal U}$." $\Box$
 
 Let us comment on the above proof. In spite of the difference between Cantor's and Hilbert's definitions of set (e.g., the empty set $\empty set$ is not in ${\cal U}$, and other non-intersections) their proofs and arguments are similar and contain the same mistake. Let us see it in Hilbert's paradox. Hilbert has proved, using Cantor's diagonalization method, only that his supposition that ${\cal F}\subseteq{\cal U}$ is false, nothing more. Thus, there is no contradiction in such an unfortunate supposition. If he had proved that ${\cal F}\subseteq{\cal U}$ was true, there would be a real paradox to be named after him. What he has essentially proved is that there is an element $g$ of ${\cal F}$ such that it was not an element of ${\cal U}$. Then by his first principle ${\cal U}\cup{\cal F}$ is a set; moreover, ${\cal U}$ is a proper subset of ${\cal U}\cup{\cal F}$ because $b$ is not in ${\cal U}$ and it is in ${\cal U}\cup{\cal F}$. Thus, by the {\it universality} of ${\cal U}$, (i.e., ${\cal U}$ is the set of {\it all} Hilbert's sets), ${\cal U}$ is not a Hilbert set, nor is ${\cal F}={\cal U}^{{\cal U}}$ a Hilbert set. There is no paradox here at all. There is only a proof by contradiction, and the supposition that ${\cal U}$ is a set is simply false. $\Box$

 {\bf 6.4 Mirimanoff's paradox.} 
 
 In \cite{l2}, p. 43, Mirimanoff was actually concerned with the collection ${\bf WF}$ of all well-founded sets, and presented his paradox:

  \begin{equation}
  \label{f7}
    {\bf WF}\in {\bf WF}\Longleftrightarrow {\bf WF}\notin {\bf WF} 
  \end{equation}
as an analogue of Russell's paradox. And he concluded that the set ${\bf WF}$ did not exist, i.e., in Cantor's terminology, it is an inconsistent system, or in modern terminology, it is not a set. His argument is the following. It is clear that  ${\bf WF}\notin{\bf WF}$ because if it were not so and ${\bf WF}\in{\bf WF}$, then there would exist an infinite descending chain in ${\bf WF}$, e.g.,  ${\bf WF}\owns{\bf WF}\owns{\bf WF}\owns...$ and we would obtain that ${\bf WF}$ is non-well-founded and hence ${\bf WF}\notin{\bf WF}$. 

On the other hand, if ${\bf WF}$ is non-well-founded; there exists an infinite chain $x_0\owns x_1\owns x_2...$ in ${\bf WF}$ for some $x_0\in{\bf WF}$, and then $x_0$ is evidently non-well-founded, which is wrong. Therefore, ${\bf WF}$ is well-founded. Hence ${\bf WF}\in{\bf WF}$; contradiction. 

The last conclusion is false because ${\bf WF}\in{\bf WF}$ is indeterminate as we saw above more than once. 

Here is a correct proof of the following 

{\bf Proposition 9.} {\it ${\bf WF}$ is  a proper class, not a set.} 

{\bf Proof.}
Assume that ${\bf WF}$ is a set. Then by $({\bf I})$ the singleton $\{{\bf WF}\}$ is also a set. Since all elements of ${\bf WF}$ are well-founded we conclude that  ${\bf WF}$ is not an element of ${\bf WF}$ because if it were, then there would be an evident infinite descending $\in$-sequence ${\bf WF}\ni{\bf WF}\ni {\bf WF}\ni...$ Thus ${\bf WF}$ and $\{{\bf WF}\}$ are two sets without common elements. Then, by $({\bf II})$,   $X'={\bf WF}\cup \{{\bf WF}\}$ is also a set. Moreover, its elements are well-founded and ${\bf WF}$ is a proper subset of $X'$ because ${\bf WF}$ is a member of $X'$ but not a member of ${\bf WF}$. This is in contradiction with the maximality of ${\bf WF}$ because it is the well-founded universe, since ${\bf WF}$ is the collection of {\it all} well-founded sets. Consequently, the assertion that ${\bf WF}$ is a set is false. By the {\it well-founded predicate} and the Existence Axiom for a class,  ${\bf WF}$ is a proper class. $\Box$

 In $NBG^{(1)}$  we have the following  collection:
  \begin{equation}
  \label{f6}
    {\bf T}=\{X|\,\,X\in X\} 
  \end{equation}
which is symmetric to Russell's \grqq set" $R=\{X|\,\,X\notin X\}$ \cite{l931}, p. 277, and on the assumption that it is a set it does not look paradoxical or inconsistent.   We are not saying that ${\bf T}\in {\bf T}$, because otherwise, the relation ${\bf T}\in {\bf T}$ cannot even be called {\it indeterminate} as above; it is {\it meaningless} as Russell tells us in \cite{l55}, p. 80. Indeed, ${\bf T}$ contains, e.g., all skands of arbitrary lengths; ${\bf T}\in{\bf T}$, if it were the case, would be a skand of a fixed length, and other skands of different lengths would be different from this fixed one.   In other words, ${\bf T}$ cannot be a skand of {\it all lengths at once}.  Nevertheless, 

{\bf Proposition 10.} {\it ${\bf T}$ is a proper class, not a set.}

{\bf Proof.} Consider for an arbitrary well-founded set $X$, i.e., $X\in{\bf WF}$, a self-similar skand $X_{(0,\omega)}$ with the components $X_i=X$, $0\leq i<\omega$, and denote by  ${\bf L}$ the subclass $\{X_{(0,\omega)}|\,\,X\in{\bf WF},\,\, X_i=X,\,\,0\leq i<\omega\}$ which consists of all such self-similar skands $X_{(0,\omega)}$ with $X_i=X$, $0\leq i<\omega$. Clearly, ${\bf L}\subset{\bf T}$. Since ${\bf WF}$ is a proper class then ${\bf L}$, which is bijective to it, is also a proper class. Consequently, ${\bf T}$ is also a proper class because on the contrary, if  ${\bf T}$ were a set  the relation ${\bf L}\subset{\bf T}$ would imply that ${\bf L}$ were a set. $\Box$

All these results are well-known. What is new is the simple proofs of them and the unexpected discovery  that the classical set-theoretical paradoxes are not paradoxes at all. Moreover, we maintain that the substance of all  such \grqq paradoxial sets", actually proper classes,  is in the following simplest lemma of such a kind of proposition.

{\bf Lemma 1.} The collection ${\bf S}$ of all singletons in ${\bf WF}$ is a proper class, not a set.

{\bf Proof.} Suppose the contrary, and ${\bf S}$ is a set. By $({\bf I})$, $\{{\bf S}\}$ and $\{\{{\bf S}\}\}$ are also well-founded  sets. Moreover, $\{{\bf S}\}\notin {\bf S}$ because all elements of ${\bf S}$ are well-founded. Then, by $({\bf II})$, the disjoint union ${\bf S}\cup\{\{{\bf S}\}\}$ is also a set and  ${\bf S}$ is a proper subset of ${\bf S}\cup\{\{{\bf S}\}\}$ which contradicts the  maximality of ${\bf S}$. $\Box$

{\bf Corollary 1.} That propositions 2, 6-10 together with Generalized Lemma 1, i.e., the collection of all singletons in ${\bf V}[{\cal U}]^-$ form a proper class, not a set, is a consequence of Lemma 1.

The {\bf Proof} is evident. One can easily find in all these cases  that ${\bf S}$ is a  subclass  of the supposed \grqq sets" or that there is a subclass of the supposed \grqq sets" which is in bijection with ${\bf S}$. Then the assumption that the wider class is a set implies, by an axiom in $NBG$ which says that the intersection of any set $X$ with an abitrary class $Y$ is a set (see, \cite{l7}, Chap. 4, Axiom S), that ${\bf S}$ is a set, contrary to Lemma 1. Indeed, the simplest proof is of the Generalized Lemma 1, Proposition 7, Proposition 9, because evidently the class of all well-founded singletons is a subclass of all singletons as well as ${\bf S}\subset {\bf V}^-$ and ${\bf S}\subset {\bf WF}$, and the statements made of them hold. Since, by von Neumann's axiom ${\bf N}$, ${\bf V}^-\approx{\bf On}$ there is a subclass of ${\bf N}$, ${\bf On}$ which is equivalent to ${\bf S}$, and by the injection $j:{\bf On}\rightarrow{\bf Card}$ (Remark 9), the statements of Propositions 6 and 7 hold.  As to Russell's set $R$ (Proposition 2), it contains ${\bf S}$ as a subclass because, for each well-founded set $X\in R$, the singleton $\{X\}$ is also well-founded and $\{X\}\in R$; otherwise, $\{X\}\in \{X\}$ and hence $\{X\}=X$ and the latter is a non-well-founded set.  In ${\bf T}$ there is a subclass of self-similar skands $X_{(0,\omega)}$ whose components $X_i=X$, $0\leq n<\omega$,  $X\in{\bf S}$, and thus Proposition 10 holds. $\Box$

{\bf Remark 10.} In \cite{l98} there is another, also very short proof of a generalized Lemma 1, but \grqq from the top to the bottom", whereas in our proof we are going \grqq from  the bottom to the top". Here it is.  If ${\bf S}$ in ${\bf V}[{\cal U}]^-$ were a set, its union would be a set. But $\bigcup{\bf S}$ contains the proper class ${\bf V}[{\cal U}]^-$ of all sets. This is because for all sets $a$, $\{a\}\in{\bf S}$. Since $\bigcup{\bf S}$ is both a proper class and a set, we have a contradiction. (See \cite{l98}, p. 338, and note that in \cite{l98} it has already been  proved that the class of all sets is a proper class, although via Russell's paradox; see Ibidem, p. 16). $\Box$

\bigskip

\parindent=0 cm
{\bf 7. Comparison of $NBG^{(1)}$ with some other approaches to non-well-founded sets}

\parindent=0,5cm
First of all, notice that objects $X^{(1)}$ in ${\bf V}[{\cal U}]^{(1)}$ are essential extensions of some but not all extraordinary sets in the sense of Mirimanoff \cite{l2}, \cite{l3}, \cite{l4} as well as non-well-founded sets in \cite{l25}. In \cite{l4}, p. 33, Mirimanoff considered extraordinary sets in the following form:
\begin{equation}
\label{f55}
E=\{y,z,...a,b,c,...\},
\end{equation}
where sets $y,z,...$ depend themselves on $E$; in particular, he considered sets in the form 
\begin{equation}
\label{f56}
E=\{E,a,b,c,...\}.
\end{equation}

In  $NBG^{(1)}$ we consider not only the latter form and those which are self-similar skands, or circular skands of period $1$, but also skands or extraordinary sets of the form 
\begin{equation}
\label{f57}
E=\{a,b,c,...,\{E, x,y,z,...\}\},
\end{equation}
 i.e.,  circular skands of period $2$, and of course we consider also  circular skands of arbitrary finite period $n$. But we {\it do not consider} extraordinary sets, or self-similar skands, e.g., of the form 
 \begin{equation}
\label{f57}
E=\{E,a,b,c,...,\{E\}\}
\end{equation}
 which is a particular case of $(\ref{f55})$. This restriction of ours is done knowingly and it is similar to the restriction of well-founded sets in comparison with non-well-founded sets which can be called $0$-rank self-similar skands (well-founded sets), $1$-rank self-similar skands (circular sets of period $1$), and we could postulate the existence of self-similar skands of arbitrary finite rank.  For our purposes we need such a restriction. So, skands are not extensions of  {\it all} of Mirimanoff's extraordinary sets of the form $(\ref{f55})$.

On the other hand, the definition $(\ref{f55})$ provides the  possibility of the existence of skands whose lengths are greater than $\omega$. However nothing was said in \cite{l4} about such things, and it appears that Mirimanoff tacitly supposed  that the  length of skands was equal to $\omega$. Only a single footnote on page 12 in \cite{l25} says that \grqq $\in_\mu$-constituent of relation (i.e., $\in$-chains of length $\mu$) can be defined for an arbitrary ordinal number $\mu$, but we consider in this paper only finite numbers $\mu$". Nevertheless, this footnote, given by Elkund, does not determine the fact that non-well-founded objects, presented in \cite{l25}, are skands of $\omega$-length in our terminology, and hence non-well-defined. So, skands in $NBG^{(1)}$ are essential extensions of extraordinary sets of length $\omega$. $\Box$
$$
$$

Now we want to compare $NBG^{(1)}$ with Aczel's model or theory of non-well-founded sets \cite{l93} which was motivated by Robin Milner's  work \cite{l70} on computer science modeling of concurent processes. Aczel's model of so-called {\it hypersets} was successfully applied to the treatment of the Liar paradox \cite{l8} and various other vicious circle phenomena \cite{l98}.

Starting out with the Zermelo-Frankel set axiomatic system ${\bf ZFC}$, which includes the axiom of choice and the axiom of foundation, Aczel rejected the latter and proposed his own {\bf Anti-Foundation Axiom} ({\bf AFA} for short); and with a natural correction of the {\bf Axiom of Extensionality} he obtained the set axiomatic system ${\bf ZFA}^-+{\bf AFA}$, which is consistent on the assumption that ${\bf ZFC}$ is consistent.

The binary relation $x\in y$ Aczel pictures as an element of a graph $x\longrightarrow y$ with nodes $x,y$ and edge $(x,y)$ between them; finite $x_0\longrightarrow x_1\longrightarrow...\longrightarrow x_{n-1}\longrightarrow x_n$ and infinite sequences finite $x_0\longrightarrow x_1\longrightarrow x_2\longrightarrow...$ as well as pointed and accesible graphs, labelled graphs, children, decorations of graphs, etc., are understood in the usual way. So, every set with an $\in$-relation can be pictured by labelled graphs. 

Then the {\bf Labelled Anti-Foundation Axiom} says:  {\it Every labelled graph has a unique labelled decoration}. 

This axiom is a natural extension of Mostowski's Collapsing Lemma that tells us that every well-founded labeled graph has a unique decoration in ${\bf WF}$. On the other hand, ${\bf AFA}$ is a strong restriction of ${\bf V}[{\cal U}]^-$ because of Scott's axiom: {\it For any extensional relation $R$ on $A$ and for any not $R$-well-founded $x\in A$ there is no set containing all possible images of $x$ under isomorphisms between $(A,R)$ and a transitive structure $(T,\in)$}. (Unpublished paper of 1960). As a consequence of this axiom the uniqueness of a \grqq Mostowski collapse" cannot be consistently postulated for non-well-founded sets. (See details in \cite{l94}.) 

There is an equivalent formulation of ${\bf AFA}$ which is nearer to our paper called the {\it  Solution Lemma}. 

 {\bf The Anti-Foundation Axiom}: {\it Every flat system of equations has a unique solution}, \cite{l98}, p. 72. 

It is enough for us because a particular example of a flat equation is $({\ref{f1}})$. Thus, axiomatic systems ${\bf ZFA}^-+{\bf AFA}$ and $NBG^{(1)}$ are incompatible, although they have mutual objects such as the empty skand ${\bf e}_{(0,\omega)}=\{\{...\}\}$ in our paper and ${\bf \Omega}=\{\{...\}\}$ in \cite{l93}. Moreover, the hyperset ${\Omega}^*=\{{\bf\Omega},{\bf \Omega}^*\}$ is equal to ${\bf \Omega}$ because ${\bf \Omega}=\{{\bf \Omega}\}=\{{\bf \Omega},{\bf \Omega}\}$ and we apply the Solution Lemma. We will see below that the generalized skand $\{{\bf e}_{(0,\omega)},\{{\bf e}_{(0,\omega)},\{...\}\}\}$ is not equal to ${\bf e}_{(0,\omega)}$.

Analysing ${\bf AFA}$ we see that a unique solution of a flat equation $({\ref{f1}})$ is actually a proper class of all its solutions in ${\bf V}[{\cal U}]^{(1)}$. Thus quotienting ${\bf V}[{\cal U}]^{(1)}$ by the solution lemma relation, we obtain a subclass of the hyperset universe. That is why Russell's paradox and other set-theoretical paradoxes (except Cantor's  and Hilbert's paradoxes) are {\it possible} in the hyperset universe (and only in such systems) of course, on the assumption that $R$ or other considered classes are sets, and, what is especially important, that one has to refer to ${\bf AFA}$ in implications such as, e.g., $R\notin R\Longrightarrow R\in R$. $\Box$

\bigskip

\parindent=0 cm
{\bf 8. Generalized skands and the universe ${\bf U}^{(\Omega)}$}

\parindent=0,5 cm
We have already enriched ${\bf V}[{\cal U}]={\bf V}[{\cal U}]^{(0)}$ by new  objects and have extended it to the class ${\bf V}[{\cal U}]^{(1)}$ of sets $X$ which can be well-founded and can be non-well-founded, i.e., have skands of length greater than or equal to $\omega$ as their elements. We can successively enrich ${\bf V}[{\cal U}]^{(1)}$ by defining    {\it generalized skands}. Moreover,  there are sequential embeddings 
${\bf V}[{\cal U}]^{(0)}\subset{\bf V}[{\cal U}]^{(1)}\subset{\bf V}[{\cal U}]^{(2)}\subset...\subset{\bf V}[{\cal U}]^{(\lambda)}\subset...\subset{\bf V[{\cal U}]^-}$, where  
$\lambda$ is any ordinal number in ${\bf On}$. 

This can be done by the following transfinite induction. 

Let  us suppose that for each ordinal number $1\leq\nu<\lambda$ we have already defined  the corresponding extensions ${\bf V}[{\cal U}]^{(\nu)}$, ${\bf U}^{(\nu)}$, and $NBG^{(\nu)}$ of ${\bf V}[{\cal U}]$, ${\bf U}$, and $NBG$, respectively, including the classes of \grqq $\nu$-pseudo-individuals" ${\cal U}^{(\nu)}$, axioms of strong extensionality in $NBG^{(\nu)}$,  the $\nu$-pseudo-well-foundness of elements of ${\bf V}[{\cal U}]^{(\nu)}$, and the following embeddings 
${\bf V}[{\cal U}]^{(0)}\subset{\bf V}[{\cal U}]^{(1)}\subset{\bf V}[{\cal U}]^{(2)}\subset...\subset{\bf V}[{\cal U}]^{(\nu)}\subset...\subset{\bf V[{\cal U}]^-}$. We are going to define the universe ${\bf U}^{(\lambda)}$ and the set theory $NBG^{(\lambda)}$.

{\bf Definition 5.}  
 By a $\lambda$-{\it generalized skand} $X^{(\lambda)}_{(\alpha_0,\alpha)}$ of length $l=\alpha-\alpha_0$ we understand   a system of embedded curly braces $\{_{\alpha_0}\{_{\alpha_0+1}...\{_{\alpha'}..._{\alpha'}\}..._{\alpha_0+1}\}_{\alpha_0}\}$, indexed by ordinals $\alpha'\in(\alpha_0,\alpha)$, whose component $X^{(\lambda)}_{\alpha'}$, for each fixed $\alpha_0\leq\alpha'<\alpha$, is either trivial (empty) or equal to $\{_{\alpha'} X^{(\nu_1)},X^{(\nu_2)},...,X^{(\nu_\mu)},...,\{_{\alpha'+1}$ or in the case, when $\alpha'=\alpha-1$, to $\{_{\alpha'} X^{(\nu_1)},X^{(\nu_2)},...,X^{(\nu_\mu)},..._{\alpha'}\}$,  for some  elements $X^{(\nu_1)},X^{(\nu_2)},$ $...,X^{(\nu_\mu)},...$ of  $\bigcup\limits_{\nu<\lambda}{\bf  V}[{\cal U}]^\nu\cup{\cal U}\cup{\cal P}\bigcup\limits_{\nu<\lambda}{\bf  V}[{\cal U}]^\nu$ such that 
 \begin{equation}
 \label{f999} \{X^{(\nu_1)},X^{(\nu_2)},...,X^{(\nu_\mu)},...\}\in\bigcup\limits_{\nu<\lambda}{\bf V[{\cal U}]^{(\nu)}}\cup{\cal P}\bigcup\limits_{\nu<\lambda}{\bf  V}[{\cal U}]^\nu,
 \end{equation}
where  $0\leq\nu_\mu<\lambda$,  $\mu\in{\bf On}$, and ${\cal P}\bigcup\limits_{\nu<\lambda}{\bf  V}[{\cal U}]^\nu$ is the class of all {\it subsets} (not proper subclasses) of $\bigcup\limits_{\nu<\lambda}{\bf  V}[{\cal U}]^\nu$.

{\bf Definition 6.}
Two $\lambda$-generalized skands $X^{(\lambda)}_{(\alpha_0,\alpha)}$ and $Y^{(\lambda)}_{(\beta_0,\beta)}$ are {\it equal} if the segments $({\alpha_0,\alpha})$ and $({\beta_0,\beta})$ are isomorphic as well-ordered sets and the corresponding components  
\begin{equation}
\label{f100}
X^{(\lambda)}_{\alpha'}=\{_{\alpha'} X^{(\nu_1)},X^{(\nu_2)},...,X^{(\nu_\mu)},...,\{_{\alpha'+1}
\end{equation}
and 
\begin{equation}
\label{f200}
Y^{(\lambda)}_{\beta'}=\{_{\beta'}Y^{(\nu_1)},Y^{(\nu_2)},...,Y^{(\nu_\mu)},...,\{_{\beta'+1},
\end{equation}
where this unique isomorphism is given by $\varphi:(\alpha_0,\alpha)\rightarrow(\beta_0,\beta)$ and $\beta'=\varphi(\alpha')$, $\alpha_0\leq\alpha'<\alpha$, $\beta_0\leq\beta'<\beta$, are equal as sets of ${\bigcup\limits_{\nu<\lambda}\bf V[{\cal U}]}^{(\nu)}\cup{\cal P}\bigcup\limits_{\nu<\lambda}{\bf  V}[{\cal U}]^\nu$; i.e., 
\begin{equation}
\label{f300}
\{X^{(\nu_1)},X^{(\nu_2)},...,X^{(\nu_\mu)},...\}=\{X^{(\nu_1)},X^{(\nu_2)},...,X^{(\nu_\mu)},...\}.
\end{equation}

 We want to   define the class of sets ${\bf V}[{\cal U}]^{(\lambda)}$, the universe  ${\bf U}^{(\lambda)}$, and the set theory $NBG^{(\lambda)}$, respectively. For this purpose we denote by ${\cal U}^{(\lambda)}$ the class of all $\lambda$-generalized skands $X^{(\lambda)}_{(\alpha_0,\alpha)}$ such that $X^{(\lambda)}_{(\alpha_0,\alpha)}\notin \bigcup\limits_{\nu<\lambda}{\bf V}[{\cal U}]^{(\nu)}$, e.g.,
 
 \parindent=0 cm
  $X^{(\lambda)}_{(0,\lambda)}=\{X^{(0)},\{X^{(1)},\{...\{X^{(\nu)},\{...\}\}\}\}\}$ or
  
   $X^{(\lambda)}_{(0,\omega)}=\{X^{(0)},X^{(1)},...,X^{(\nu)},...,\{X^{(0)},X^{(1)},...,X^{(\nu)},...,\{...\}\}\}$,
   
    where in both cases $X^{(\nu)}\in{\cal U}^{(\nu)}$, $0\leq\nu<\lambda$, and $\lambda$ is a limit ordinal; or
   
    $X^{(\lambda)}_{(\alpha_0,\alpha)}=\{X^{(\lambda-1)}_{(\alpha_0,\alpha)},\{X^{(\lambda-1)}_{(\alpha_0,\alpha)},\{...\}\}\}$, 
    
    where $X^{(\lambda-1)}\in{\cal U}^{(\lambda-1)}$, and is not a limit ordinal. 
    
\parindent=0,5 cm    
    Note that as above we consider a  $\lambda$-generalized skand $X^{(\lambda)}_{(\alpha_0,\alpha)}$ as the set whose elements are elements of the $\alpha_0$-component $X^{(\lambda)}_{\alpha_0}$ of $X^{(\lambda)}_{(\alpha_0,\alpha)}$ and one more element $X^{(\lambda)}_{(\alpha_0+1,\alpha)}$.  
    
    All $\lambda$-generalized skands of finite length as well as those $\lambda$-generalized skands, which are as sets elements of $\bigcup\limits_{\nu<\lambda}{\bf V}[{\cal U}]^{(\nu)}\cup{\cal P}\bigcup\limits_{\nu<\lambda}{\bf  V}[{\cal U}]^\nu$ do not enrich $\bigcup\limits_{\nu<\lambda}{\bf V}[{\cal U}]^{(\nu)}\cup{\cal P}\bigcup\limits_{\nu<\lambda}{\bf  V}[{\cal U}]^\nu$; moreover, on the class of {\it all such} $\lambda$-generalized skands  we define the following {\it equivalence relation}: $X^{(\lambda)}_{(\alpha_0,\alpha)}\sim Y^{(\lambda)}_{(\beta_0,\beta)}$ iff sets $X^{(\lambda)}_{\alpha_0}\cup\{X^{(\lambda)}_{(\alpha_0+1,\alpha)}\}$ and $Y^{(\lambda)}_{\beta_0}\cup\{Y^{(\lambda)}_{(\beta_0+1,\beta)}\}$ are equal as sets in $\bigcup\limits_{\nu<\lambda}{\bf V[{\cal U}]^{(\nu)}}\cup{\cal P}\bigcup\limits_{\nu<\lambda}{\bf  V}[{\cal U}]^\nu$.

Then denote by ${\bf V}[{\cal U}]^{(\lambda)}$ the class of sets {\it generated} by $\bigcup\limits_{\nu<\lambda}{\bf V[{\cal U}]^{(\nu)}}\cup{\cal P}\bigcup\limits_{\nu<\lambda}{\bf  V}[{\cal U}]^\nu$ and ${\cal U}^{(\lambda)}$. More precisely, 

$a$) all  sets of $\bigcup\limits_{\nu<\lambda}{\bf V}[{\cal U}]^{(\nu)}\cup{\cal P}\bigcup\limits_{\nu<\lambda}{\bf  V}[{\cal U}]^\nu$ are members of ${\bf V}[{\cal U}]^{(\lambda)}$, i.e., $\bigcup\limits_{\nu<\lambda}{\bf V[{\cal U}]^{(\nu)}}\cup{\cal P}\bigcup\limits_{\nu<\lambda}{\bf  V}[{\cal U}]^\nu\subset{\bf V}[{\cal U}]^{(\lambda)}$; 

$b$) all $\lambda$-generalized skands $X^{(\lambda)}_{(\alpha_0,\alpha)}\in{\cal U}^{(\lambda)}$ considered as  sets  are elements of ${\bf V}[{\cal U}]^{(\lambda)}$; 

$c$) the results of all set-theoretic operations on sets  (i.e., unions, intersections, power-sets, images of sets under mappings, products, coproducts, function sets, inverse and direct limits, etc.) and all possible constructions or implicit sets, whose constituents are elements of $\bigcup\limits_{\nu<\lambda}{\bf V}[{\cal U}]^{(\nu)}\cup{\cal P}\bigcup\limits_{\nu<\lambda}{\bf  V}[{\cal U}]^\nu$ or of
 ${\cal U}^{(\lambda)}$, are elements of ${\bf V}[{\cal U}]^{(\lambda)}$.
 
Note that $\bigcup\limits_{\nu<\lambda}{\bf V[{\cal U}]^{(\nu)}}\cup{\cal P}\bigcup\limits_{\nu<\lambda}{\bf  V}[{\cal U}]^\nu=\bigcup\limits_{\nu<\lambda}{\bf V[{\cal U}]^{(\nu)}}={\bf V[{\cal U}]^{(\lambda-1)}}$ if $\lambda$ is not a limit ordinal, but $\bigcup\limits_{\nu<\lambda}{\bf V[{\cal U}]^{(\nu)}}\cup{\cal P}\bigcup\limits_{\nu<\lambda}{\bf  V}[{\cal U}]^\nu\not=\bigcup\limits_{\nu<\lambda}{\bf V[{\cal U}]^{(\nu)}}$ if $\lambda$ is  a limit ordinal.

One can see that for sets and subclasses of ${\bf V}[{\cal U}]^{(\lambda)}$ all axioms (including the axiom ${\bf N}$) of $NBG^{(\nu)}$, $0\leq\nu<\lambda$, (except the Axiom of Extensionality for sets and classes) are satisfied. To obtain $NBG^{(\lambda)}$, we add to $NBG^{(\nu)}$, $0\leq\nu<\lambda$, the {\bf Axiom of $\lambda$-Generalized Skand Existence},  which is accepted as Definition 5, and refine the axiom of extensionality for sets and classes  in the following way:
 
 {\bf The Axiom of Strong Extensionality} in $NBG^{(\lambda)}$. Two sets (resp., subclasses) $X^{(\lambda)}$ and $Y^{(\lambda)}$ in ${\bf V}[{\cal U}]^{(\lambda)}$ (resp., of ${\bf V}[{\cal U}]^{(\lambda)}$) are equal if for each element $x\in X^{(\lambda)}$ there is an element $y\in Y^{(\lambda)}$ such that $x=y$ and for each element $y\in Y^{(\lambda)}$  there is an element $x\in X^{(\lambda)}$  such that $y=x$, where \grqq$=$" means the following:

 $1)$ the equality of individuals, when $x,y\in {\cal U}\stackrel{def}{=}{\cal U}^{(0)}$; 
 
 $2$) the equality of sets, when $x,y\in\bigcup\limits_{\nu<\lambda}{\bf V}[{\cal U}]^{(\nu)}\cup{\cal P}\bigcup\limits_{\nu<\lambda}{\bf  V}[{\cal U}]^\nu$;
 
 $3)$ the equality of skands, when $x,y\in {\cal U}^{(\lambda)}$; 
 
 $4)$ the (iterative) equality  of sets in ${\bf V}[{\cal U}]^{(\lambda)}$, i.e., the equalities of members of members of $x$ and $y$, respectively, etc. 
 
 This axiom is correct and complete since ${\bf V}[{\cal U}]^{(\lambda)}$ turns out to be {\it $\lambda$-pseudo-well-founded}; i.e., for each element $x\in{\bf V}[{\cal U}]^{(\lambda)}$  an arbitrary finite $\in$-chain starting with $x$ terminates after a finite number of steps in the following sense: $x\owns x_1\owns...\owns x_{n-1}\owns x_n$, where $x_n=\emptyset$,   or $x_n\in{\cal U}$,    or $x_n\in {\cal U}^{(\nu)}$  and $x_{n-1}\notin {\cal U}^{(\nu)}$, $1\leq\nu\leq\lambda$. Thus, the iteration (a descending $\in$-chain) for each element $x$ is finite (modulo \grqq $\nu$-pseudo-individuals", $1\leq\nu\leq\lambda$) and we definitely know that $X^{(\lambda)}=Y^{(\lambda)}$ and vice versa.
 
 That is why we can call elements of ${\cal U}^{(\lambda)}$ \grqq $\lambda$-pseudo-individuals" as above, because they along with the empty set $\emptyset$,  individuals ${\cal U}$ and $\nu$-pseudo-individuals ${\cal U}^{(\nu)}$, $1\leq\nu<\lambda$, are constituents of sets in ${\bf V}[{\cal U}]^{(\lambda)}$.
 At last, we put ${\bf U}^{(\lambda)}={\bf V}[{\cal U}]^{(\lambda)}\cup{\cal U}$. 
 
 This construction gives the following theorem.

 {\bf Theorem 2.} {\it The set theory $NBG^{(\lambda)}$ is consistent on the assumption that $NBG^{(\nu)}$ is consistent, for each $0\leq\nu<\lambda$.}

The {\bf  Proof} is similar to that for Theorem 1. 
$\Box$

The sequential embeddings 

${\bf V}[{\cal U}]^{(0)}\subset{\bf V}[{\cal U}]^{(1)}\subset{\bf V}[{\cal U}]^{(2)}\subset...\subset{\bf V}[{\cal U}]^{(\lambda)}\subset...\subset{\bf V}[{\cal U}]^{(\Omega)}\subset{\bf V[{\cal U}]^-}$, 

\parindent=0 cm where
${\bf V}[{\cal U}]^{(\Omega)}=\bigcup\limits_{\lambda\in{\bf On}}{\bf V}[{\cal U}]^{(\lambda)}$, are evident. We also put ${\bf U}^{(\Omega)}={\bf V}[{\cal U}]^{(\Omega)}\cup{\cal U}$ and consider the corresponding $NBG^{(\Omega)}$-type theory as the \grqq union" of $NBG^{(\lambda)}$, $\lambda\geq 0$, $\lambda\in {\bf On}$.

\parindent=0,5 cm 
 
{\bf Proposition 11}. Let $X^{(\Omega)}_{(\alpha_0,\alpha)}$ be a generalized skand whose components $X^{(\Omega)}_{\alpha'}$, $\alpha_0\leq\alpha'<\alpha$, are elements of ${\bf U}^{(\Omega)}$. Then there exists an ordinal $\lambda\in{\bf On}$ such that $X^{(\Omega)}_{(\alpha_0,\alpha)}\in{\bf V}[{\cal U}]^{(\lambda)}$.

{\bf Proof}. Under  assumptions on components $X^{(\Omega)}_{\alpha'}$ of  $X^{(\Omega)}_{(\alpha_0,\alpha)}$, $\alpha_0\leq\alpha'<\alpha$, we have
\begin{equation}
\label{f888}
\begin{array}{cc} X^{(\Omega)}_{(\alpha_0,\alpha)}=\{X^{(\nu_{(\alpha_0,0)})},X^{(\nu_{(\alpha_0,1)})},...,
X^{(\nu_{(\alpha_0,\mu)})}...,\{X^{(\nu_{(\alpha_0+1,0)})},X^{(\nu_{(\alpha_0+1,1)})},...,\\ X^{(\nu_{(\alpha_0+1,\mu)})},...,
\{X^{(\nu_{(\alpha',0)})},X^{(\nu_{(\alpha',1)})},...,X^{(\nu_{(\alpha',\mu)})},...,\{...\}\}\}\},
\end{array}
\end{equation}
 where $X^{(\nu_{(\alpha',\mu)})}$ are elements of ${\bf U}^{(\nu_{(\alpha',\mu)})}$,  for some $\nu_{(\alpha',\mu)}\in{\bf On}$, where $0\leq\mu<\kappa_\alpha'$ for some $\kappa_{\alpha'}\in{\bf On}$, $\alpha_0\leq\alpha'<\alpha$ (see Remark 1).
 Since for each $\alpha_0\leq\alpha'<\alpha$ components $X^{(\Omega)}_{\alpha'}$ are sets there exists $\lambda_{\alpha'}=\sup\limits_{0\leq\mu<\kappa_{\alpha'}}\nu_{(\alpha',\mu)}$, otherwise, $X^{(\Omega)}_{\alpha'}$ should be not a set, for some $\alpha_0\leq\alpha'<\alpha$. Since $(\alpha_0,\alpha)$ is a set there exists $\lambda'=\sup\limits_{\alpha_0\leq\alpha'<\alpha}\lambda_{\alpha'}$. It is clear that $\nu_{(\alpha',\mu)}\leq\lambda$, for all $0\leq\mu<\kappa_{\alpha'}$ and $\alpha_0\leq\alpha'<\alpha$. Consequently, $X^{(\nu_{(\alpha',\mu)})}\in{\bf V}[{\cal U}]^{(\lambda')}$, for all $0\leq\mu<\kappa_{\alpha'}$ and $\alpha_0\leq\alpha'<\alpha$, and hence $X^{(\Omega)}_{\alpha'}\in{\bf V}[{\cal U}]^{(\lambda)}$ and thus $X^{(\Omega)}_{(\alpha_0,\alpha)}\in{\bf V}[{\cal U}]^{(\lambda)}$, where $\lambda=\lambda'+1$. $\Box$

 {\bf Theorem 3.} {\it The set theory $NBG^{(\Omega)}$ is consistent on the assumption that $NBG^{(\lambda)}$ is consistent, for every $\lambda\in{\bf On}$.}

{\bf Proof.} We have to show only the correction of the {\bf Axiom of Strong Extensionality} in $NBG^{(\Omega)}$. For any two elements $X^{(\Omega)}$ and $Y^{(\Omega)}$ in ${\bf V}[{\cal U}]^{(\Omega)}$, by Proposition 11, there exists an ordinal $\lambda\in{\bf On}$ such that $X^{(\Omega)},Y^{(\Omega)}\in{\bf V}[{\cal U}]^{(\lambda)}$. But we know what the equality  $X^{(\Omega)}=Y^{(\Omega)}$ in $NBG^{(\lambda)}$, for each $\lambda\in{\bf On}$, means. The pseudo-well-foundedness of any element of ${\bf V}[{\cal U}]^{(\Omega)}$ is nothing other than the pseudo-well-foundedness of this element in ${\bf V}[{\cal U}]^{(\lambda)}$, for some $\lambda\in{\bf On}$. $\Box$

\parindent=0,5 cm 

{\bf Remark 11.} Proposition 11 tells us that we cannot enrich ${\bf V}[{\cal U}]^{(\Omega)}$, by adding new possible $\Omega$-pseudo-individuals ${\cal U}^{(\Omega)}$, because they are hyper-skands (see Definition 12 below), not sets. On the other hand,  ${\bf V}[{\cal U}]^{(\Omega)}\not={\bf V[{\cal U}]^-}$ because in ${\bf V[{\cal U}]^-}\setminus{\bf V}[{\cal U}]^{(\Omega)}$ there are a lot of more complicated skands, e.g., of the form $X=\{a,b,c,\{X\},X\}$, which are not considered in this paper. $\Box$

\bigskip

\parindent=0 cm
{\bf 9. Skand operations and categories of generalized skands}

\parindent=0,5cm 
Let $\alpha\geq\omega$ be an arbitrary fixed ordinal number. We consider now all generalized skands $X_{(0,\alpha)}$ whose components $X_{\alpha'}$, $0\leq\alpha'<\alpha$, are elements of ${\bf U}^{(\Omega)}$. (Actually, $X^{(\Omega)}_{(0,\alpha)}$ and $X^{(\Omega)}_{\alpha'}$, but we omit the index for simplicity.) 

It is clear, that all such skands form a proper class, not a set. We denote it here by ${\bf S}^{(\alpha)}$. Inside ${\bf S}^{(\alpha)}$ we define the following so called {\it skand-operations}.

{\bf Definition 7}. By the union $X_{(0,\alpha)}\cup Y_{(0,\alpha)}$, intersection $X_{(0,\alpha)}\cap Y_{(0,\alpha)}$, difference $X_{(0,\alpha)}\setminus Y_{(0,\alpha)}$ of two elements
$X_{(0,\alpha)}$ and $Y_{(0,\alpha)}$ of ${\bf S}^{(\alpha)}$ we understand a generalized skand $Z_{(0,\alpha)}$ whose  components $Z_{\alpha'}$ are $X_{\alpha'}\cup Y_{\alpha'}$, $X_{\alpha'}\cap Y_{\alpha'}$ and $X_{\alpha'}\setminus Y_{\alpha'}$, respectively, $0\leq\alpha'<\alpha$.
By the power  of a generalized skand $X_{(0,\alpha)}$ we understand a generalized skand 
${\cal P}X_{(0,\alpha)}$ whose  components $({\cal P}X)_{\alpha'}$ are ${\cal P}X_{\alpha'}$.

{\bf Definition 8}. By the mapping $f_{(0,\alpha)}:X_{(0,\alpha)}\rightarrow Y_{(0,\alpha)}$ of two generalized skands
$X_{(0,\alpha)}$ and $Y_{(0,\alpha)}$ we understand a generalized skand $f_{(0,\alpha)}$ whose components $f_{\alpha'}$ are the mappings $f_{\alpha'}:X_{\alpha'}\rightarrow Y_{\alpha'}$, where $X_{\alpha'}$ and $Y_{\alpha'}$ are components of $X_{(0,\alpha)}$ and $Y_{(0,\alpha)}$, respectively, $0\leq\alpha'< \alpha$. 
Moreover, $f_{(0,\alpha)}$ is an injection, surjection and bijection, respectively, if the same are all mappings $f_{\alpha'}$, $0\leq\alpha'<\alpha$. In particular, a generalized skand $X_{(0,\alpha)}$ is a {\it subskand} of $Y_{(0,\alpha)}$, if there exists an identical injection $1_{(0,\alpha)}:X_{(0,\alpha)}\rightarrow Y_{(0,\alpha)}$. 

It is clear that one can define other skand-theoretical operations such as {\it products and coproducts, inverse and direct limits, pull-backs and push-outs, the equivalence relation and quotients}, etc., by the  set-operations on the corresponding components of skands, respectively. We shall not use these constructions here, and thus omit details.

One can easily verify that the elements of ${\bf S}^{(\alpha)}$ as objects and mappings between such objects as morphisms form a category which we denote by ${\bf Sk}^{(\alpha)}$.

Of special interest for us here are categories ${\bf Sk}^{(\omega^\kappa)}(P)$, $\kappa\geq 0$, whose   objects are self-similar generalized skands $X_{(0,\alpha)}(X)$, i.e., when $\alpha=\omega^\kappa$, $\kappa\geq 0$ is fixed, and whose components are {\it permanent}, i.e., $X_{\alpha'}=X$, for each $0\leq\alpha'<\alpha$,  $X\in{\bf V}[{\cal U}]^{(\Omega)}$ and whose morphisms are self-similar  generalized skand-mappings $f_{(0,\alpha)}(f)$ of such generalized skands whose components are permanent, i.e., mappings $f_{\alpha'}=f:X\rightarrow Y$, for each $0\leq\alpha'<\alpha$. 

For $\kappa=0$, we consider the usual category of sets whose objects are elements of ${\bf V}[{\cal U}]^{(\Omega)}$; i.e., sets and morphisms are mappings of sets.

{\bf Proposition 12.} For each $X_{(0,\omega^\kappa)}(X)\in{\bf Sk}^{(\omega^\kappa)}(P)$, $\kappa\geq 1$, there is an index $\lambda$, $1\leq\lambda$, such that $X_{(0,\omega^\kappa)}(X)\in {\cal U}^{(\lambda)}$.

{\bf Proof.} Since $X\in {\bf V}[{\cal U}]^{(\Omega)}$ there exists an ordinal number $\nu\geq 0$ such that $X\in {\bf V}[{\cal U}]^{(\nu)}$. Among all such $\nu$ there is a minimal number and we denote it also by $\nu$. Then $X_{(0,\omega^\kappa)}(X)\in {\bf V}[{\cal U}]^{(\nu+1)}$ and we put $\lambda=\nu+1$.

 It is clear that {\it all} skand-operations on self-similar generalized skands give us a self-similar generalized skand as a result. We shall often use the following skand-operation:
 
 {\bf Definition 9.} Let $X_{(0,\omega^\kappa)}(X)$ ($X_{(0,\omega^\kappa)}$ for short) be a self-similar generalized skand, where $\kappa\geq 1$ is fixed. By a {\it singleton-skand} we understand the following self-similar generalized skand $S_{(0,\omega^\kappa)}(X_{(0,\omega^\kappa)})=\{X_{(0,\omega^\kappa)},\{X_{(0,\omega^\kappa)},\{...\}\}\}$.

It is clear that, by Definition 6,  the empty skand ${\bf e}_{(0,\omega^\kappa)}$ is not equal to the singleton-skand $S_{(0,\omega^\kappa)}({\bf e}_{(0,\omega^\kappa)})$ in contrast to their equality in ${\bf ZFA}^-+{\bf AFA}$ (see \cite{l93}, p. 8) as we promised to show in paragraph 7.

There is an evident functor $F_{\kappa'\kappa''}:{\bf Sk}^{(\omega^{\kappa'})}(P)\rightarrow{\bf Sk}^{(\omega^{\kappa''})}(P)$, $1\leq\kappa'\leq\kappa''$, which associates with every object $X_{(0,\omega^{\kappa'})}(X)$ the object $X_{(0,\omega^{\kappa''})}(X)$ and with every self-similar morphism $f_{(0,\omega^{\kappa'})}(f):X_{(0,\omega^{\kappa'})}(X)\rightarrow Y_{(0,\omega^{\kappa'})}(Y)$ the self-similar morphism $f_{(0,\omega^{\kappa''})}(f):X_{(0,\omega^{\kappa''})}(X)\rightarrow Y_{(0,\omega^{\kappa''})}(Y)$. A verification of the category's axioms is trivial. Moreover, all categories ${\bf Sk}^{(\omega^\kappa)}(P)$, $\kappa\geq 1$, are isomorphic to each other and hence to the usual category of sets, which are elements of ${\bf U}^{(\Omega)}$ and morphisms are maps of such sets. 

{\bf Proposition 13.}  The skand ${\bf V}_{(0,\omega^\kappa)}$ which is a skand-union of all self-similar generalized skands of length $\omega^\kappa$, i.e., ${\bf V}_{(0,\omega^\kappa)}=\bigcup\limits_{X\in{\bf V}[{\cal U}]^{(\Omega)}}X_{(0,\omega^\kappa)}(X)$, where $\kappa\geq 1$ and fixed, is not an element of ${\bf V}[{\cal U}]^{(\Omega)}$; i.e., it is a proper class, not a set.

{\bf Proof.}
Suppose the contrary, and  ${\bf V}_{(0,\omega^\kappa)}$ is a set.  Consequently, there exists  a generalized singleton-skand ${\bf S}_{(0,\omega^\kappa)}({\bf V}_{(0,\omega^\kappa)})$ of ${\bf V}_{(0,\omega^\kappa)}$ which is a set, too. Thus, there exists a singleton skand ${\bf S}_{(0,\omega^\kappa)}({\bf S}_{(0,\omega^\kappa)}({\bf V}_{(0,\omega^\kappa)}))$. Since 
${\bf S}_{(0,\omega^\kappa)}({\bf V}_{(0,\omega^\kappa)})$ is not an element of a permanent component of ${\bf V}_{(0,\omega^\kappa)}$ (otherwise, we should have a skand of the form ${\bf V}_{(0,\omega^\kappa)}=\{a,b,c,...,\{{\bf V}_{(0,\omega^\kappa)}\},{\bf V}_{(0,\omega^\kappa)}\}$ which is impossible) we obtain that ${\bf V}_{(0,\omega^\kappa)}$ is a proper subskand of

\parindent=0 cm
${\bf V}_{(0,\omega^\kappa)}\cup {\bf S}_{(0,\omega^\kappa)}({\bf S}_{(0,\omega^\kappa)}({\bf V}_{(0,\omega^\kappa)}))$ which is in contradiction with the maximality of ${\bf V}_{(0,\omega^\kappa)}$. Consequently, ${\bf V}_{(0,\omega^\kappa)}$ is not a set, but a proper class, and there is no a generalized singleton-skand ${\bf S}_{(0,\omega^\kappa)}({\bf V}_{(0,\omega^\kappa)})$ of it, which is of a great importance for us. $\Box$

\bigskip

\parindent=0 cm
{\bf 10. A new representation of ordinal and cardinal numbers}

\parindent=0,5cm

We recall once more that, due to Cantor, an ordinal number $\alpha$ (shortly an ordinal or an element of ${\bf On}$) is \grqq the ordinal type of a well-ordered set" \cite{l132}, p. 152. Likewise, a cardinal number (shortly a cardinal or an element of ${\bf Card}$) was defined by Cantor as \grqq the power type of equivalent sets" \cite{l132}, p. 87. In other words, an ordinal number $\alpha$ is a common symbol for the class of all isomorphic (\grqq similar") well-ordered sets and a cardinal number ($n\geq 0$, for natural or finite cardinal numbers and $\aleph_\nu$, $\nu\in{\bf On}$, $\nu\geq 0$, for transfinite or infinite cardinal numbers) is a common symbol for the class of equivalent sets, i.e., which are into a one-to-one correspondence. There is a natural binary relation $<$ between two ordinals (resp., cardinals) $\alpha$ and $\beta$ iff there exist well-ordered sets (resp., sets) $A$ and $B$ of the ordinal (resp., cardinal) types  $\alpha$ and $\beta$, respectively, and an initial segment $B'\subset B$ such that $A$ and $B'$ are similar (resp, equivalent). (Note that we identify here cardinal numbers with the natural and initial ordinal numbers.) 

There are diffferent eliminations of ordinals numbers as special well-founded well-ordered sets, e.g., initial segments $(0,\alpha)$ of ${\bf On}$  ordered by inclusion (Cantor \cite{l1}); the canonical representation of ordinals by pure sets $\emptyset$, $\{\emptyset\}$, $\{\emptyset,\{\emptyset\}\}$, $\{\emptyset,\{\emptyset\},\{\emptyset,\{\emptyset\}\}\}$, etc. (Mirimanoff \cite{l2}); the elimination of ordinal numbers by well-founded sets which are ordered by $\in$ and transitive (von Neumann \cite{l090}). There is a series of similar constructions which represent the class of all ordinal numbers (and hence all cardinal numbers) by non-well-founded sets which are well-ordered by $\in$ and transitive. The state of affairs is that these constructions give one more detail concerning the last ordinal and cardinal number, called here the eschaton.

{\bf Definition 10}. By an ordinal number $\alpha$ we understand an element of ${\bf V}[{\cal U}]^{(\Omega)}$ which is well-ordered by the relation $\in$ between its elements and transitive, i.e., if $X,Y,Z\in\alpha$ and $Z\in Y\in X$, then $Z\in X$.

{\bf Remark 12.} Definition 10 differs from the analogous classical definition of an ordinal number in $NBG$ because in the latter case sets $X,Y,Z$ are well-founded, as opposite to the former; one more remark: if $\alpha\in\beta$, then $\alpha<\beta$ and there is no relation $\alpha\in\alpha$ in the classical case, contrary to Definition 10, where $\alpha\in\beta$ implies in general $\alpha\leq\beta$, since an equality may exist in the case $\beta=\alpha$.

The classes ${\bf On}_{\omega^{\kappa}}$, $\kappa\geq 1$, of all ordinals in the sense of Definition 10 can be defined by the following transfinite induction.

 Let $\kappa>0$ be a fixed usual ordinal number, i.e., an element of ${\bf On}$. We begin with 
the empty skand ${\bf e}_{(0,\omega^\kappa)}=\{\{\{...\}\}\}$ and call it the first element of ${\bf On}_{\omega^{\kappa}}$,  denoting it by ${\bf e}^{(1)}$. Using Definition 9, we put successively

${\bf e}^{(1)}={\bf e}_{(0,\omega^\kappa)}$,

 ${\bf e}^{(2)}=\{{\bf e}^{(1)},\{{\bf e}^{(1)},\{...\}\}\}=\{{\bf e}^{(1)},{\bf e}^{(2)}\}$, 
 
 ${\bf e}^{(3)}=\{{\bf e}^{(1)}, {\bf e}^{(2)},\{{\bf e}^{(1)}, {\bf e}^{(2)},\{...\}\}\}=\{{\bf e}^{(1)}, {\bf e}^{(2)},{\bf e}^{(3)}\}$,
 
 .............................
 
 ${\bf e}^{(n)}=\{{\bf e}^{(1)}, {\bf e}^{(2)},...,{\bf e}^{(n-1)},\{{\bf e}^{(1)}, {\bf e}^{(2)},...,{\bf e}^{(n-1)},\{...\}\}\}=\{{\bf e}^{(1)}, {\bf e}^{(2)},...,{\bf e}^{(n-1)},{\bf e}^{(n)}\}$,
 
 .............................
 
 ${\bf e}^{(\omega)}=\{{\bf e}^{(1)}, {\bf e}^{(2)},{\bf e}^{(3)},...,\{{\bf e}^{(1)}, {\bf e}^{(2)},{\bf e}^{(3)},...,\{...\}\}\}=\{{\bf e}^{(1)}, {\bf e}^{(2)},{\bf e}^{(3)},...,{\bf e}^{(\omega)}\}$,

 ${\bf e}^{(\omega+1)}=\{{\bf e}^{(1)}, {\bf e}^{(2)},{\bf e}^{(3)},...,{\bf e}^{(\omega)},\{{\bf e}^{(1)}, {\bf e}^{(2)},{\bf e}^{(3)},...,{\bf e}^{(\omega)},\{...\}\}\}=$
 
 $\{{\bf e}^{(1)}, {\bf e}^{(2)},{\bf e}^{(3)},...,{\bf e}^{(\omega)},{\bf e}^{(\omega+1)}\}$,
 
 .....................................\,\,.

It is clear that

${\bf e}^{(1)}\in {\bf e}^{(1)}$ and ${\bf e}^{(1)}\subseteq {\bf e}^{(1)}$, ${\bf e}^{(1)}\in {\bf e}^{(2)}$ and ${\bf e}^{(1)}\subset {\bf e}^{(2)}$, 
${\bf e}^{(1)}\in {\bf e}^{(3)}$ and ${\bf e}^{(1)}\subset {\bf e}^{(3)}$,
...
${\bf e}^{(1)}\in {\bf e}^{(\omega)}$ and ${\bf e}^{(1)}\subset {\bf e}^{(\omega)}$,
...;

${\bf e}^{(2)}\in {\bf e}^{(2)}$ and ${\bf e}^{(2)}\subseteq {\bf e}^{(2)}$,
${\bf e}^{(2)}\in {\bf e}^{(3)}$ and ${\bf e}^{(2)}\subset {\bf e}^{(3)}$,
${\bf e}^{(2)}\in {\bf e}^{(4)}$ and ${\bf e}^{(2)}\subset {\bf e}^{(4)}$,...
${\bf e}^{(2)}\in {\bf e}^{(\omega)}$ and ${\bf e}^{(2)}\subset {\bf e}^{(\omega)}$,...;

and so forth.

One can see that each ordinal ${\bf e}^{(\alpha)}$ is the ordinal type of the $0$-component ${\bf e}^{(\alpha)}_0$ of the skand ${\bf e}^{(\alpha)}$. Moreover, there is a one-to-one correspondence between the class of all ordinals ${\bf On}_{\omega^{\kappa}}=\{{\bf e}^{(1)},{\bf e}^{(2)},...,{\bf e}^{(\alpha)},...\}$, $\alpha\geq 1$, and the class of all  $\alpha\in{\bf On}$, because the ordinal type of $0$-component ${\bf e}^{(\alpha)}$ is the same as the initial segment $(0,\alpha)$ of ${\bf On}$. $\Box$
$$
$$

Consider now a more interesting description of the class of all ordinals in the sense of Definition 10. It is the skand-union of all ordinal numbers, i.e., ${\bf \Omega}_{(0,\omega^\kappa)}=\bigcup\limits_{\alpha\geq 1}{\bf e}^{(\alpha)}$.

{\bf Proposition 14.}  The self-similar skand ${\bf \Omega}_{(0,\omega^\kappa)}({\bf e}^{(1)}),{\bf e}^{(2)},...)$, where $\kappa\geq 1$ and is fixed, is not an element of ${\bf V}[{\cal U}]^{(\Omega)}$; i.e., it is a proper class, not a set.

{\bf Proof.}
Suppose the contrary, and  ${\bf \Omega}_{(0,\omega^\kappa)}({\bf e}^{(1)}),{\bf e}^{(2)},...)$ (${\bf \Omega}_{(0,\omega^\kappa)}$ for short) is a set.  Consequently, there exists  a generalized singleton-skand $S_{(0,\omega^\kappa)}({\bf \Omega}_{(0,\omega^\kappa)})=\{{\bf \Omega}_{(0,\omega^\kappa)},\{{\bf \Omega}_{(0,\omega^\kappa)},\{...\}\}\}$ of ${\bf \Omega}_{(0,\omega^\kappa)}$ which is a set, too. Thus, there exists a singleton skand $S_{(0,\omega^\kappa)}(S_{(0,\omega^\kappa)}({\bf \Omega}_{(0,\omega^\kappa)}))$. Since $S_{(0,\omega^\kappa)}({\bf \Omega}_{(0,\omega^\kappa)})$ is not an element of a permanent component of ${\bf \Omega}_{(0,\omega^\kappa)}$ (otherwise, we should have a skand of the form ${\bf \omega}_{(0,\omega^\kappa)}=\{a,b,c,...,\{{\bf \Omega}_{(0,\omega^\kappa)}\},{\bf \Omega}_{(0,\omega^\kappa)}\}$ which is impossible)
 we obtain that ${\bf \Omega}_{(0,\omega^\kappa)}$ is a proper subskand of
${\bf \Omega}_{(0,\omega^\kappa)}\cup S_{(0,\omega^\kappa)}(S_{(0,\omega^\kappa)}({\bf \Omega}\}_{(0,\omega^\kappa)}))$ which is in contradiction with the maximality of ${\bf \Omega}_{(0,\omega^\kappa)}$. Consequently, ${\bf \Omega}_{(0,\omega^\kappa)}$ is not a set, but a proper class, and there is not a generalized singleton-skand $S_{(0,\omega^\kappa)}({\bf \Omega}_{(0,\omega^\kappa)})$ of it, which is of great importance for us. $\Box$

 Thus ${\bf \Omega}_{(0,\omega^\kappa)}$ with a fixed $\kappa\geq 1$  is a {\it generalized skand-class} whose permanent component is a well-ordered class ${\bf On}_{\omega^\kappa}$. 
 
Moreover, formally ${\bf \Omega}_{(0,\omega^\kappa)}\in{\bf \Omega}_{(0,\omega^\kappa)}$, which is not in contradiction with an agreement that classes are not elements of classes. In our case, with a specific definition of skand-operations, ${\bf \Omega}_{(0,\omega^\kappa)}$ cannot be an element of any class or set, e.g., the singleton $S_{(0,\omega^\kappa)}({\bf \Omega}_{(0,\omega^\kappa)})$ which does not exist at all as an element of ${\bf V}[{\cal U}]^{(\Omega)}$, but, by our natural construction, {\it it is an element of itself}. 

We see also that ${\bf \Omega}_{(0,\omega^{\kappa'})}$ and ${\bf \Omega}_{(0,\omega^{\kappa''})}$ are isomorphic for every $1\leq\kappa'<\kappa''$ and we can omit indexes; i.e., we write ${\bf \Omega}$ and call this generalized skand-class the {\it last ordinal number or the eschaton} in the sense of Definition 10, because it is well-ordered by $\in$, and is class-transitive. In other words, ${\bf\Omega}$ is a common symbol for the ordinal type of well-ordered proper classes whose all initial segments are sets. It is the last, indeed, because there are no {\it more units}, i.e., generalized singleton-skands one could add to ${\bf\Omega}$.

 It is clear that ${\bf\Omega}$ is the initial class-ordinal number because it is not equinumerous to any smaller ordinal number. Indeed, any $\alpha<{\bf\Omega}$ is a set and hence is not equivalent to ${\bf\Omega}$. By definition, the  cardinality $|A|$ of any proper class $A$ is defined as the unique class-cardinal ${\bf\Omega}$ which is equinumerous to $A$ (the existence of such equinumerousness follows from the well-ordering theorem). 
 
{\bf Proposition 15.} ${\bf\Omega}$ is   a {\it strongly inaccessible class-cardinal}, not a set cardinal.

{\bf Proof.} Let ${\bf e}^{({\omega_\nu})}<{\bf\Omega}$, where $\omega_\nu$ is the initial ordinal number. Then the power-skand ${\cal P}{\bf e}^{({\omega_\nu})}$ is a set, and hence its permanent component is not in one-to-one correspondence with the permanent  component ${\bf \Omega}_0$ of ${\bf \Omega}_{(0,\omega^\kappa)}$, $\kappa\geq 1$ is fixed; thus, ${\cal P}{\bf e}^{({\omega_\nu})}<{\bf\Omega}$. Moreover, for any ${\bf e}^{({\omega_\alpha})}<{\bf{\Omega}}$ and $\beta<{\bf\Omega}$, the sum of cardinals, i.e., the initial ordinals, $\sum\limits_{\alpha<\beta}{\bf e}^{({\omega_\alpha})}$ is a set and hence its permanent component is not in one-to-one correspondence with the permanent component ${\bf \Omega}_0$ of ${\bf \Omega}_{(0,\omega^\kappa)}$, $\kappa\geq 1$ is fixed  and hence ${\bf\Omega}$ is not its ordinal type. $\Box$

{\bf Remark 13.} The eschaton ${\bf\Omega}$ looks like the  initial ordinal $\omega$ which is also a {\it strongly inaccessible cardinal} with respect to all finite numbers, and is the first transfinite ordinal. The same can be said for ${\bf\Omega}$, which is a  strongly inaccessible class-cardinal with respect to all infinite numbers, i.e., all transfinite numbers, and is the first trans-infinite class-ordinal, or the first {\it trans-definite ordinal}, as it was called in \cite{l20}.

{\bf Definition 11.} By a {\it proper generalized skand-class} we understand $X_{(\alpha_0,\alpha)}$, at least one component $X_{\alpha'}$, $\alpha_0\leq\alpha'<\alpha$, $\alpha\in{\bf On}$ of which is a proper class, in particular, a proper generalized {\it self-similar} skand-class.

{\bf Definition 12.} By a {\it   hyper-skand} and {\it  generalized hyper-skand} we understand $X_{(\alpha_0,{\bf \Omega})}$, whose  components $X_{\alpha'}$, $\alpha_0\leq\alpha'<{\bf \Omega}$, are elements of ${\bf V}[{\cal U}]$ and ${\bf V}[{\cal U}]^{(\Omega)}$, respectively, in particular, a {\it self-similar} hyper-skand and  generalized {\it self-similar} hyper-skand $X_{(\alpha_0,{\bf \Omega})}(X)$.

Thus ${\bf O}_{(0,\omega)}=\{0,1,2,...,\omega,\omega+1,...,\{0,1,2,...,\omega,\omega+1,...,\{...\}\}\}$, 

${\bf C}_{(0,{\bf \Omega})}=\{0,\{1,\{...,
\{n,\{...\{\omega,\{\omega+1,\{...\{\alpha,\{...\}\}\}\}\}\}\}\}\}\}$ and 

${\bf E}_{(0,{\bf \Omega})}=\{{\bf e}^{(1)},\{{\bf e}^{(2)},...,
{\bf e}^{(n)},\{...\{{\bf e}^{(\omega)},\{{\bf e}^{(\omega+1)},\{...\{{\bf e}^{(\alpha)},\{...\}\}\}\}\}\}\}\}\}\}$, $\alpha\in{\bf On}$, are examples of a skand-class, a hyper-skand and a generalized hyper-skand of all ordinals, respectively.

{\bf Definition 13.} By a {\it proper $($generalized$)$ hyper-skand-class} we understand $X_{(\alpha_0,{\bf \Omega})}$, at least one component $X_{\alpha'}$, $\alpha_0\leq\alpha'<{\bf \Omega}$ of which is a proper class, in particular, a proper $($generalized$)$ {\it self-similar} hyper-skand-class.

{\bf Remark 14.} Definitions 11, 12, 13 are very general in the sense that defined objects are outside of the ${\bf U}^-$-world and are not even subclasses of it. Nevertheless, there are operations similar to operations on proper classes which are subclasses of ${\bf V}[{\cal U}]^{(\Omega)}$. On the other hand, there are no such operations as the power-skand, singleton-skand or other  set-theoretic operations. It is indeed true: \grqq The content of a concept diminishes as its extension increases; if its extension becomes all-embracing, its 
content must vanish altogether". We need these definitions for descriptions of some aspects of  the Skand Theory considered here. $\Box$

\bigskip

\parindent=0 cm
{\bf 10. Applications to $\varepsilon$-numbers}

\parindent=0,5cm
What we now want   to show are applications of skands outside of the Skand Theory considered above.

{\bf Remark 15.} The concept of a skand, i.e., objects $X_{(\alpha_0,\alpha)}$ above, is wider than its concrete realization as a system of embedded  braces and the components thereof; it is not rigidly attached to curly brackets and it may also be  a system of embedded round brackets, e.g., {\it streams} 
\begin{equation}
\label{f001}
\begin{array}{c}
s=(a_1,(a_2,(a_3,(....(a_\omega,(a_{\omega+1},(...(a_\lambda,(a_{\lambda+1},(...)))))))))), 
\end{array}
\end{equation}

\parindent=0 cm
where $a_\lambda\in A$, $A$ is a set, $1\leq\lambda<\Lambda\in{\bf On}$, (see, in particular, the  case of a countable system of embedded round brackets in \cite{l98}, p. 34-35, 197-208), or a system of embedded angle  brackets for the set theoretic operation modeling the operation of ordered pairs, or a system of embedded square brackets, e.g., {\it propositions} $f=[Fa\,\,[Fa\,\,f]]$ (the Liar proposition) or $t=[Tr\,\,[Tr\,\,t]]$ (the Truth-teller proposition) as well as $p_1=[Tr\,\,[Tr,\,\,...[Fa\,\,p_1]]]$ (the long Liar cycle of propositions $p_1,...,p_n,q$, where each proposition claims that the next one is true, except for $q$, which claims that $p_1$ is false). Note that actually there is a proper class of all the above propositions; in particular, there is a huge number of Liars, contrary to a special $\omega$-case, or inside a set-axiomatic system with ${\bf AFA}$, where there is a unique Liar, see \cite{l8}, p. 64-65. Finally, skands $X_{(\alpha_0,\alpha)}$ can be interpreted as a {\it limit power} or {\it continued exponential} with a basis which is an ordinal $\gamma_0\geq 1$ and is at the $\alpha_0$th place of $X_{(\alpha_0,\alpha)}$ and  with the exponent $X_{(\alpha_0+1,\alpha)}$, i.e., 
\begin{equation}
\label{f222}
X_{(\alpha_0,\alpha)}=\gamma_0^{X_{(\alpha_0+1,\alpha)}}, 
\end{equation}
and so on, i.e.,
a transfinite sequence of basis-exponents
\begin{equation}
\label{f002}
\begin{array}{c}
X_{(\alpha_0,\alpha)}\stackrel{sign}{=}\gamma_0^\wedge\gamma_1^\wedge\gamma_2^\wedge...^\wedge\gamma_\omega^\wedge\gamma_{\omega+1}^
\wedge...^\wedge\gamma_{\alpha'}...\stackrel{def}{=}\gamma_0^{X_{(\alpha_0+1,\alpha)}}=...=\gamma_0^{\gamma_1^{\gamma_2^{.^{.^{.^{\gamma_{\alpha'}^{.^{.^{.}}}}}}}}}, 
\end{array}
\end{equation}
where  $\gamma_{\alpha'}\in {\bf On}$, $\alpha_0+1\leq\alpha'<\alpha$, $\gamma_{\alpha'}>0$. We shall use  such a {\it skand-exponent} in application to the theory of $\varepsilon$-numbers in the sense of Cantor. For convenience in the further designation of such a skand-exponent we prefer braces $\{\{$, i.e., $X_{(\alpha_0,\alpha)}=\{\gamma_0\{\gamma_1\{...\}\}\}$ to the power-sign $^\wedge$ in $(\ref{f002})$ or expressions like $[\gamma_0,\gamma_1,\gamma_2,...\gamma_\omega,\gamma_{\omega+1},...,\gamma_{\alpha'},...]$ in \cite{l190}, and want to avoid confusion with the usual skand-set $X_{(\alpha_0,\alpha)}=\{\gamma_0,\{\gamma_1,\{...\}\}\}$, which differs from the former by commas before the second brace of a pair of opening braces, and has a different meaning.
$$
$$

\parindent=0,5 cm
Let $\gamma$, $\xi$ be arbitrary ordinal numbers such that $\gamma>0$ and $\xi\geq 0$. We recall that the $\xi$th power  of $\gamma$, i.e., $\gamma^\xi$, is defined by the following transfinite induction:
\begin{equation}
\label{f010}
\gamma^0=1,
\end{equation}

\begin{equation}
\label{f011}
\gamma^{\alpha+1}=\gamma^\alpha\gamma,
\end{equation}

\begin{equation}
\label{f012}
\gamma^\lambda=\lim\limits_{\alpha<\lambda}\gamma^\alpha,
\end{equation}
for limit-numbers $\lambda=\lim\limits_{\alpha<\lambda}\alpha$. 

And as in arithmetic, $\gamma$ is called the {\it basis}, $\xi$, the {\it exponent}, of the power $\gamma^\xi$.

Consider now the following equation:
\begin{equation}
\label{f013}
\gamma^\xi=\xi
\end{equation}
with indeterminate $\xi$.

The roots $\xi=\alpha$ of the equation $(\ref{f013})$ in the case  $\gamma=\omega$ and $\alpha<\omega_1$, where $\omega_1$ is the smallest non-denumerable ordinal, Cantor called {\it epsilon-numbers}. More precisely, \grqq to distinguish them from all other numbers I call them the \grqq$\varepsilon$-numbers of the second number-class" (\cite{l111}, \S 20). 

We can omit  these  restrictions of Cantor's, since all his results on the $\varepsilon$-numbers are valid in general cases. Here we repeat Cantor's construction in a generalized form. 

If $\alpha>0$ is any ordinal number which does not satisfy the equation $(\ref{f013})$, it determines an increasing sequence  $\alpha_n$, $0\leq n<\omega$, by means of the equalities
\begin{equation}
\label{f014}
\alpha_0=\alpha,\,\,\,\,\alpha_1=\gamma^\alpha\,\,\,\,\alpha_2=\gamma^{\alpha_1},\,\,\,\,...,\,\,\,\,\alpha_n=\gamma^{\alpha_{n-1}},\,\,\,\,...\,\,\,.
\end{equation}
Then $\lim\limits_n\alpha_n\stackrel{def}{=}\sup\limits_n\alpha_n=E(\alpha)$ of this increasing sequence always exists because ${\bf On}$ is well-ordered, and we call it an $\varepsilon$-number, too. 

Indeed, in the trivial case, when $\gamma=1$, the only root of the equation $(\ref{f013})$ is evidently $\xi=1$, and hence this actually increasing sequence in  $(\ref{f014})$ is the constant sequence $\alpha_n=1$, $0\leq n<\omega$, and thus $\lim\limits_n\alpha_n=\sup\limits_n\alpha_n=E(\alpha)=1$, for every $\alpha>1$.

If $\gamma>1$, then $(\ref{f014})$ is an ascending sequence
(in Cantor's terminology, an \grqq ascending fundamental series")
because
\begin{equation}
\label{f115}
\gamma>1\Longrightarrow\gamma^\alpha\geq\alpha,
\end{equation}
for every $\alpha\geq 0$ (see, e.g., \cite{l5}, Chap. VII, \S 6); in our case, when $\alpha>1$ and does not satisfy the equation $(\ref{f013})$, we have $\gamma^\alpha>\alpha$ and, by $(\ref{f115})$, $\gamma^{\gamma^{\alpha}}>\gamma^\alpha$, and so on. Consequently, $(\ref{f014})$ is an ascending sequence; indeed,  for all $0<n<\omega$, $\alpha_1>\alpha_0$, $\alpha_2>\alpha_1$, $\alpha_3>\alpha_2$,$\,\,$..., $\alpha_n>\alpha_{n-1}$,$\,\,$...$\,\,$. 

Put now $E(\alpha)=\lim\limits_n\alpha_n=\sup\limits_n\alpha_n$ which is a limit-ordinal number and always exists, because for the set $A=\{\alpha_0,\alpha_1,...,\alpha_n,...\}$, there exists an ordinal $\beta<\beta_0$ (for some fixed $\beta_0$), which is greater than each element of $A$ (see \cite{l5}, Chap. VII, \S 2, Theorem 6). Consequently, the least of such $\beta$, which always exists in the well-ordered set $(0,\beta_0)$ since each of its subsets has the smallest element, is the desired ordinal number $E(\alpha)$.  

By $(\ref{f011})$, $(\ref{f012})$,  the function $f(\alpha)=\gamma^\alpha$ is ascending and continuous; therefore, we have $\gamma^{E(\alpha)}=\gamma^{\lim\limits_n\gamma^{\alpha_n}}=\lim\limits_n\gamma^{\alpha_n}=\lim\limits_n\alpha_{n+1}=E(\alpha)$; i.e., $E(\alpha)$ satisfies $(\ref{f013})$.

Cantor considered the case $\gamma=\omega$, $\alpha=1$ and proved that 
\begin{equation}
\label{f215}
E(1)=\lim\limits_n\omega_n,
\end{equation}
 where
\begin{equation}
\label{f015}
\omega_1=\omega\,\,\,\,\omega_2=\omega^{\omega_1},\,\,\,\,...,\,\,\,\,\omega_n=\omega^{\omega_{n-1}},\,\,\,\,...\,\,\,\,,
\end{equation}
is an $\varepsilon$-number \cite{l111}, \S 20, [Theorem] A. Moreover, $\varepsilon_0=E(1)=\lim\limits_n\omega_n$
is the least of all the $\varepsilon$-numbers (\cite{l111}, \S 20, [Theorem] B). (This, of course, is true in his own sense; in our general construction there are two more $\varepsilon$-numbers $1$, when $\gamma=1$ and $\omega$, when $1<\gamma<\omega$ and $\alpha=\gamma$, which are evidently smaller than $\varepsilon_0$.)

He also showed that after the least $\varepsilon$-number, $\varepsilon_0$, there follows then the next greater one:
\begin{equation}
\label{f016}
\varepsilon_1=E(\varepsilon_0+1),
\end{equation}
and so on; i.e., there is the following formula of recursion:
\begin{equation}
\label{f017}
\varepsilon_n=E(\varepsilon_{n-1}+1),
\end{equation}
$1\leq n<\omega$. (\cite{l111}, \S 20, [Theorem] D).

The limit $\lim\limits_n\varepsilon_{\nu_n}$ of any ascending sequence $\varepsilon_{\nu_0},\varepsilon_{\nu_1},...,\varepsilon_{\nu_n},..$
of $\varepsilon$-number $\varepsilon_{\nu_n}$, $0\leq n<\omega$, is an $\varepsilon$-number, too. (\cite{l111}, \S 20, [Theorem] E).
Finally, all the totality of $\varepsilon$-numbers of the second number-class is a well-ordered set 
\begin{equation}
\label{f018}
\varepsilon_0,\,\,\,\varepsilon_1,\,\,\,...,\,\,\,\varepsilon_n\,\,\,...\,\,\,\varepsilon_\omega,\varepsilon_{\omega+1},\,\,\,...\,\,\,\varepsilon_{\alpha'},\,\,\,...
\end{equation}
of the second number-class type, and has thus the power $\aleph_1$, where
$0\leq\alpha'<\Omega$, (\cite{l111}, \S 20, [Theorem] F); i.e.,  $\Omega$ is the initial ordinal, in recent terminology, which is the first after the initial ordinal $\omega$. (Here we quote notations $\omega_n$, $n\geq 1$, and $\Omega$ in \cite{l132}, pp. 196, 199, literally; thus, do not confuse them with $\omega_\kappa$, $\kappa=1,2,...$, and ${\bf\Omega}$ below, respectively, which denote absolutely different objects.)

We also mention two more of Cantor's results. If $\varepsilon'$ is any $\varepsilon$-number, $\varepsilon''$ is the next greater $\varepsilon$-number, and $\alpha$ is any number which lies between them:
\begin{equation}
\label{f019}
\varepsilon'<\alpha<\varepsilon'',
\end{equation}
then $E(\alpha)=\varepsilon''$ (\cite{l111}, \S 20, [Theorem] C);
and if $\varepsilon$ is any $\varepsilon$-number and $\alpha$ is any number such that $1<\alpha<\varepsilon$, then $\varepsilon$ satisfies the three equations:
\begin{equation}
\label{f020}
\alpha+\varepsilon=\varepsilon;\,\,\,\alpha\varepsilon=\varepsilon;\,\,\,\alpha^\varepsilon=\varepsilon
\end{equation}
(\cite{l111}, \S 20, [Theorem] G).

Note that for the  $\varepsilon$-number $\omega$ in our general construction, for each $\alpha$ such that $1<\alpha <\omega$, $(\ref{f020})$ also holds.

All these  results of Cantor are valid in the general case of course, with natural corrections: ordinal types, the cardinality of the initial ordinals, etc. But we want more. We want with the help of self-similar skands to clarify this general situation.

 Cantor's formula $(\ref{f215})$ of the least (in Cantor's sense) $\varepsilon$-number $\varepsilon_0$ can be symbolically written in the following form: 
\begin{equation}
\label{f021}
\varepsilon_0=\omega^{\omega^{\omega^{.^{.^{.^{\omega^{.^{.^{.}}}}}}}}};
\end{equation}
i.e., by misuse of language, \grqq$\varepsilon_0$ is the power of $\omega$ whose exponent is the power of  $\omega$, whose exponent is the power $\omega$, etc., more precisely, $\omega$ times \grq the power of'$\,$", or \grqq the first \grq limit power' of $\omega$ whose exponents at each $n$th place of the skand-exponent $X_{(0,\omega)}$ is $\omega$, $1\leq n<\omega$". How else? 

Let us denote $(\ref{f021})$ by a bit shorter formula:
\begin{equation}
\label{f022}
\varepsilon_0=\omega^{\omega^{\omega^{.^{.^{.^{\omega^{.^{.^{.}}}}}}}}}=E_{(0,\omega)}(\omega)\stackrel{def}{=}\omega^{E_{(1,\omega)}(\omega)},
\end{equation}
where $E_{(0,\omega)}(\omega)$ denotes the self-similar skand-exponent of length $\omega$ whose components $E_n$, $0\leq n<\omega$, are $\omega$ (see Remarks 1, 15). Actually, we want to generalize this notation and the notion of the \grqq limit power" to the following one: 
\begin{equation}
\label{f023}
\varepsilon_\lambda=E_{(\alpha_0,\alpha)}(\gamma)\stackrel{def}{=}\gamma^{E_{(\alpha_0+1,\alpha)}(\gamma)}=\gamma^{\gamma^{\gamma^{.^{.^{.^{\gamma^{.^{.^{.^{\gamma{.^{.^{.^{\gamma^{.^{.^{.}}}}}}}}}}}}}}}}},
\end{equation}
for arbitrary ordinal numbers $\alpha=\omega^\kappa$, $\kappa\geq 1$, $\gamma>1$,  and for some possible $\lambda$. 

Actually, we have to explain the meaning of the symbol $E_{(\alpha_0,\alpha)}(\gamma)=\gamma^{\gamma^{\gamma^{.^{.^{.^{\gamma^{.^{.^{.^{\gamma{.^{.^{.^{\gamma^{.^{.^{.}}}}}}}}}}}}}}}}}$ as  a good way of describing {\it all} possible $\varepsilon$-numbers  because, by Definition 2, $E_{(\alpha_0,\omega^\kappa)}(\gamma)=E_{(\alpha_0+1,\omega^\kappa)}(\gamma)$ and therefore,
\begin{equation}
\label{f009}
\gamma^{E_{(\alpha_0,\omega^\kappa)}(\gamma)}=\gamma^{E_{(\alpha_0+1,\omega^\kappa)}(\gamma)}=E_{(\alpha_0,\omega^\kappa)}(\gamma); 
\end{equation}
i.e., $E_{(\alpha_0,\omega^\kappa)}(\gamma)$ in $(\ref{f023})$ is a root of the equation $(\ref{f013})$, for an arbitrary ordinal $\kappa\geq 1$, and hence is an $\varepsilon$-number. 

Here is an explicit explanation of this idea.

{\bf Definition 14.} By a  {\it skand-exponent}  $E_{(\alpha_0,\alpha)}$ of length $l=\alpha-\alpha_0$ we understand a system of embedded curly braces, indexed by $\alpha'\in(\alpha_0,\alpha)$, all of whose components  are one-element, moreover, for each $\alpha_0\leq\alpha'<\alpha$, $E_{\alpha'}=\{\gamma_{\alpha'}\{$ or $E_{\alpha'}=\{\gamma_{\alpha'}\}$, if $\alpha'=\alpha-1$, where   $\gamma_{\alpha'}\not=0$ and $\gamma_{\alpha'}\in{\bf On}$. (Notice that we do not put commas before the second brace of a pair of opening braces to point out that a skand-exponent is not considered as a two-element set.) If all components are equal to $\gamma$, then we write $E_{(\alpha_0,\alpha)}(\gamma)=\{\gamma\{\gamma\{...\}\}\}$.

{\bf Definition 15}.
 Two skand-exponents $E^1_{(\alpha_0,\alpha)}$ and $E^2_{(\beta_0,\beta)}$ are called {\it equal} if the segments $({\alpha_0,\alpha})$ and $({\beta_0,\beta})$ are isomorphic as well-ordered sets, where $\varphi:(\alpha_0,\alpha)\rightarrow(\beta_0,\beta)$ is this isomorphism, and the corresponding  components  $E^1_{\alpha'}$ and $E^2_{\beta'}$ are equal, for each $\beta'=\varphi(\alpha')$, $\alpha_0\leq\alpha'<\alpha$.

{\bf Definition 16.} By an {\it $\omega$-limit power} of the skand-exponent
\begin{equation}
\label{f008}
E_{(\alpha_0,\omega)}=\{\gamma_{\alpha_0}\{\gamma_{\alpha_0+1}\{...\}\}\}=\gamma_{\alpha_0}^{\gamma_{\alpha_0+1}^{.^{.^{.}}}}
\end{equation}
 we understand $1$, if $\gamma_{\alpha_0}=1$; if $\gamma_{\alpha_0}\not=1$, then  we understand  $\lim\limits_nE_{(\alpha_0,\alpha_0+n)}=\sup\limits_nE_{(\alpha_0,\alpha_0+n)}$ of the following $\omega$-sequence
\begin{equation}
\label{f007}
E_{(\alpha_0,\alpha_0+1)},E_{(\alpha_0,\alpha_0+2)},\,\,...\,\,E_{(\alpha_0,\alpha_0+n+1)},\,\,...
\end{equation}
 where $E_{(\alpha_0,\alpha_0+1)}=\gamma_{\alpha_0}$ and
\begin{equation}
\label{f017}
E_{(\alpha_0,\alpha_0+n+1)}=\gamma_{\alpha_0}^{\gamma_{\alpha_0+1}^{.^{.^{.^{\gamma_{\alpha_0+n-1}^{\gamma_{\alpha_0+n}}}}}}},
\end{equation}
for each $1\leq n<\omega$, is understood in the usual way: we descend, beginning with $\gamma_{\alpha_0+n-1}^{\gamma_{\alpha_0+n}}$,  $\gamma_{\alpha_0+n-2}^{\gamma_{\alpha_0+n-1}^{\gamma_{\alpha_0+n}}}$,...,$\gamma_{\alpha_0+1}^{\gamma_{\alpha_0+2}^{.^{.^{.^{\gamma_{\alpha_0+n}}}}}}$ $\gamma_{\alpha_0+1}^{\gamma_{\alpha_0+2}^{.^{.^{.^{\gamma_{\alpha_0+n}}}}}}$, ..., up to $\gamma_{\alpha_0}^{\gamma_{\alpha_0+1}^{\gamma_{\alpha_0+2}^{.^{.^{.^{\gamma_{\alpha_0+n}}}}}}}=E_{(\alpha_0,\alpha_0+n+1)}$.

Now, first of all, we shall prove the following lemmas.

{\bf Lemma 2.} Let $\varepsilon_0$ and $\varepsilon_1$ be the first and the second $\varepsilon$-numbers in   Cantor's sense. Then for each $\gamma$, $\varepsilon_0\leq\gamma<\varepsilon_1$, we have
\begin{equation}
\label{f30}
E_{(0,\omega)}(\gamma)=\varepsilon_1;
\end{equation}
in particular,
\begin{equation}
\label{f31}
E_{(0,\omega)}(\varepsilon_0)=\varepsilon_0^{\varepsilon_0^{\varepsilon_0^{.^{.^{.}}}}}=\varepsilon_1.
\end{equation}
Moreover, for each $\omega$-sequence of ordinals $\gamma_0,\gamma_1,...,\gamma_n,...$ such that $\varepsilon_0\leq\gamma_n<\varepsilon_1$, $0\leq n<\omega$, we have
\begin{equation}
\label{f32}
E_{(0,\omega)}\{\gamma_0\{\gamma_1\{\gamma_2\{...\}\}\}\}=\gamma_0^{\gamma_1^{\gamma_2^{.^{.^{.}}}}}=\varepsilon_1.
\end{equation}

{\bf Proof.} Since by the third equation in $(\ref{f020})$, $\gamma_0^{\varepsilon_1}=\varepsilon_1$, we have $\gamma_0^{\gamma_1}<\gamma_0^{\varepsilon_1}=\varepsilon_1$ as well as $\gamma_1^{\gamma_2}<\gamma_1^{\varepsilon_1}=\varepsilon_1$, and hence $\gamma_0^{\gamma_1^{\gamma_2}}<\varepsilon_1$. The same argument says that $\gamma_0^{\gamma_1^{\gamma_2^{.^{.^{.^{\gamma_n}}}}}}<\varepsilon_1$, for each $0\leq n<\omega$. Consequently, 
\begin{equation}
\label{f018}
E_{(0,\omega)}=\{\gamma_0\{\gamma_1\{\gamma_1...\}\}\}=\gamma_0^{\gamma_1^{.^{.^{.}}}}=\lim\limits_nE_{(0,n)}=\lim\limits_n\gamma_0^{\gamma_1^{\gamma_2^{.^{.^{.^{\gamma_n}}}}}}\leq\varepsilon_1.
\end{equation}
In particular, for $\gamma=\gamma_n$, $0\leq n<\omega$, $E_{(0,\omega)}(\gamma)=\gamma^{\gamma^{.^{.^{.}}}}\leq\varepsilon_1$. 

In spite of the fact that
$E_{(0,\omega)}(\gamma)$ satisfies the equation $(\ref{f013})$ and seems to be an $\varepsilon$-number, it is not a definition in  Cantor's sense, and it might be something different from  Cantor's classical  $\varepsilon$-numbers and $X_{(0,\omega)}(\gamma)=\gamma^{\gamma^{.^{.^{.}}}}<\varepsilon_1$. Why not? We shall show now that this is not the case. Clearly, by $\varepsilon_0\leq\gamma$,  $E_{(0,\omega)}(\varepsilon_0)\leq E_{(0,\omega)}(\gamma)$ and, by $\varepsilon_0\leq\gamma_n$, $0\leq n<\omega$,  $E_{(0,\omega)}(\varepsilon_0)\leq E_{(0,\omega)}=\gamma_0^{\gamma_1^{\gamma_2^{.^{.^{.}}}}}$. On the other hand, by  $\omega<\varepsilon_0$, we obtain
\begin{equation}
\label{f118} E_{(0,n+1)}=\omega^{\omega^{.^{.^{.^{\omega^{\varepsilon_0^{\varepsilon_0}}}}}}}<\varepsilon_0^{\varepsilon_0^{.^{.^{.^{\varepsilon_0^{\varepsilon_0^{\varepsilon_0}}}}}}}=X_{(0,n+1)}(\varepsilon_0).
\end{equation}
  Since $\varepsilon_0+1<\varepsilon_0^{\varepsilon_0}<\varepsilon_1$, by $(\ref{f016})$ and [Theorem] C, we obtain  
\begin{equation}
\label{f119}  
\begin{array}{c}
\varepsilon_1=\lim\limits_nE_{(0,n+1)}=\lim\limits_n\omega^{\omega^{.^{.^{.^{\omega^{\varepsilon_0^{\varepsilon_0}}}}}}}\leq\lim\limits_n\varepsilon_0^{\varepsilon_0^{.^{.^{.^{\varepsilon_0^{\varepsilon_0^{\varepsilon_0}}}}}}} =\lim\limits_nE_{(0,n+1)}(\varepsilon_0)\leq\\\leq\lim\limits_n\{\gamma_0\{\gamma_1\{\gamma_2\{...\{\gamma_n\}\}\}\}\}=\gamma_0^{\gamma_1^{\gamma_2^{.^{.^{.^{\gamma_n}}}}}}\leq\varepsilon_1.
\end{array}
\end{equation}
 Consequently, $E_{(0,\omega)}=\gamma_0^{\gamma_1^{\gamma_2^{.^{.^{.}}}}}=\varepsilon_1$, which completes the proof of Lemma 2. $\Box$

{\bf Lemma 3.} Let $\varepsilon'$ and $\varepsilon''$ be neighboring $\varepsilon$-numbers in   Cantor's sense. Then for each $\gamma$, $\varepsilon'\leq\gamma<\varepsilon''$, we have
\begin{equation}
\label{f300}
E_{(0,\omega)}(\gamma)=\varepsilon'';
\end{equation}
in particular,
\begin{equation}
\label{f310}
E_{(0,\omega)}(\varepsilon')=\varepsilon'^{\varepsilon'^{\varepsilon'^{.^{.^{.}}}}}=\varepsilon''.
\end{equation}
Moreover, for each $\omega$-sequence of ordinals $\gamma_0,\gamma_1,...,\gamma_n,...$ such that $\varepsilon'\leq\gamma_n<\varepsilon''$, $0\leq n<\omega$, we have
\begin{equation}
\label{f320}
E_{(0,\omega)}\{\gamma_0,\{\gamma_1\{\gamma_2\{...\}\}\}\}=\gamma_0^{\gamma_1^{\gamma_2^{.^{.^{.}}}}}=\varepsilon''.
\end{equation}

The {\bf Proof} is absolutely similar to the proof of Lemma 2.

{\bf Lemma 4.} Let $\varepsilon_0,\varepsilon_1,\varepsilon_2,...$ be an ascending $\omega$-sequence of 
$\varepsilon$-numbers. Then $E_{(0,\omega)}=\varepsilon_0^{\varepsilon_1^{\varepsilon_2^{.^{.^{.}}}}}=\lim\limits_n\varepsilon_n$, and hence is an $\varepsilon$-number in Cantor's sense.

{\bf Proof.} Since $\varepsilon_0<\varepsilon_1<\varepsilon_2<...$ we have, by $(\ref{f020})$, $E_{(0,n+1)}=\varepsilon_0^{\varepsilon_1^{\varepsilon_2^{.^{.^{.^{\varepsilon_n}}}}}}=\varepsilon_n$. Consequently, by Cantor's Theorem E, $E_{(0,\omega)}=\lim\limits_nE_{(0,n+1)}=\lim\limits_n\varepsilon_n$ is an $\varepsilon$-number. $\Box$

{\bf Lemma 5.} For any ordinal number $\gamma\geq\omega$, $E_{(0,\omega)}(\gamma)=\gamma^{\gamma^{\gamma^{.^{.^{.}}}}}$ is an $\varepsilon$-number in the sense of Cantor.
Moreover, for each increasing $\omega$-sequence of ordinals $\gamma_0,\gamma_1,...,\gamma_n,...$ such that $\omega\leq\gamma_0$, $E_{(0,\omega)}=\gamma_0^{\gamma_1^{\gamma_2^{.^{.^{.}}}}}$ is an $\varepsilon$-number in   Cantor's  sense.

{\bf Proof.} By $(\ref{f115})$, we have $\omega^\gamma\geq\gamma$ and hence $\gamma<E(\gamma)$. This is the well-known fact that for each ordinal $\gamma$, there is an $\varepsilon$-number greater than $\gamma$ (\cite{l66}, p. 327). Let $\varepsilon''$ be the least of such $\varepsilon$-numbers. Take the preceeding $\varepsilon'$ which always exists because $\varepsilon''$ cannot be a limit of $\varepsilon$-numbers, otherwise, $\varepsilon''$ could not be the smallest $\varepsilon$-number greater than $\gamma$. Obviously, $\varepsilon'\leq\gamma<\varepsilon''$. Then we apply Lemma 3 and obtain $E_{(0,\omega)}(\gamma)=\varepsilon''$. 

Let now $\gamma_0,\gamma_1,...,\gamma_n,...$ be an increasing $\omega$-sequence, i.e., $\gamma_0\leq\gamma_1\leq...\leq\gamma_n\leq...$. If it is stable, i.e., $\gamma_{n_0}=\gamma_{n_0+1}=\gamma_{n_0+2}=...$, we put $\gamma=\gamma_{n_0}$ and apply the first assertion of Lemma 5, i.e. $E_{(n_0,\omega)}=\gamma_{n_0}^{\gamma_{n_0+1}^{\gamma_{n_0+2}^{.^{.^{.}}}}}=\gamma^{\gamma^{\gamma^{.^{.^{.}}}}}=\varepsilon''$. Then, clearly, $E_{(0,\omega)}=\gamma_0^{\gamma_1^{\gamma_2^{.^{.^{.^{\gamma_{n_0-1}^{\varepsilon''}}}}}}}=\varepsilon''$.

If it is not stable, then without loss of generality we can assume that $\gamma_0<\gamma_1<...<\gamma_n<...$.
Suppose now that there are only finite $\varepsilon$-numbers between these $\varepsilon$-numbers, e.g., $\gamma_0$ and $\gamma_{n_0}$, and  denote the greatest of them  by $\varepsilon'$. Then, clearly, $E_{(n_0,\omega)}=\varepsilon''$ and $E_{(0,\omega)}=\varepsilon''$ as in the previous case. Finally, suppose that there is a ascending $\omega$-sequence $\varepsilon_{\nu_0}<\varepsilon_{\nu_1}<...$ with is cofinal to a ascending $\omega$-sequence $\gamma_0<\gamma_1<...$. Then, evidently, $\varepsilon''=\lim\limits_n\gamma_n=\lim\limits_n\varepsilon_n$ and the proof that $E_{(0,\omega)}=\varepsilon''$ is similar to that for Lemma 4. $\Box$

{\bf Remark 16.} It is not true that for an arbitrary $\omega$-sequence of ordinals
 
\parindent=0 cm $\gamma_0,\gamma_1,...,\gamma_n,...$ the limit power $E_{(0,\omega)}=\gamma_0^{\gamma_1^{.^{.^{.^{\gamma_n^{.^{.^{.}}}}}}}}$ 
 is an $\varepsilon$-number. 
 
\parindent=0,5 cm Indeed, consider the following $\omega$-sequence $\omega,2,2,...\,\,$, then 

$\omega^{2^{2^{.^{.^{.}}}}}=\lim\limits_n\omega^{2^{2^{.^{.^{.^{2}}}}}}=\omega^{\lim\limits_n{2^{2^{.^{.^{.^{2}}}}}}}=\omega^\omega$. 
And $\omega^\omega$ is not an $\varepsilon$-number. Indeed, $2^{\omega^{\omega}}=2^{\lim\limits_n\omega^n}=\lim\limits_n2^{\omega^n}=\lim\limits_n\omega^{\omega^{n-1}}=\omega^{\lim\limits_n\omega^{n-1}}=\omega^{\omega^\omega}>\omega^\omega$. The moral is that there are many
$\omega$-sequences of ordinals and  corresponding skand-exponents whose limit-powers {\it are equal} to $\varepsilon$-numbers; on the other hand, there are also a lot of $\omega$-sequences of ordinals and the corresponding skand-exponents whose limit-powers {\it are not} $\varepsilon$-numbers at all.

Remark 16 allows us to give the following

{\bf Definition 17.} By an {\it $\varepsilon$-number} we understand any ordinal number of the form
\begin{equation}
\label{f77}
\varepsilon=E_{(0,\omega)}(\gamma),\,\,\gamma\geq 1,
\end{equation}
and also   for an arbitrary set ${\cal E}$ of such $\varepsilon$-numbers in $(\ref{f77})$ its supremum, i.e.,
\begin{equation}
\label{f771}
\varepsilon'=\sup\limits_{\varepsilon\in{\cal E}}{\cal E}.
\end{equation}

It is clear that all $\varepsilon$-numbers in Definition 17 are $\varepsilon$-numbers in the sense of Cantor except two numbers: $1$ because $1^1=1$ and $(\ref{f013})$ holds, and $\omega$, because $2^\omega=\omega$ and $(\ref{f013})$ holds, too. Let us denote these first two $\varepsilon$-numbers by $\varepsilon_0$ and $\varepsilon_1$, respectively. Since $\varepsilon_2=E_{(0,\omega)}(\omega)$, by Lemma 2, we obtain that it is the least $\varepsilon$-number in the sense of Cantor; and all finite indexes of $\varepsilon$-numbers in the sense of Cantor are shifted by adding $2$. 

{\bf Definition 18.} Let $\alpha=\omega^\kappa$, $\kappa\geq 1$. Then by a {\it limit-power} with the basis $\gamma>1$ and the same exponents we understand  the skand-exponent $E_{(0,\alpha)}(\gamma)$, given by the following transfinite recursion:
\begin{equation}
\label{f51}
E_{(0,\omega)}(\gamma)=\bar\varepsilon_1, \,\,\,\kappa=1;
\end{equation}
\begin{equation}
\label{f52}
E_{(0,\omega^2)}(\gamma)=E_{(0,\omega)}(\bar\varepsilon_1)=\bar\varepsilon_2, \,\,\,\kappa=2;
\end{equation}
$$
..................................
$$
\begin{equation}
\label{f53}
E_{(0,\omega^n)}(\gamma)=E_{(0,\omega)}(\bar\varepsilon_{n-1})=\bar\varepsilon_n, \,\,\,\kappa=n;
\end{equation}
$$
..................................
$$
\begin{equation}
\label{f54}
E_{(0,\omega^\omega)}(\gamma)=\sup\limits_nE_{(0,\omega^n)}(\gamma)=\sup\limits_nE_{(0,\omega)}(\bar\varepsilon_{n-1})=\bar\varepsilon_\omega, \,\,\,\kappa=\omega;
\end{equation}
$$
..................................
$$
In the general case
\begin{equation}
\label{f551}
E_{(0,\omega^\kappa)}(\gamma)=E_{(0,\omega)}(\bar\varepsilon_{\kappa-1})=\bar\varepsilon_\kappa, \,\,\,\kappa-1<\kappa;
\end{equation}
\begin{equation}
\label{f552}
E_{(0,\omega^\kappa)}(\gamma)=\sup\limits_{\lambda<\kappa}E_{(0,\omega^\lambda)}(\gamma)=\sup\limits_{\lambda<\kappa}E_{(0,\omega^\lambda)}(\gamma)=\bar\varepsilon_\kappa, \,\,\,\not\exists\,\,\,\kappa-1.
\end{equation}

In accordance with Remark 16 and Definition 18, we are going to describe in canonical form all $\varepsilon$-numbers.

{\bf Theorem 4.} {\it There is a one-to-one correspondence between all ordinal numbers $\omega^\kappa$, $0\leq\kappa\in{\bf On}$ and all $\varepsilon$-numbers $\varepsilon_{\kappa}$, defined in Definition $17$, as follows: $\varepsilon_0=E_{(0,1)}(1)$ and $\varepsilon_\kappa=E_{(0,\omega^\kappa)}(2)$, $\kappa\geq 1$.}

{\bf Proof.} If $\kappa=0$, we put $\varepsilon_0=1$. If $\kappa>0$, putting in Definition 18 $\gamma=2$, $\bar\varepsilon_\kappa=\varepsilon_\kappa$, we obtain, by Lemmas 2,3,4,5, a successive enumeration of all $\varepsilon$-numbers in the sense of Definition 17 and thus all $\varepsilon$-numbers in the sense of Cantor. $\Box$

{\bf Corollary 2.} The set $E_{(\alpha)}=\{\varepsilon_\kappa|\,\,\varepsilon_\kappa<\omega_\alpha\}$, i.e., the set of all $\varepsilon$-numbers $\varepsilon_\kappa$ such that $\varepsilon_\kappa$ less than the  initial number $\omega_\alpha$ has a power greater than or equal to $\omega_\alpha$.

{\bf Proof.} If the power of $E_{(\alpha)}$ were less than the power of $\omega_\alpha$, then $\varepsilon=\sup\limits_{\varepsilon_\kappa\in E_{(\alpha)}}\varepsilon_\kappa$  would be greater than all of the elements in $E_{(\alpha)}$, and it would be an $\varepsilon$-number whose  power would be less than $\omega_\alpha$, because the latter is not the limit of a transfinite sequence of smaller powers. Thus $E_{(\alpha)}\cup\{\varepsilon\}$ would be larger than $E_{(\alpha)}$, which is in  contradiction with the maximality of $E_{(\alpha)}$. Consequently, $|E_{(\alpha)}|\geq|\omega_\alpha|$ ($|E_{(\alpha)}|$ and $|\omega_\alpha|$ mean the cardinality of $E_{(\alpha)}$ and $\omega_\alpha$, respectively). 

{\bf Corollary 3.} The set $E_{(\alpha)}=\{\varepsilon_\kappa|\,\,\varepsilon_\kappa<\omega_\alpha\}$, i.e., the set of all $\varepsilon$-numbers $\varepsilon_\kappa$ such that $\varepsilon_\kappa$ less than the  initial number $\omega_\alpha$ is a well-ordered set of the ordinal type $\omega_\alpha$ and thus has the power of $\omega_\alpha$.

{\bf Proof.} For each $\kappa<\omega_\alpha$, clearly, $\varepsilon_\kappa<\omega_\alpha$. And this is a one-to-one correspondence between the set $(0,\omega_\alpha)=\{\kappa|\,\,\kappa<\omega_\alpha\}$ and the set $E_{(\alpha)}$. Moreover, if $\kappa'<\kappa''<\omega_\alpha$, then $\varepsilon_{\kappa'}<\varepsilon_{\kappa''}$, i.e., the ordinal type $(0,\omega_\alpha)$ is the same as $E_{(\alpha)}$. Since by Corollary 2, $|E_{(\alpha)}|\geq|\omega_\alpha|$, we obtain that $| E_{(\alpha)}|=|\omega_\alpha|$.

{\bf Corollary 4.} Every  initial number $\omega_\alpha$ is an $\varepsilon$-number.

{\bf Proof.} Clearly, $\omega_\alpha=\lim\limits_{\varepsilon_\kappa\in E_{(\alpha)}}\varepsilon_\kappa$ and thus, by Definition 17, it is an  $\varepsilon$-number.

{\bf Remark 17.} The assertion of Corollary 4 was given in \cite{l66} and proved in the particular case $\omega_1$. It is well-known that each initial number $\omega_\alpha$ has a form $\omega_\alpha=\omega^\lambda$ for some ordinal $\lambda$ (see, e.g., \cite{l5}, Chap. VIII, \S 3, Theorem 9), but usually no one points out that $\lambda=\omega_\alpha$.

{\bf Corollary 5.} The well-ordered proper class of all $\varepsilon$-numbers is isomorphic to ${\bf On}$.

 The {\bf Proof} is similar to that for Corollary 3. $\Box$

{\bf Remark 18.} The limit power $E_{(0,\omega^\kappa)}(\gamma)$, $\kappa\geq 1$, in Definition 18 should be better called the {\it quantified power}, because actually we quantify exponents by $\omega$-sequences of the same exponents: $1<\gamma<\omega$, then $\omega$, then $\varepsilon_2$, $\varepsilon_3$,..., $\varepsilon_\omega$,..., just  to simplify the algorithm of transfinite recursion. We could do the same thing by a {\it continued exponentials} process, e.g., $E_{(0,\omega_1)}(2)$:

 $E_{(0,2)}(2)=2^2$,..., $E_{(0,n)}(2)=2^{2^{.^{.^{.^2}}}}$,...$E_{(0,\omega)}(2)=2^{2^{.^{.^{.}}}}=\omega=\varepsilon_1$;

$E_{(0,\omega +1)}(2)=\omega^2$, $E_{(0,\omega +2)}(2)=\omega^{2^2}$,... $E_{(0,\omega 2)}(2)=\omega^{2^{2^{.^{.^{.}}}}}=\omega^\omega=\varepsilon_1^{\varepsilon_1}$;

$E_{(0,\omega 2 +1)}(2)=\omega^{\omega^2}$, $E_{(0,\omega +2)}(2)=\omega^{\omega^{2^2}}$,... $E_{(0,\omega 3)}(2)=\omega^{\omega^{2^{2^{.^{.^{.}}}}}}=\omega^{\omega^{\omega}}$,...,

$E_{(0,\omega n)}(2)=\omega^{\omega^{^{.^{.^{.^\omega}}}}}$,...$E_{(0,\omega^2)}(2)=
\omega^{\omega^{\omega^{\omega^{.^{.^{.}}}}}}=\varepsilon_1^{\varepsilon_1^{\varepsilon_1^{\varepsilon_1{.^{.^{.}}}}}}=\varepsilon_2$;

.................................

$E_{(0,\omega^2+1)}(2)=\varepsilon_2^2$, $E_{(0,\omega^2+2)}(2)=\varepsilon_2^{2^2}$,..., $E_{(0,\omega^2+\omega)}(2)=\varepsilon_1^{2^{2^{.^{.^{.}}}}}=\varepsilon_2^\omega$,

$E_{(0,\omega^2+\omega+1)}(2)=\varepsilon_2^{\omega^2}$, $E_{(0,\omega^2+\omega+2)}(2)=\varepsilon_2^{\omega^{2^2}}$,...$E_{(0,\omega^22)}(2)=\varepsilon_2^{\varepsilon_2}$;

.................................

$E_{(0,\omega^3)}(2)=\varepsilon_2^{\varepsilon_2^{\varepsilon_2^{\varepsilon_2^{.^{.^{.}}}}}}=\varepsilon_3$;

.................................

$E_{(0,\omega^\omega)}(2)=\varepsilon_1^{\varepsilon_2^{\varepsilon_3^{\varepsilon_4^{.^{.^{.}}}}}}=\varepsilon_\omega$;

..................................

$E_{(0,\omega_1)}(2)=\varepsilon_1^{\varepsilon_2^{.^{.^{.^{\varepsilon_\omega^{\varepsilon_{\omega+1}^{.^{.^{.^{\varepsilon_\alpha^{.^{.^{.}}}}}}}}}}}}}=2^{E_{(1,\omega_1)}(2)}=2^{\omega_1}=\omega_1$, where $1\leq\alpha<\omega_1$.

There is a similar process in the case $\alpha=\omega_\lambda$, $\lambda\geq 2$, i.e., in the calculation of $E_{(0,\omega_\lambda)}$.

{\bf Proposition 16.} For an arbitrary ordinal number $\alpha>0$ and any $1<\gamma<\omega$, $E_{(0,\alpha)}(\gamma)$ can be expressed by the following unique formula:
\begin{equation}
\label{f888}
E_{(0,\alpha)}(\gamma)=\gamma^{E_{(1,\alpha)}(\gamma)}=\varepsilon_{\eta_1}^{\varepsilon_{\eta_1}^{.^{.^{.^{\varepsilon_{\eta_1}^{\varepsilon_{\eta_2}^{.^{.^{.^{\varepsilon_{\eta_n}^{\gamma^{.^{.^{.^{\gamma}}}}}}}}}}}}}}},
\end{equation}
where $\varepsilon_{\eta_1}>\varepsilon_{\eta_2}>...>\varepsilon_{\eta_n}$ are  $\varepsilon$-numbers, $\eta_1>\eta_2>...>\eta_n>0$ are ordinal numbers, and the quantization of them, and of $\gamma$ in the exponents is given by $\beta_1,\beta_2,...,\beta_n,\beta_{n+1}$, respectively, $0\leq\beta_i<\omega$, $i=1,2,...,n,n+1$.

{\bf Proof.} By  Cantor's normal form of $\alpha$, it may be represented uniquely as
\begin{equation}
\label{f881}
\alpha=\omega^{\eta_1}\beta_1+\omega^{\eta_2}\beta_2+...+\omega^{\eta_n}\beta_n+\beta_{n+1},
\end{equation}
where $\eta_1,\eta_2,...,\eta_n$ is a descending sequence of  ordinal numbers $>0$ and natural numbers $\beta_1,\beta_2,...,\beta_n,\beta_{n+1}$ are  $\geq 0$ (see, e.g., \cite{l5}, Chap. VII, \S 7, Theorem 2). Then, by Theorem 4 and Remark  18, we obtain $(\ref{f888})$. $\Box$

{\bf Proposition 17.} For every ordinal number $\gamma\geq 2$, the hyper-skand $E_{(0,{\bf\Omega})}(\gamma)$ is the greatest $\varepsilon$-class-number, i.e.,
\begin{equation}
\label{f882}
E_{(0,{\bf \Omega})}(\gamma)=\gamma^{E_{(1,{\bf \Omega})}(\gamma)}=\lim\limits_{\alpha\in{\bf On}}\varepsilon_{\alpha}=\varepsilon_{\bf\Omega}={\bf\Omega}=\gamma^{\bf\Omega}.
\end{equation}

{\bf Proof.} It is clear that for each ordinal number $\alpha>0$, $E_{(\alpha,{\bf \Omega})}(\gamma)$ is equal to $E_{(0,{\bf \Omega})}(\gamma)$. In other words, the skand $E_{(0,{\bf\Omega})}(\gamma)$ is self-similar, because each remainder $(\alpha,{\bf\Omega})$ as an ordered class is equal to $(0,{\bf\Omega})$, i.e., $(\alpha,{\bf\Omega})$ and $(0,{\bf\Omega})$ are isomorphic as ordered classes, in particular, for $\alpha=1$. Since, by Corollary 5, the class of all $\varepsilon$-numbers has of the same ordinal type as ${\bf\Omega}$  we obtain that $\lim\limits_{\alpha\in{\bf On}}\varepsilon_{\alpha}=\lim\limits_{\alpha\in{\bf On}}\alpha={\bf\Omega}$ and consequently, $\gamma^{\bf\Omega}=\gamma^{\lim\limits_{\alpha\in{\bf On}}\varepsilon_{\alpha}}=\lim\limits_{\alpha\in{\bf On}}\gamma^{\varepsilon_{\alpha}}=\lim\limits_{\alpha\in{\bf On}}\varepsilon_{\alpha}={\bf\Omega}$, i.e., ${\bf\Omega}$ satisfies $(\ref{f013})$, and thus ${\bf\Omega}$  is an $\varepsilon$-class-number. $\Box$

\parindent=0 cm
{\bf 11. Applications to generalized real fractions}

\parindent=0,5cm 
Likewise, we denote here a  transfinite $\alpha$-sequence $x_{\alpha_0},x_{\alpha_0+1},...,x_{\alpha'},...$, where $\alpha_0\leq\alpha'<\alpha<{\bf \Omega}$ and $x_{\alpha'}\in{\bf U}^{\Omega}$, as a special skand 
\begin{equation}
\label{f555}
X_{(\alpha_0,\alpha)}=\{x_{\alpha_0}\{x_{\alpha_0+1}\{...\{x_{\alpha'}\{...\}\}\}\}\}
\end{equation}
  without commas before the second brace of the opening pairs of braces.

\parindent=0,5cm
Let $\alpha$ be an ordinal of the $2$nd type, i.e. that having no predecessor; in particular, $0$ is an ordinal of the $2$nd type; its form $\alpha=\omega\nu$, where $\nu\geq 0$, is known. Recall also that ordinals of the $1$st type are those having predecessors. (We have already used other terminology above as a limit ordinal number and an ordinal number which is not a limit number; an ordinal of the 2nd type and of the 1st type are shorter, and we now prefer the latter terminology.)

 Consider now  for a fixed $\alpha=\omega\nu$, $\nu\geq 1$, the set $A_\alpha$  of all special skands (henceforth, in short: skands)
$X_{(0,\alpha)}$ whose components $X_{\alpha'}=0$ or $X_{\alpha'}=1$ without commas before the second brace of the opening pairs of braces as well.

We endow $A_\alpha$ with the following lexicographic linear ordering: $X_{(0,\alpha)}<Y_{(0,\alpha)}$ iff there is an $\alpha'$, $0\leq\alpha'<\alpha$, such that $X_{\alpha'}=0$ and $Y_{\alpha'}=1$ and at whichever $\beta$th place, $\beta<\alpha'$, the elements are equal; i.e., $X_\beta=Y_\beta$. If in addition   $X_\beta=1$ and  $Y_\beta=0$, for all $\beta>\alpha'$, then for  such  pairs {\it only}, there are no skands $Z_{(0,\alpha)}$ in $A_\alpha$ with  $X_{(0,\alpha)}<Z_{(0,\alpha)}<Y_{(0,\alpha)}$. We call those pairs of neighboring  skands  {\it twins}. We shall identify them and denote  the obtained new element in the canonical form, i.e., of a greater  $Y_{(0,\alpha)}$,  not forgetting that there is a different form of it, i.e., of a smaller $X_{(0,\alpha)}$, and using it when it is convenient. 

{\bf Definition 18}.
By $R_\alpha|_{[0,1]}=[0_{(0,\alpha)},1_{(0,\alpha)}]$  we denote the quotient  set $A_\alpha/_\sim$ of $A_\alpha$ ($\sim$ identifies each pair of twins as one element) with the quotient  linear ordering  and call it a {\it generalized real number (more precisely, fractional) unit interval } of the power $2^{|\alpha|}$. 

\parindent=0,5cm Here $0_{(0,\alpha)}$ and $1_{(0,\alpha)}$ are  minimal and maximal elements (integers) of $R_\alpha|_{[0,1]}$, i.e., skands with $0$ and $1$ at all places, respectively. 

{\bf Definition 19}.
  By $Q_\alpha|_{[0,1]}$ we denote the subset of $R_\alpha|_{[0,1]}$ of all skands $X_{(0,\alpha)}$  which are eventually $0$ or $1$.
  
  In particular, we distinguish in $Q_\alpha|_{[0,1]}$ {\it dyadic fractions}, i.e. $\frac{1}{2^{\alpha'}}$ as $X_{(0,\alpha)}$ such that  $X_{\alpha'-1}=1$ and $X_\beta=0$ for all $\beta\not=\alpha'-1$, for each ordinal number $\alpha'$ of the $1$st kind, $1\leq\alpha'<\alpha$, (which is a twin to $Y_{(0,\alpha)}$ with $Y_\beta=0$ for $0\leq\beta<\alpha'$, and $Y_\beta=1$, for $\alpha'\leq\beta<\alpha$) and also $\frac{1}{2^{\alpha'}}$ as $X_{(0,\alpha)}$ such that  $X_{\beta}=0$, for $0\leq\beta<\alpha'$, and $X_\beta=1$ for all $\beta\geq\alpha'$, for each ordinal number $\alpha'$ of the $2$nd kind, $0\leq\alpha'<\alpha$. In other words, in short, $\frac{1}{2^{\alpha'}}$ are skands $X_{(0,\alpha)}$ which are eventually $1$.

  {\bf Proposition 18.} $R_\alpha|_{[0,1]}$ and $Q_\alpha|_{[0,1]}$ are the dense linear orderings and  $Q_\alpha|_{[0,1]}$ is dense in $R_\alpha|_{[0,1]}$.
  
 The {\bf Proof} is an immediate consequence of Definitions 18 and 19 together with the definition of the ordering on $A_\alpha/_\sim$.
   
{\bf Theorem 5}. {\it The space $R_\alpha|{_{[0,1]}}$ is continuous: i.e., every non-empty subset $S$ of $R_\alpha|{_{[0,1]}}$ has a smallest upper bound $M_{(0,\alpha)}=\sup\,S$ and a greatest lower bound $m_{(0,\alpha)}=\inf\,S$ in $R_\alpha|{_{[0,1]}}$.}

{\bf Proof.}  
If there exists a  maximal element $\max\,S$ in $S$, then   $\sup\,S=\max\,\,S$,  if there exists a minimal  element $\min\,S$ in $S$, then   $\inf\,S=\min\,\,S$.

Consider now the case when $S$ has no maximal element and prove that  there exists $M_{(0,\alpha)}=\sup\,S$ in $R_\alpha|{_{[0,1]}}$, i.e.,  $M_{(0,\alpha)}\in R_\alpha|{_{[0,1]}}$
such that for all $X_{(0,\alpha)}\in S$ we have $X_{(0,\alpha)}<M_{(0,\alpha)}$ and for each $Y_{(0,\alpha)}\in R_\alpha|{_{[0,1]}}$ such that $Y_{(0,\alpha)}<M_{(0,\alpha)}$ there is $Z_{(0,\alpha)}\in S$ such that  $Y_{(0,\alpha)}<Z_{(0,\alpha)}$.

Indeed, there exists the smallest ordinal $\alpha_1\geq 0$ such that there is an element $X^1_{(0,\alpha)}\in S$ with  $X^1_{\alpha_1}=1$ and for every  $X_{(0,\alpha)}\in S$, $X_\beta=0$, for each  $0=\alpha_0\leq\beta<\alpha_1$ (if $\alpha_1=0$,  then conditions $X_\beta=0$, $\beta<\alpha_1$, are absent). Otherwise, $S$ should be  $\{0_{(0,\alpha)}\}$ or the empty set $\emptyset$, which is impossible by assumption.

We shall define $M_{(0,\alpha)}\in R_\alpha|{_{[0,1]}}$ by induction on its non-trivial components, putting at the beginning  $M_{\alpha_1}=1$ and $M_\beta=0$, for each $0\leq\beta<\alpha_1$, and
then define  the following subset $S^1=S\setminus\{X_{(0,\alpha)}|\,\,X_{(0,\alpha)}\in S,\,\,X_{(0,\alpha)}<X^1_{(0,\alpha)}\}$ of $S=S^0$.
 
Since $S$ has no maximal element and $X^1_{(0,\alpha)}\in S$, there exists the smallest ordinal $\alpha_2>\alpha_1$ such that there is an  element   $X^2_{(0,\alpha)}>X^1_{(0,\alpha)}$ in $S^1$ with $X^2_{\alpha_2}=1$,    and for every  $X_{(0,\alpha)}\in S^1$, $X_\beta=0$, for each  $\alpha_1<\beta<\alpha_2$. 

We continue to define $M_{(0,\alpha)}$ for the next series of indexes by putting  $M_{\alpha_2}=1$ and $M_\beta=0$, for each $\alpha_1<\beta<\alpha_2$, and define now $S^2=S^1\setminus\{X_{(0,\alpha)}|\,\,X_{(0,\alpha)}\in S^1,\,\,X_{(0,\alpha)}<X^2_{(0,\alpha)}\}$. 

Suppose that, for each $1\leq k\leq n<\omega$, we have already  found the smallest  ordinal $\alpha_k>\alpha_{k-1}$ and elements  $X^k_{(0,\alpha)}>X^{k-1}_{(0,\alpha)}$ such that $X^k_{\alpha_k}=1$ and in addition for every $X_{(0,\alpha)}\in S^{k-1}$ its components $X_\beta=0$, where   $\alpha_{k-1}<\beta<\alpha_k$. Suppose also that we have already defined the next components of $M_{(0,\alpha)}$, by putting  $M_{\alpha_k}=1$ and $M_\beta=0$, for each $\alpha_{k-1}<\beta<\alpha_k$, as well as the  set $S^k=S^{k-1}\setminus\{X_{(0,\alpha)}|\,\,X_{(0,\alpha)}\in S^{k-1},\,\,X_{(0,\alpha)}<X^k_{(0,\alpha)}\}$. Notice that in this induction we put formally  $X^0_{(0,\alpha)}=0_{(0,\alpha)}$. 

Since $S^n$ has no maximal element and $X^n_{(0,\alpha)}\in S$, there exists the smallest ordinal $\alpha_{n+1}>\alpha_n$ such that there is an element $X^{n+1}_{(0,\alpha)}>X^n_{(0,\alpha)}$ in $S^n$ with  $X^{n+1}_{\alpha_{n+1}}=1$ and for every $X_{(0,\alpha)}\in S^n$, $X_\beta=0$, for each   $\alpha_n<\beta<\alpha_{n+1}$. We put $M_{\alpha_{n+1}}=1$ and  $M_\beta=0$, for each $\alpha_n<\beta<\alpha_{n+1}$, and define $S^{n+1}=S^n\setminus\{X_{(0,\alpha)}|\,\,X_{(0,\alpha)}\in S^n,\,\,X_{(0,\alpha)}<X^{n+1}_{(0,\alpha)}\}$.

Thus $S^n$ are defined for all $0\leq n<\omega$  and we can consider their intersection $\bigcap\limits_n S^n$. If $\bigcap\limits_n S^n=\emptyset$, i.e., there are no more elements $X_{(0,\alpha)}$ in $S$ with $X_\beta=1$, $\beta>\alpha_n$, for each $0\leq n<\omega$, then we put $M_\beta=0$, for all $\alpha_\omega\leq\beta<\alpha$, where $\alpha_\omega=\lim\limits_n\alpha_n$. Since for each $0\leq\alpha'<\alpha$, $M_{\alpha'}$ has been already defined, we obtain an element $M_{(0,\alpha)}\in R_\alpha|{_{[0,1]}}$ and we claim that it is $\sup\,S$. Indeed, by construction, for every $X_{(0,\alpha)}\in S$, we have $X_{(0,\alpha)}<M_{(0,\alpha)}$. If $Y_{(0,\alpha)}\in R_\alpha|{_{[0,1]}}$ and $Y_{(0,\alpha)}<M_{(0,\alpha)}$, then there is a minimal index $\alpha_n$, $0\leq n<\omega$, such that $Y_{\alpha_n}=0$ and $M_{\alpha_n}=1$. Take $X^n_{(0,\alpha)}\in S$. It is clear that $Y_{(0,\alpha)}<X^n_{(0,\alpha)}$. Thus in this case the existence of $\sup\,S$ is proved.

If $\bigcap\limits_n S^n\not=\emptyset$, then we define $S^\omega=\bigcap\limits_n S^n$. Since $S^\omega\not=\emptyset$  there exists the smallest ordinal    $\alpha_{\omega+1}\geq\alpha_\omega$ such that there is an element  $X^{\omega+1}_{(0,\alpha)}$ in $S^\omega$ with $X^{\omega+1}_{(0,\alpha)}>X^n_{(0,\alpha)}$, for each $0\leq n<\omega$,  $X^{\omega+1}_{\alpha_{\omega+1}}=1$, and for every $X_{(0,\alpha)}\in S^\omega$, $X_\beta=0$, for each   $\alpha_\omega\leq\beta<\alpha_{\omega+1}$ (if $\alpha_{\omega+1}=\alpha_\omega$,  then conditions $X_\beta=0$, $\beta<\alpha_{\omega+1}$, are absent). We put $M_{\alpha_{\omega+1}}=1$ and $M_\beta=0$, for each $\alpha_\omega\leq\beta<\alpha_{\omega+1}$ 
and define the following set $S^{\omega+1}=S^\omega\setminus\{X_{(0,\alpha)}|\,\,X_{(0,\alpha)}\in S^\omega,\,\,X_{(0,\alpha)}<X^{\omega+1}_{(0,\alpha)}\}$. Then we continue our algorithm as above.

Since each step of our inductive construction enlarges the index $\alpha_\nu$, $1\leq\nu$, at least by $1$, we shall exhaust all of $0\leq\alpha'<\alpha$ and obtain an element $M_{(0,\alpha)}\in R_\alpha|{_{[0,1]}}$ such that $X^\nu_{(0,\alpha)}<M_{(0,\alpha)}$, for each $\nu\geq 1$. Since $\bigcap\limits_\nu S^\nu=\emptyset$ (otherwise, $M_{(0,\alpha)}$ should be an elements of $\bigcap\limits_\nu S^\nu$ and thus the greatest element of $S$) we conclude that  $M_{(0,\alpha)}=\sup\,S$. Indeed, for each $X_{(0,\alpha)}\in S$, $X_{(0,\alpha)}<M_{(0,\alpha)}$ and if $Y_{(0,\alpha)}\in R_\alpha|{_{[0,1]}}$ and $Y_{(0,\alpha)}<M_{(0,\alpha)}$, then there is the smallest ordinal $\alpha_\nu$ condidered above such that $Y_{\alpha_\nu}=0$ and $M_{\alpha_\nu}=1$. Take $X^\nu_{(0,\alpha)}\in S$ and by construction $Y_{(0,\alpha)}<X_{(0,\alpha)}$.  
 $\Box$

The proof of the existence of $\inf\,S$ in the case when $S$ has no minimal element is  absolutely similar. We omit details, but there is another proof of it. Putting $S^*=\{1_{(0,\alpha)}-X_{(0,\alpha)}|\,\,X_{(0,\alpha)}\in S\}$, we obtain $S^*\subset R_\alpha|{_{[0,1]}}$ and $S^*$ has no maximal element. By the above proof, there exists a smallest upper bound $M_{(0,\alpha)}$ of $S^*$. If we put $m_{(0,\alpha)}=1_{(0,\alpha)}-M_{(0,\alpha)}$, then it is a greatest lower bound of $S$. (For the meaning of $1_{(0,\alpha)}-X_{(0,\alpha)}$ see  the next paragraph.) $\Box$

{\bf Theorem 6.} {\it The covering dimension  of a topological space $R_\alpha|{_{[0,1]}}$  in the order topology is equal to $1$, i.e., $dim\,\, R_\alpha|{_{[0,1]}}=1$. }

{\bf Proof.}
It is well known that every linearly ordered space $X$ is hereditarily normal (\cite{l37} Bourbaki, [1948]). 

For every normal space $X$ $dim\,X\leq Ind\, X$ (\cite{l38} Vedenissoff, [1939]). 

For every space $Y$ the properties $dim\,Y=0$ and $Ind\,Y=0$ are equivalent and have as their consequence the normality of $Y$ (\cite{l193} Chap. II, \S 3, Proposition 1, p. 170). 

It is also known that if every hereditarily normal space $X$ is a union of two spaces $Y$ and $Z$ such that $Ind\,Y= 0$ and $Ind\,Z= 0$, then
$Ind\,X\leq 1$ (\cite{l39} Kat\v{e}tov, [1951]). 

Putting $Y=Q_\alpha|_{[0,1]}$ and $Z=R_\alpha|_{[0,1]}\setminus Y$, we notice that since $Y$ and $Z$ are dense in themselves and in $R_\alpha|_{[0,1]}$, we conclude that $ind\, Y=0$ and $ind\, Z=0$; 
therefore, they are hereditarily disconnected (\cite{l40} Hausdorff, [1914]), and being linearly ordered, are strongly zero-dimensional (\cite{l36} Herrlich, [1965]), which is equivalent to $Ind\, Y=0$ and $Ind\, Z=0$ (\cite{l41} Engelking, [1977]). Then for $R_\alpha|_{[0,1]}=Y\cup Z$ we obtain $Ind\,R_\alpha|_{[0,1]}\leq 1$. 

Clearly, since $R_\alpha|_{[0,1]}$ is continuous then $ind\, R_\alpha|_{[0,1]}=1$ and hence, by $1=ind\,R_\alpha|_{[0,1]}\leq Ind\, R_\alpha|_{[0,1]}=1$ and $dim\,R_\alpha|_{[0,1]}\leq Ind\, R_\alpha|_{[0,1]}$, we obtain a desired equality $dim\, R_\alpha|_{[0,1]}=1$. Otherwise, if $dim\, R_\alpha|_{[0,1]}=0$, then $Ind\, R_\alpha|_{[0,1]}=0$, which is false, because $Ind\, R_\alpha|_{[0,1]}=1$. $\Box$ 

The nature of the $1$-dimensional manifold $R_\alpha|_{[0,1]}$ (a generalized real number unit interval) which comes from the  real number unit interval $R_\omega|_{[0,1]}$ is illustrated by the following figures, where $\alpha=\omega\cdot 2$. The first figure is a bifurcation of a rational number, for a example, $x=\frac{1}{2}$:
$$
$$
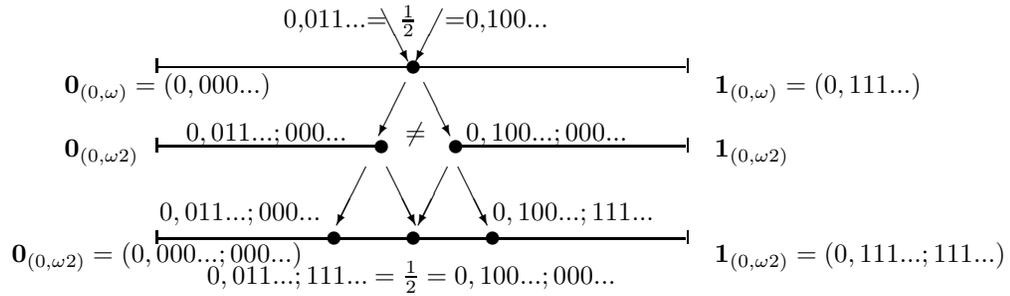
\begin{picture}(100,80)(-110,0)
\put(52,60){\circle*{5}}
\put(52,-5){\circle*{5}}
\put(22,-5){\circle*{5}}
\put(82,-5){\circle*{5}}
\put(40,30){\circle*{5}}
\put(68,30){\circle*{5}}
\put(-45,60){\line(1,0){200}}
\put(-45,58){\line(0,1){5}}
\put(156,58){\line(0,1){5}}
\put(-80,50){${\bf 0}_{(0,\omega)}=(0,000...)$}
\put(166,50){${\bf 1}_{(0,\omega)}=(0,111...)$}
\put(-45,-5){\line(1,0){200}}
\put(-45,-7){\line(0,1){5}}
\put(156,-7){\line(0,1){5}}
\put(-100,-14){${\bf 0}_{(0,\omega 2)}=(0,000...;000...)$}
\put(166,-14){${\bf 1}_{(0,\omega 2)}=(0,111...;111...)$}
\put(-45,30){\line(1,0){85}}
\put(70,30){\line(1,0){85}}
\put(-45,28){\line(0,1){5}}
\put(156,28){\line(0,1){5}}
\put(-80,26){${\bf 0}_{(0,\omega 2)}$}
\put(166,26){${\bf 1}_{(0,\omega 2)}$}
\put(40,82){\vector(1,-2){10}}
\put(56,54){\vector(1,-2){10}}
\put(63,82){\vector(-1,-2){10}}
\put(49,54){\vector(-1,-2){10}}
\put(42,22){\vector(1,-2){11}}
\put(65,22){\vector(-1,-2){11}}
\put(69,22){\vector(1,-2){11}}
\put(34,22){\vector(-1,-2){11}}
\put(47,75){$\frac{1}{2}$}
\put(-26,-22){$0,011...;111...=\frac{1}{2}=0,100...;000...$}
\put(-34,32){$0,011...;000...$}
\put(-44,2){$0,011...;000...$}
\put(49,32){$\not=$}
\put(72,32){$0,100...;000...$}
\put(82,2){$0,100...;111...$}
\put(3,75){0,011...=}
\put(64,75){=0,100...}
\end{picture}

$$
$$
\begin{center}
{\bf Fig. 2}
\end{center}
$$
$$
$$
$$
$$
$$
$$
$$

\parindent=0 cm
and the second figure is a bifurcation of an irrational number, for example, $\frac{1}{\pi}$:

\begin{picture}(200,100)(-110,0)
\put(30,60){\circle*{5}}
\put(62,-5){\circle*{5}}
\put(6,-5){\circle*{5}}
\put(20,30){\circle*{5}}
\put(48,30){\circle*{5}}
\put(-45,60){\line(1,0){200}}
\put(-45,58){\line(0,1){5}}
\put(156,58){\line(0,1){5}}
\put(-80,56){${\bf 0}_{(0,\omega)}$}
\put(166,56){${\bf 1}_{(0,\omega)}$}
\put(-45,-5){\line(1,0){200}}
\put(-45,-7){\line(0,1){5}}
\put(156,-7){\line(0,1){5}}
\put(-80,-10){${\bf 0}_{(0,\omega 2)}$}
\put(166,-10){${\bf 1}_{(0,\omega 2)}$}
\put(-45,30){\line(1,0){65}}
\put(50,30){\line(1,0){106}}
\put(-45,28){\line(0,1){5}}
\put(156,28){\line(0,1){5}}
\put(-80,26){${\bf 0}_{(0,\omega 2)}$}
\put(166,26){${\bf 1}_{(0,\omega 2)}$}
\put(36,54){\vector(1,-2){10}}
\put(29,54){\vector(-1,-2){10}}
\put(50,22){\vector(1,-2){11}}
\put(18,22){\vector(-1,-2){11}}
\put(30,75){$\frac{1}{\pi}=0,01101...$}
\put(-60,32){$0,01101...;000...$}
\put(-64,2){$0,01101...;000...$}
\put(30,32){$\not=$}
\put(68,32){$0,01101...;111...$}
\put(78,2){$0,01101...;111...$}
%\put(60,75){$\pi=$}
\end{picture}

$$
$$
\begin{center}
{\bf Fig. 3}
\end{center}

\parindent=0,5cm
We can also demonstrate a canonical embedding $i_\omega^{\omega\cdot 2}:R_\omega|_{[0,1]}\rightarrow R_{\omega\cdot 2}|_{[0,1]}$ which preserves the linear ordering, given by $i_\omega^{\omega\cdot 2}(X_{(0,\omega)})=Y_{(0,\omega\cdot 2)}$, where $Y_\beta=X_\beta$, $0\leq\beta<\omega$, and $Y_\beta=0$, $\omega\leq\beta<\omega\cdot 2$.

Moreover, for each $\alpha=\omega\cdot\nu$ and $\beta=\omega\cdot\mu$, where $\nu\leq\mu$, there is a canonical embedding $i_\alpha^\beta:R_\alpha|_{[0,1]}\rightarrow R_\beta|_{[0,1]}$ which preserves the linear ordering, given by $i_\alpha^\beta(X_{0,\alpha})=Y_{(0,\beta)}$, where $Y_\gamma=X_\gamma$, $0\leq\gamma<\alpha$, and $Y_\gamma=0$, $\alpha\leq\gamma<\beta$.

We  can also formally consider  the case where $\alpha={\bf \Omega}=\omega\cdot{\bf \Omega}$ and the rational skands $X_{(0,{\bf \Omega})}$, which are eventually $0$ or $1$, form a proper class $Q_{\bf \Omega}|_{[0,1]}$ in $NBG^-$. As to $R_{\bf \Omega}|_{[0,1]}$, i.e., the union of rational and irrational skands $X_{(0,{\bf \Omega})}$, it is not an object in $NBG^-$ because irrational skands are elements of a super-class and outside of $NBG^-$-type theory. Nevertheless, one can consider a canonical embedding $i_\alpha^{\bf \Omega}:R_\alpha|_{[0,1]}\rightarrow R_{\bf \Omega}|_{[0,1]}$ which preserves the linear ordering, given by $i_\alpha^{\bf \Omega}(X_{(0,\alpha)})=Y_{(0,{\bf \Omega})}$, where $Y_\gamma=X_\gamma$, $0\leq\gamma<\alpha$, and $Y_\gamma=0$, $\alpha\leq\gamma<{\bf \Omega}$. Moreover, $Q_{\bf \Omega}|_{[0,1]}=\bigcup\limits_{\alpha<{\bf \Omega}}i^{\bf \Omega}_\alpha(R_\alpha|[0,1])$, where $\alpha=\omega\cdot\nu$, $1\leq\nu<{\bf \Omega}$.

One can also prove that $Q_{\bf \Omega}|_{[0,1]}$ and $R_{\bf \Omega}|_{[0,1]}$ are dense in themselves; $Q_{\bf \Omega}|_{[0,1]}$ is dense in $R_{\bf \Omega}|_{[0,1]}$ and $R_{\bf \Omega}|_{[0,1]}$ is continuous.  It is also clear that $ind\,Q_{\bf \Omega}|_{[0,1]}=0$ and $ind\,R_{\bf \Omega}|_{[0,1]}=0$; therefore, they are hereditarily disconnected, and being linearly ordered, are strongly zero-dimensional,  which is equivalent to $Ind\, Y=0$ and $Ind\, Z=0$. Then for $R_{\bf \Omega}|_{[0,1]}=Y\cup Z$ we obtain $Ind\,R_{\bf \Omega}|_{[0,1]}\leq 1$.  
Clearly, since $R_{\bf \Omega}|_{[0,1]}$ is continuous, then $ind\, R_{\bf \Omega}|_{[0,1]}=1$ and hence, by $1=ind\,R_{\bf \Omega}|_{[0,1]}\leq Ind\, R_{\bf \Omega}|_{[0,1]}=1$ and $dim\,R_{\bf \Omega}|_{[0,1]}\leq Ind\, R_{\bf \Omega}|_{[0,1]}$, we obtain a desired equality $dim\, R_{\bf \Omega}|_{[0,1]}=1$. Otherwise, if $dim\, R_{\bf \Omega}|_{[0,1]}=0$, then $Ind\, R_{\bf \Omega}|_{[0,1]}=0$, which is false, because $Ind\, R_{\bf \Omega}|_{[0,1]}=1$. 

Unfortunately, all these natural arguments cannot be applied, because all references are to results valid only for {\it topological spaces} which are {\it sets}, and there are no similar results (even definitions of topology, dimensions, normality, etc.) for point {\it proper classes}, and it is not clear how to pass this gap. We are in a situation where something is evident, but there are no resources to prove it. We can only formulate the following conjecture.

{\bf Conjecture 1.} {\it $R_{\bf \Omega}|_{[0,1]}$ is a maximally dense one-dimensional continuum whose elements $($mostly hyper-classes$)$ are limits of $Q_{\bf \Omega}|_{[0,1]}$, i.e., limits of ${\bf \Omega}$-sequences whose terms are sets.}
$\Box$

\bigskip

\parindent=0 cm
{\bf 12. Additive and multiplicative operations on generalized fractions}

\parindent=0,5cm

In the next paragraph we shall formally extend for each ordinal $\alpha=\omega\nu$, $1\leq\nu\leq{\bf \Omega}$, the closed unit $R_\alpha|_{[0,1]}$ to the set (hyper-class in the case $\nu={\bf\Omega}$) $R_\alpha$ of generalized real numbers which is a one-dimensional linearly ordered continuous dense homogeneous point set (a dense homogeneous point  manifold in the sense of Cantor). Unfortunately, only for $\nu=1$, i.e. $\alpha=\omega$, $R_\omega$ is  supplied with addition and multiplication which are associative, commutative and distributive; moreover, $R_\omega={\mathbb R}$ is the field of real numbers. If $\nu>1$, then it is not easy    to completely define  addition and multiplication in $R_\alpha$, although it has been tried in \cite{l20}. We introduce here our version, which is different both in notations and in intention from that in \cite{l20}.

For ordinal numbers there is a commutative and associative {\it natural sum} and a {\it natural product} in the sense of Hessenberg (\cite{l24}, 591-594); i.e., if ordinal numbers $\xi$ and $\eta$ are represented in the form of normal expansion $\xi=\omega^{\xi_1}n_1+\omega^{\xi_2}n_2+...+\omega^{\xi_r}n_r$ and $\eta=\omega^{\xi_1}m_1+\omega^{\xi_2}m_2+...+\omega^{\xi_r}m_r$, respectively, where $\xi_1>\xi_2>...>\xi_r$ are ordinal numbers, $n_1,n_2,...,n_r$  and $m_1,m_2,...,m_r$ are integers $\geq 0$, then by definition the natural sum is the following ordinal number:
\begin{equation}
\label{f090}
\xi\oplus\eta=\omega^{\xi_1}(n_1+m_1)+\omega^{\xi_2}(n_2+m_2)+...+\omega^{\xi_r}(n_r+m_r).
\end{equation}

In order to define the natural product $\alpha\odot\beta$ of the ordinal numbers $\alpha$ and $\beta$ one multiplies their normal expansions as if they were polynomials of variable $\omega$; multiplying two powers of number $\omega$, one forms the natural sum of the exponents and arranges the terms obtained from the multiplication according to decreasing exponents. 

(Notice that instead of $\omega$ we can consider any ordinal number $\gamma>1$ as a basis of normal expansion of $\xi$ with $0\leq n_i<\gamma$, $1\leq i\leq r$, see \cite{l5}, Chap. XII, \S 7, e.g., $\gamma=2$; we use shall it below.)

Thus, the problem of the existence of the desired algebraic operations on $R_\alpha$ concerns only generalized real fractions of $R_\alpha|_{[0,1]}$.

So {\it addition} \grqq$+$", {\it subtraction} \grqq$-$", {\it multiplication} \grqq$\cdot$" and {\it division} \grqq$/$" are not defined for {\it  all} generalized real fractions $X_{(0,\alpha)}$ and $Y_{(0,\alpha)}$ in $R_\alpha|_{[0,1]}$ but are defined for {\it some} of them. Here are typical cases of such a possibility:

 $1$). For every $X_{(0,\alpha)},Y_{(0,\alpha)}\in R_\alpha|_{[0,1]}$, such that $X_{\alpha'}\leq Y_{\alpha'}$, $0\leq\alpha'<\alpha$, we put  $Y_{(0,\alpha)}-X_{(0,\alpha)}=Z_{(0,\alpha)}$, where $Z_{\alpha'}=Y_{\alpha'}- X_{\alpha'}$, $0\leq\alpha'<\alpha$. E.g., $X_{(0,\alpha)}-X_{(0,\alpha)}=0_{(0,\alpha)}$ or, for each $X_{(0,\alpha)}\in R_\alpha|_{[0,1]}$, $1_{[0,1]}- X_{(0,\alpha)}$ is always well defined.
 
If $Y_{(0,\alpha)}=\frac{1}{2^{\alpha'}}$ and $X_{(0,\alpha)}=\frac{1}{2^{\alpha''}}$, $0\leq\alpha'<\alpha''<\alpha$, the subtraction $\frac{1}{2^{\alpha'}}-\frac{1}{2^{\alpha''}}$  is well defined. Indeed, in the case when $\alpha'$, $\alpha''$ are ordinal numbers of the 2nd kind $\frac{1}{2^{\alpha'}}-\frac{1}{2^{\alpha''}}=Z_{(0,\alpha)}$, where $Z_\beta=1$, for $\alpha'\leq\beta<\alpha''$, and $Z_\beta=0$, for $\beta<\alpha'$ or $\beta\geq\alpha''$; in the case when $\alpha'$, $\alpha''$ are ordinal numbers of the $1$st kind $\frac{1}{2^{\alpha'}}-\frac{1}{2^{\alpha''}}=Z_{(0,\alpha)}$, where $Z_\beta=1$, for $\alpha'<\beta\leq\alpha''$, and $Z_\beta=0$, for $\beta\leq\alpha'$ or $\beta>\alpha''$; and at last in the case when $\alpha'$ is an ordinal of the $1$st kind and $\alpha''$ is an ordinal of the 2nd kind $\frac{1}{2^{\alpha'}}-\frac{1}{2^{\alpha''}}=Z_{(0,\alpha)}$, where $Z_\beta=1$, for $\alpha'<\beta<\alpha''$, and $Z_\beta=0$, for $\beta\leq\alpha'$ or $\beta\geq\alpha''$. (Note that we use another notation in an appropriate case  for twins.) 

We will say that the interval $[X_{(0,\alpha)},Y_{(0,\alpha)}]$ is of length $\frac{1}{2^{\alpha'}}$,  if $Y_{(0,\alpha)}-X_{(0,\alpha)}=\frac{1}{2^{\alpha'}}$, $0\leq\alpha'<\alpha$.
 
  $2$). We can define $X_{(0,\alpha)}+Y_{(0,\alpha)}$ in the case where, for each $0\leq\alpha'<\alpha$, $X_{\alpha'}$ and $Y_{\alpha'}$ are both $0$, or one of them is $1$ and another one is $0$. Then the result is $Z_{(0,\alpha)}=X_{(0,\alpha)}+Y_{(0,\alpha)}$ such that $Z_{\alpha'}=X_{\alpha'}+Y_{\alpha'}$, $0\leq\alpha'<\alpha$. In particular, if $Z_{(0,\alpha)}=Y_{(0,\alpha)}-X_{(0,\alpha)}$, then $X_{(0,\alpha)}+Z_{(0,\alpha)}=Z_{(0,\alpha)}+X_{(0,\alpha)}=Y_{(0,\alpha)}$.  
  
  If $\alpha'$ is an ordinal of the $1$st kind, then $\frac{1}{2^{\alpha'}}+\frac{1}{2^{\alpha'}}=\frac{1}{2^{\alpha'-1}}$, $0\leq\alpha'<\alpha$; and if  $1\leq\alpha'<\alpha''<\alpha$, are arbitrary ordinal numbers of the $1$st kind, then $\frac{1}{2^{\alpha'}}+\frac{1}{2^{\alpha''}}=\frac{1}{2^{\alpha''}}+\frac{1}{2^{\alpha'}}$ is $Z_{(0,\alpha)}$ such that $Z_\beta=1$, for $\beta=\alpha',\alpha''$; otherwise $Z_\beta=0$. If $0\leq\alpha''<\alpha$ is an ordinal of the $2$nd kind and $0\leq\alpha'<\alpha''<\alpha$ is an ordinal of the $1$st kind, then $\frac{1}{2^{\alpha'}}+\frac{1}{2^{\alpha''}}=\frac{1}{2^{\alpha''}}+\frac{1}{2^{\alpha'}}=Z_{(0,\alpha)}$ such that $Z_\beta=1$, for $\beta=\alpha'$ and $\beta\geq\alpha''$; otherwise, $Z_\beta=0$. 
  
  On the other hand, there are no magnitudes in 
  $R_\alpha|_{[0,1]}$ such as $\frac{1}{2^{\alpha'}}+\frac{1}{2^{\alpha'}}$, for every $0\leq\alpha'<\alpha$
which is an ordinal number of the $2$nd kind, because there is no ordinal number $\alpha'-1$; as well as   $\frac{1}{2^{\alpha'}}+\frac{1}{2^{\alpha''}}$, if $\alpha'$ is an ordinal number of the $2$nd kind and $\alpha'<\alpha''$. 

Moreover, we can define addition $X_{(0,\alpha)}+Y_{(0,\alpha)}$ in all cases which avoid the addition of components of the latter cases when sums of dyadic fractions do not exist.

  $3$). For $1\leq\alpha'<\alpha$, where $\alpha'$ is an ordinal number of the $1$st kind, $\sum\limits_{\alpha'\leq\beta<\alpha'+\omega}
  \frac{1}{2^{\beta}}=\frac{1}{2^{\alpha'-1}}-\frac{1}{2^{\alpha'+\omega}}$, where by an infinite sum we understand the supremum of finite sums, if of course they are well defined. Indeed, $\frac{1}{2^{\alpha'}}+\frac{1}{2^{\alpha'+1}}+...+\frac{1}{2^{\alpha'+n}}+...=\sup\limits_n(\frac{1}{2^{\alpha'}}+\frac{1}{2^{\alpha'+1}}+...+\frac{1}{2^{\alpha'+n}})=\sup\limits_n[(\frac{1}{2^{\alpha'-1}}-\frac{1}{2^{\alpha'}})+(\frac{1}{2^{\alpha'+1}}-\frac{1}{2^{\alpha'}})+...+(\frac{1}{2^{\alpha'+n-1}}-\frac{1}{2^{\alpha'+n}})]=\sup\limits_n(\frac{1}{2^{\alpha'-1}}-\frac{1}{2^{\alpha'+n}})=\frac{1}{2^{\alpha'-1}}-\inf\limits_n\frac{1}{2^{\alpha'+n}}=\frac{1}{2^{\alpha'-1}}-\frac{1}{2^{\alpha'+\omega}}$.

  $4$). We also put $\frac{1}{2^{\alpha''}}\cdot\frac{1}{2^{\alpha'}}=\frac{1}{2^{\alpha''\bigoplus\alpha'}}$, $0\leq\alpha',\alpha''<\alpha$.

 $5$). Each $X_{(0,\alpha)}\in R_\alpha|_{[0,1]}$ can be divided by $2$. Indeed, it is an immediate consequence of the following propositions.

{\bf Proposition 19.} The following formula 
\begin{equation}
\label{f101}
\frac{1}{2^{\alpha'}}=\frac{1}{2^{\alpha'+1}}+\frac{1}{2^{\alpha'+2}}+...+\frac{1}{2^{\alpha'+\omega+1}}+...+\frac{1}{2^{\beta}}+...=\sum\limits_{\alpha'<\beta<\alpha}\frac{1}{2^{\beta}}
\end{equation}
holds, where  summing is given by all $\beta\not=\omega\eta$, $1\leq\eta<\nu$, $\alpha=\omega\nu$, $\nu\geq 1$, and each infinite  sum is by definition the supremum of successive sums of smaller powers: finite and countable (more precisely, $\omega$) terms. 

The {\bf Proof} is an immediate consequence of the above definitions, i.e., sums of different dyadic fractions. Meanwhile, the convergent series whose indexes are greater in general than countable ordinals are interesting in their own right.
\begin{equation}
\label{f102}
\begin{array}{l}
\sum\limits_{\alpha'<\beta<\alpha}\frac{1}{2^{\beta}}=\frac{1}{2^{\alpha'+1}}+\frac{1}{2^{\alpha'+2}}+...+\frac{1}{2^{\alpha'+\omega+1}}+...+\frac{1}{2^{\beta-1}}+...\stackrel{def}{=}\\\sup\limits_{\alpha'<\beta<\alpha}(\frac{1}{2^{\alpha'+1}}+\frac{1}{2^{\alpha'+2}}+...+\frac{1}{2^{\alpha'+\omega+1}}+...+\frac{1}{2^{\beta-1}})=\\\sup\limits_{\alpha'<\beta<\alpha}[(\frac{1}{2^{\alpha'}}-\frac{1}{2^{\alpha'+1}})+...+(\frac{1}{2^{\alpha'+\omega}}-\frac{1}{2^{\alpha'+\omega+1}})+...+(\frac{1}{2^{\beta-1}}-\frac{1}{2^{\beta}})]=\\\sup\limits_{\alpha'<\beta<\alpha}(\frac{1}{2^{\alpha'}}-\frac{1}{2^{\beta}})=\frac{1}{2^{\alpha'}}-\inf\limits_{\alpha'\leq\beta<\alpha}=\frac{1}{2^{\alpha'}},
\end{array}
\end{equation}
where  $\beta\not=\omega\eta$, $1\leq\eta<\nu$, because evidently $\inf\limits_{\alpha'<\beta<\alpha}\frac{1}{2^{\beta}}=0_{(0,\alpha)}$. $\Box$

{\bf Proposition 20.} For each $X_{(0,\alpha)}\in R_\alpha|_{[0,1]}$ 
\begin{equation}
\label{f103}
X_{(0,\alpha)}=\sum\limits_{\beta}\frac{1}{2^{\beta}},
\end{equation}
where the summation is taken over all $1\leq\beta<\alpha$ such that $\beta\not=\omega\eta$, $1\leq\eta<\nu$, $\alpha=\omega\nu$, $\nu\geq 1$, and  $(\frac{1}{2^\beta})_{\alpha'}=1$ iff $X_{\alpha'}=1$, $0\leq\alpha'<\alpha$.

The {\bf Proof} is an immediate consequence of the above definitions. 

Thus,  we can divide each $X_{(0,\alpha)}\in R_\alpha|_{[0,1]}$ by $2$, changing each term $\frac{1}{2^{\beta}}$ in $(\ref{f103})$ by $\frac{1}{2^{\beta+1}}$ and adding summands $\frac{1}{2^{\omega\eta+1}}$, if  $X_\beta=1$, $\alpha'\leq\beta<\alpha'+\omega\nu$, for some $0\leq\alpha'<\alpha$, where $\alpha'$,   is an ordinal number of the $2$nd kind such that $\alpha'+\omega\nu<\alpha$, $1\leq\eta<\nu$, and the sum of all changed terms is the result $Z_{(0,\alpha)}=X_{(0,\alpha)}/2=X_{(0,\alpha)}\cdot\frac{1}{2}$. $\Box$

{\bf Lemma 6.} Each interval $[X_{(0,\alpha)},Y_{(0,\alpha)}]$ of length $\frac{1}{2^{\alpha'}}$, $0\leq\alpha'<\alpha$, can be halved.

{\bf Proof.} In fact, $Y_{(0,\alpha)}-X_{(0,\alpha)}=\frac{1}{2^{\alpha'}}$ or $Y_{(0,\alpha)}=X_{(0,\alpha)}+\frac{1}{2^{\alpha'}}=X_{(0,\alpha)}+\frac{1}{2^{\alpha'+1}}+\frac{1}{2^{\alpha'+1}}$. Consequently, $Y_{(0,\alpha)}-\frac{1}{2^{\alpha'+1}}=X_{(0,\alpha)}+\frac{1}{2^{\alpha'+1}}$ and hence

\parindent=0 cm $[X_{(0,\alpha)},X_{(0,\alpha)}+\frac{1}{2^{\alpha'}}]\cup[Y_{(0,\alpha)}-\frac{1}{2^{\alpha'}},Y_{(0,\alpha)}]=[X_{(0,\alpha)},Y_{(0,\alpha)}]$. Moreover, lengths of $[X_{(0,\alpha)},X_{(0,\alpha)}+\frac{1}{2^{\alpha'}+1}]$ and $[Y_{(0,\alpha)}-\frac{1}{2^{\alpha'+1}},Y_{(0,\alpha)}]$ are equal to $\frac{1}{2^{\alpha'+1}}=\frac{1}{2^{\alpha'}}\cdot\frac{1}{2}$. $\Box$

\parindent=0,5cm
Although there are no  magnitudes in 
  $R_\alpha|_{[0,1]}$ like $\frac{n}{2^{\alpha'}}$, for each $2\leq n<\omega$, where $0\leq\alpha'<\alpha$ is an ordinal number  of the $2$nd kind, there are magnitudes in 
  $R_\alpha|_{[0,1]}$ of multiplications $\frac{1}{2^\beta}\cdot 2^{\alpha'}$ for some ordinals $1\leq\beta,\alpha'<\alpha=\omega\nu$, $\nu\geq 1$, e.g., $\frac{1}{2^\omega}\cdot 2^\omega=1_{(0,\alpha)}$, which is really unexpected.

{\bf Proposition 21.} For each $\alpha',\alpha''\in{\bf On}$, $0\leq\alpha',\alpha''<\alpha=\omega\nu$, $\nu\geq 1$, the following formula 
\begin{equation}
\label{f104}
\frac{1}{2^{\alpha'}}=\frac{1}{2^{\alpha'\oplus\omega''}}\cdot 2^{\alpha''}
\end{equation}
holds; in particular, $\frac{1}{\omega}\cdot\omega=\frac{1}{2^\omega}\cdot 2^\omega=1_{(0,\alpha)}$.

{\bf Proof}. Since $\frac{1}{2^{\alpha'}}\cdot\frac{1}{2^{\alpha''}}\stackrel{def}{=}\frac{1}{2^{\alpha'\oplus\omega''}}$ it is sufficient to show that $\frac{1}{2^{\alpha''}}\cdot2^{\alpha''}=1$. By Proposition 19, we obtain the following identities:
\begin{equation}
\label{f105}
\begin{array}{l}
1=\frac{1}{2^0}=\frac{1}{2}+\frac{1}{2^2}+\frac{1}{2^3}+...=(\frac{1}{2^2}+\frac{1}{2^3}+...)+(\frac{1}{2^2}+\frac{1}{2^3}+...)=\\
(\frac{1}{2^2}+\frac{1}{2^3}+...)\cdot 2=\frac{1}{2}\cdot 2=[(\frac{1}{2^3}+\frac{1}{2^4}+...)+(\frac{1}{2^3}+\frac{1}{2^4}+...)]\cdot 2=\\
(\frac{1}{2^3}+\frac{1}{2^4}+...)\cdot 2^2=\frac{1}{2^2}\cdot 2^2=...=
[(\frac{1}{2^{n+1}}+\frac{1}{2^{n+2}}+...)+\\
(\frac{1}{2^{n+1}}+\frac{1}{2^{n+2}}+...)]\cdot 2^{n-1}=
(\frac{1}{2^{n+1}}+\frac{1}{2^{n+2}}+...)\cdot 2^{n-1}+\\
(\frac{1}{2^{n+1}}+\frac{1}{2^{n+2}}+...)\cdot 2^{n-1}=(\frac{1}{2^{n+1}}+\frac{1}{2^{n+2}}+...)\cdot 2^n=
\frac{1}{2^n}\cdot 2^n=\\(\frac{1}{2^n}-\frac{1}{2^\omega}+\frac{1}{2^{\omega+1}}+\frac{1}{2^{\omega+2}}+...)\cdot 2^n=...=
(\frac{1}{2^\omega}-\frac{1}{2^\omega}+\frac{1}{2^{\omega+1}}+\frac{1}{2^{\omega+2}}+...)\cdot 2^\omega=\\
(\frac{1}{2^{\omega+1}}+\frac{1}{2^{\omega+2}}+...)\cdot 2^\omega=\frac{1}{2^\omega}\cdot 2^\omega=...=
(\frac{1}{2^{\alpha''+1}}+\frac{1}{2^{\alpha''+2}}+...)\cdot 2^{\alpha''}=\frac{1}{2^{\alpha''}}\cdot 2^{\alpha''}.
\end{array}
\end{equation}

 We also need the following lemma.
  
  {\bf Lemma 7.} For every $X_{(0,\alpha)}<Y_{(0,\alpha)}$ in $R_\alpha|_{[0,1]}$ there are  $X'_{(0,\alpha)}<Y'_{(0,\alpha)}$ in $Q_\alpha|_{[0,1]}$ such that $X_{(0,\alpha)}<X'_{(0,\alpha)}$, $Y'_{(0,\alpha)}<Y_{(0,\alpha)}$ and $Y'_{(0,\alpha)}-X'_{(0,\alpha)}=\frac{1}{2^{\alpha'}}$ for some $2\leq\alpha'<\alpha$; in particular, for each $Y_{(0,\alpha)}$ there is $\frac{1}{2^{\alpha'}}$ such that $\frac{1}{2^{\alpha'}}<Y_{(0,\alpha)}$.
  
  {\bf Proof.} Since $X_{(0,\alpha)}<Y_{(0,\alpha)}$ there exists $\alpha'$, $0\leq\alpha'<\alpha$ such that $X_{\alpha'}=0$, $Y_{\alpha'}=1$ and $X_\beta=Y_\beta$, for all $0\leq\beta<\alpha'$. There is a minimal $\alpha''>\alpha'$ such that $X_{\alpha''}=0$; otherwise, $X_\beta=1$, for all $\beta>\alpha'$ and hence $X_{\alpha'}$ should be equal to $1$. Consider $X'_{(0,\alpha)}$ such that $X'_{\alpha''}=1$, $X'_\beta=X_\beta$, $0\leq\beta<\alpha''$, and $X'_\beta=0$, $\beta>\alpha''$. Clearly, $X_{(0,\alpha)}<X'_{(0,\alpha)}<Y_{(0,\alpha)}$. Consider $Y'_{(0,\alpha)}$ such that $Y'_{\alpha''+1}=1$ and $Y'_\beta=X'_\beta$, $0\leq\beta\leq\alpha''$. Clearly, $X'_{(0,\alpha)}<Y'_{(0,\alpha)}<Y_{(0,\alpha)}$ and $Y'_{(0,\alpha)}-X'_{(0,\alpha)}=\frac{1}{2^{\alpha''+1}}$, i.e., $\alpha'=\alpha''+1$.
 $\Box$

{\bf Definition 20}. Let $X_{(0,\alpha)}$ and $Y_{(0,\alpha)}$ be two elements of $R_\alpha|_{[0,1]}$ such that  $X_{(0,\alpha)}\leq Y_{(0,\alpha)}$ and $Y_{(0,\alpha)}-X_{(0,\alpha)}$ is defined. Then the generalized real number $l_{(0,\alpha)}=Y_{(0,\alpha)}-X_{(0,\alpha)}$ is called a {\it length} of the closed interval $[X_{(0,\alpha)},Y_{(0,\alpha)}]$ and of the open interval $(X_{(0,\alpha)},Y_{(0,\alpha)})$.

{\bf Theorem 7.} {\it Let $[X_{(0,\alpha)},Y_{(0,\alpha)}]\supset[X^1_{(0,\alpha)},Y^1_{(0,\alpha)}]\supset...\supset[X^{\alpha'}_{(0,\alpha)},Y^{\alpha'}_{(0,\alpha)}]\supset...$ be a system of embedded closed intervals of $R_\alpha|_{[0,1]}$ such that $\inf\limits_{\alpha'}l^{\alpha'}_{(0,\alpha)}=0_{(0,\alpha)}$; then there exists a unique element $Z_{(0,\alpha)}\in R_\alpha|_{[0,1]}$ such that it belongs to all these intervals, i.e, $\bigcap\limits_{\alpha'}[X^{\alpha'}_{(0,\alpha)},Y^{\alpha'}_{(0,\alpha)}]=Z_{(0,\alpha)}$.}

{\bf Proof.} 
By Theorem 5, there exist  $\bar M_{(0,\alpha)}=\sup\limits_{0\leq\alpha'<\alpha}\{X^{\alpha'}_{(0,\alpha)}\}$ and $\bar m_{(0,\alpha)}=\inf\limits_{0\leq\alpha'<\alpha}\{Y^{\alpha'}_{(0,\alpha)}\}$  in $R_\alpha|_{[0,1]}$. Then $Z_{(0,\alpha)}=\bar M_{(0,\alpha)}=\bar m_{(0,\alpha)}$. Otherwise, by Lemma 7, for $\bar m_{(0,\alpha)}<\bar M_{(0,\alpha)}$ in $R_\alpha|_{[0,1]}$ there are  $X'_{(0,\alpha)}<Y'_{(0,\alpha)}$ in $Q_\alpha|_{[0,1]}$ such that $\bar m_{(0,\alpha)}<X'_{(0,\alpha)}<Y'_{(0,\alpha)}<\bar M_{(0,\alpha)}$ and $Y'_{(0,\alpha)}-X'_{(0,\alpha)}=\frac{1}{2^{\alpha'}}$ for some $2\leq\alpha'<\alpha$, which is in contradiction with the assumption that $\inf\limits_{\alpha'}l^{\alpha'}_{(0,\alpha)}=0_{(0,\alpha)}$, because $[X'_{(0,\alpha)},Y'_{(0,\alpha)}]\subset[X^{\alpha'}_{(0,\alpha)},Y^{\alpha'}_{(0,\alpha)}]$, for each $0\leq\alpha'<\alpha$. $\Box$

 {\bf Theorem 8.} {\it $R_\alpha|_{[0,1]}$ in the order topology is a compact Hausdorff space.}
 
 {\bf Proof.} We already know that the linearly ordered space $R_\alpha|_{[0,1]}$ is normal. Let now $\gamma=\{U_\lambda\}\,\,|\,\,\lambda\in\Lambda$ be an arbitrary covering of $R_\alpha|_{[0,1]}$, consisting of open intervals $U_\lambda$ of $R_\alpha|_{[0,1]}$. We have to prove that there exists a finite subcovering of $\gamma$, which covers $R_\alpha|_{[0,1]}$. Suppose the contrary, and we cannot choose such a finite subcovering. By Lemma 6, we can halve $R_\alpha|_{[0,1]}=[X_{(0,\alpha)},Y_{(0,\alpha)}]$ (evidently, $X_{(0,\alpha)}=0_{(0,\alpha)}$ and $Y_{(0,\alpha)}=1_{(0,\alpha)}$) and choose one $[X^1_{(0,\alpha)},Y^1_{(0,\alpha)}]$ of the parts that cannot be covered by finite elements of $\gamma$. Then we halve $[X^1_{(0,\alpha)},Y^1_{(0,\alpha)}]$ and choose one $[X^2_{(0,\alpha)},Y^2_{(0,\alpha)}]$ of the parts that cannot be covered by finite elements of $\gamma$. We continue this process and conclude that $[X^\omega_{(0,\alpha)},Y^\omega_{(0,\alpha)}]$ cannot be also covered by finite elements of $\gamma$, otherwise, $[X^n_{(0,\alpha)},Y^n_{(0,\alpha)}]$ can be covered by finite elements of $\gamma$, which is in contradiction with our choice. We halve $[X^\omega_{(0,\alpha)},Y^\omega_{(0,\alpha)}]$ and continue our choice for each $0\leq\alpha'<\alpha$. We have gotten a system  $[X_{(0,\alpha)},Y_{(0,\alpha)}]\supset[X^1_{(0,\alpha)},Y^1_{(0,\alpha)}]\supset...\supset[X^{\alpha'}_{(0,\alpha)},Y^{\alpha'}_{(0,\alpha)}]\supset...$  of embedded closed intervals of $R_\alpha|_{[0,1]}$ such that $\inf\limits_{\alpha'}l^{\alpha'}_{(0,\alpha)}=0_{(0,\alpha)}$. Then, by Theorem 7, there exists a unique element $Z_{(0,\alpha)}\in R_\alpha|_{[0,1]}$ such that it belongs to all these intervals, i.e, $\bigcap\limits_{\alpha'}[X^{\alpha'}_{(0,\alpha)},Y^{\alpha'}_{(0,\alpha)}]=Z_{(0,\alpha)}$. Since $\gamma$ is a covering of $R_\alpha|_{[0,1]}$ there exists an element $U_\lambda\in\gamma$ such that $Z_{(0,\alpha)}\in U_\lambda$. Since $\inf\limits_{\alpha'}l^{\alpha'}_{(0,\alpha)}=0_{(0,\alpha)}$ we conclude that there exists an ordinal number $\alpha'<\alpha$ such that $[X^{\alpha'}_{(0,\alpha)},Y^{\alpha'}_{(0,\alpha)}]\subset U_\lambda$ and one-element subcovering $U_\lambda$ of $\gamma$ covers $[X^{\alpha'}_{(0,\alpha)},Y^{\alpha'}_{(0,\alpha)}]$, which is in contradiction with our choice. Thus the assumption that there is no finite subcovering of $\gamma$ which covers $R_\alpha|_{[0,1]}$ is wrong. $\Box$

 \bigskip

\parindent=0 cm
{\bf 13. Generalized real numbers and generalized straight lines}

\parindent=0,5cm
Now we are going to extend $R_\alpha|_{[0,1]}$, $\alpha=\omega\nu$, $\nu\geq 1$, to the set $R_\alpha$ of all generalized real numbers. 

By $R^+_\alpha$  we denote $R_\alpha|_{[0,1]}\cup(R_\alpha|_{(0,1]})^*$, where $(L)^*$ is the backwards linear ordering of $L=R_\alpha|_{(0,1]}=R_\alpha|_{[0,1]\setminus \{0_{(1,\alpha)}}\}$, and identify $1_{(0,\alpha)}$ and  $(1_{(0,\alpha)})^*$ with the ordinal $1$, and  each ordinal $2^{\alpha'},$ $0<\alpha'<\alpha$ with $(\frac{1}{2^{\alpha'}})^*\in L^*$ with the obvious ordering, adding to  the following already defined relations: $Y<Z$ for each $Y\not=1_{(1,\alpha)}$ in $L$ and each $Z\not=(1_{(0,\alpha)})^*$ in $(L)^*$. Putting $R^-_\alpha=(R^+_\alpha)^*$ and denoting $(X_{(0,\alpha)})^*$ by $-X_{(0,\alpha)}$, for every $X_{(0,\alpha)}\in R^+_\alpha$,  identifying $0_{(1,\alpha)}$ and $-0_{(1,\alpha)}$, we define $R_\alpha=R^-_\alpha\cup R^+_\alpha$  as {\it generalized real numbers} with the obvious ordering,  adding to  the following already defined relations: $X_{(0,\alpha)}<Y_{(0,\alpha)}$, for each $X_{(0,\alpha)}\not=-0_{(1,\alpha)}$ in $R^-_\alpha$ and each $Y_{(0,\alpha)}\not=0_{(1,\alpha)}$ in $R^+_\alpha$. It is clear that $R_\alpha$ is a set of the power $2^{|\alpha|}$.

Similarly, we can extend $Q_\alpha|_{[0,1]}$ to the dense subordering subset $Q_\alpha$ of $R_\alpha$ and call it {\it generalized rational numbers}. Its cardinality is $\sum\limits_{\alpha'<\alpha}2^{|\alpha'|}$.

One can easily prove that $Q_\alpha$ is dense in $R_\alpha$; $dim\,\,Q_\alpha=0$, $dim\,\,(R_\alpha\setminus Q_\alpha)=0$, $dim\,\,R_\alpha=1$; $R_\alpha$ is continuous, i.e., for every  bounded set $X_\alpha\subset R_\alpha$,  there exists an interval $[\alpha_0,\alpha_1]$ such that $X_\alpha\subseteq [\alpha_0,\alpha_1]$ has a smallest  upper bound and a greatest lower bound; every closed bounded set is compact; each Dedekind section in $R_\alpha$ has no gap.

In the case $\alpha=\omega^\kappa$, $\kappa\geq 1$ we can represent $R_\alpha$ in a more natural form.

 {\bf Definition 21}. By a set $R^+_\alpha$ of all non-negative generalized numbers we understand
an extended system of embedded curly braces 
\begin{equation}
\label{f999}
...\{_{-\alpha'}...\{_{-1}\{_0\{_1...\{_{\alpha'}...\}\}\}\}\}
\end{equation}  filled by $0$ or $1$ such that for each $X_{(-\alpha,\alpha)}\in R^+_\alpha$, $X_{\alpha'}=1$ only for {\it finite number} indexes, $-\alpha<\alpha'<0$. We consider on $R^+_\alpha$ the lexicographic ordering identifying twins as above. If all non-negative places are filled by $0$ we have the usual ordinal numbers in Cantor's normal form with base $2$ and the lexicographic ordering of them which coincides with the usual one. By a set $R^-_\alpha$ of all non-positive generalized numbers we understand the backwards linear ordering $(R^+_\alpha)^*$, and denote $X^*_{(-\alpha,\alpha)}\in (R^+_\alpha)^*$ by $-X_{(-\alpha,\alpha)}$. At last, by a set $R_\alpha$ of all  generalized numbers we understand $R^-_\alpha\cup R^+_\alpha$ with the natural identification $0^+_{(-\alpha,\alpha)}$ and
$0^-_{(-\alpha,\alpha)}$ and a clear linear ordering. Further we denote by $X_{(-\alpha,\alpha)}$ an arbitrary element of $R_\alpha$ which can be positive, negative or zero.

{\bf Theorem 9.} {\it If $\alpha=\omega^\kappa$, $\kappa\geq 1$, then $R_\alpha$ is topologically and order isomorphic to $R_\alpha$ in the sense of Definition 22.}

{\bf Proof.} We give only a sketch of a proof. It is enough to prove that $R^+_\alpha\setminus R_\alpha|_{[0,1)}$ is isomorphic to $R^+_\alpha\setminus R_\alpha|_{[0,1)}$ in the sense of Definition 22, because $R_\alpha|_{[0,1)}$ are isomorphic in both senses: notice that $R_\alpha|_{(0,1]}$ are also isomorphic in both senses.

Since $R_\alpha|_{(0,1]}$ and $(R_\alpha|_{(0,1]})^*$ are evidently topologically isomorphic it is enough to show that $(R_\alpha|_{(0,1]})^*$ and $\{X_{(-\alpha,\alpha)}\in R_\alpha\,\,\,|\,\,X_{(-\alpha,\alpha)}\geq 1_{(-\alpha,\alpha)}\}$  or $R_\alpha|_{(0,1]}$ and $\{X_{(-\alpha,\alpha)}\in R_\alpha\,\,\,|\,\,X_{(-\alpha,\alpha)}\geq 1_{(-\alpha,\alpha)}\}$ are topologically isomorphic, respectively.

This isomorphism can be defined by the following transfinite induction: the first step is to show that $R_\alpha|_{[\frac{1}{2},1]}$ is isomorphic to $R_\alpha|_{[1,2]}$. It can be done by putting in  correspondence $\frac{1}{2}$ to $2$ and $1$ to $1$, respectively (we obviously simplify the notation). By halving the intervals $[\frac{1}{2},1]$ and $[1,2]$ we put their centers in correspondence to each other, i.e., $\frac{1}{2}+\frac{1}{2^2}$ to $1+\frac{1}{2}$, and do the same (i.e., halving the intervals) with each corresponding interval $[\frac{1}{2},\frac{1}{2}+\frac{1}{2^2}]$ and $[1+\frac{1}{2},2]$ as well as $[\frac{1}{2}+\frac{1}{2^2},1]$ and $[1,1+\frac{1}{2}]$, respectively. Of course, the limit ends will be in this natural correspondence and we continue halving further and further, i.e., $\alpha$ times. It is clear that the closures of the corresponding isomorphic sets of halving are also isomorphic and they coincide with   $R_\alpha|_{[\frac{1}{2},1]}$ and $R_\alpha|_{[1,2]}$, respectively. 

We can do the same with $R_\alpha|_{[\frac{1}{2^2},1]}$ and $R_\alpha|_{[1,2^2]}$, respectively, and notice that a new isomorphism restricted on $R_\alpha|_{[\frac{1}{2},1]}$ and $R_\alpha|_{[1,2]}$, respectively, will coincide with the previous one. The limit isomorphism between $R_\alpha|_{(\frac{1}{2^\omega},1]}$ and $R_\alpha|_{[1,2^\omega)}$ is obvious, which we extend to the isomorphism between $R_\alpha|_{[\frac{1}{2^\omega},1]}$ and $R_\alpha|_{[1,2^\omega]}$. 

In the same manner we show that $R_\alpha|_{[\frac{1}{2^{\omega+n}},\frac{1}{2^\omega}]}$ and $R_\alpha|_{[2^\omega,2^{\omega+n}]}$, $1\leq n<\omega$, are isomorphic. (Only notice here that if $\alpha\not=\omega^\kappa$, $\kappa\geq 1$, then this step would be wrong.) And extend it to the isomorphism between $R_\alpha|_{[\frac{1}{2^{\omega 2}},\frac{1}{2^\omega}]}$ and $R_\alpha|_{[2^\omega,2^{\omega 2}]}$ and hence between $R_\alpha|_{[\frac{1}{2^{\omega 2}},1]}$ and $R_\alpha|_{[1,2^{\omega 2}]}$. The further steps are similar. We omit the details. The resultant ordering isomorphism between $R_\alpha$ in  both senses will be an isomorphism, too. $\Box$

If each element of $R_{(-\alpha,\alpha)}$, $\alpha=\omega_\xi$, $\xi\geq 0$, is considered as a geometric point, then we denote this point set by $L_{(-\alpha,\alpha)}$ and call it a {\it generalized straight line}. If $\xi=0$, then we obtain a classical Euclidean line, which as we know  is uniquely defined by Hilbert's axioms. Moreover, by  one of Euclid's definitions: \grqq The straight line is a line such that it is uniformly arranged towards all its points". It is not so in the cases where $\xi>0$; i.e. there are different kinds of points, e.g., $\xi=1$ and $\frac{1}{2^\omega}$ and $\frac{1}{2^{\omega+1}}$ in $R_{(-\omega_1,\omega_1)}$, which have a different structure to the right of them, though the same structure to the left. We are going to state the following conjecture.

{\bf Conjecture 2.} {\it There exists a system of axioms including a generalized Archimedian axiom, a generalized Cantor's axiom of continuity, which uniquely defines the generalized straight line $L_{(-\omega_\xi,\omega_\xi)}$, $\xi\geq 1$, with an isomorphism $\varphi:L_{(-\omega_\xi,\omega_\xi)}\rightarrow R_{(-\omega_\xi,\omega_\xi)}$.}

In favour of this conjecture says the following corollary of Proposition 21:

{\bf Corollary 6.} For each generalized dyadic fraction $a=\frac{1}{2^{\alpha'}}$, $\alpha'\in{\bf On}$, $1\leq\alpha'<\alpha=\omega_\xi$, $\xi\geq 1$, and arbitrary positive generalized real number $b\in R^+_{(-\omega_\xi,\omega_\xi)}$ there exists an ordinal number $\nu$, $1\leq\nu<\omega_\xi$, such that $a\cdot\nu>b$.

{\bf Proof.} Indeed, there exists an ordinal number $\beta$, $0\leq\beta<\alpha=\omega_\xi$, such that $2^\beta>[b]$, where $[b]$ is an integral part of $b$. Putting $\alpha''=\alpha'\oplus\beta$, we obtain $\nu=2^{\alpha'\oplus\beta}$ which, by Proposition 21, satisfies the desired condition, i.e., $a\cdot\nu=\frac{1}{2^{\alpha'}}\cdot\nu=\frac{1}{2^{\alpha'}}\cdot 2^{\alpha'\oplus\beta}=\frac{1}{2^{\alpha'}}\cdot 2^{\alpha'}\odot 2^\beta=(\frac{1}{2^{\alpha'}}\cdot 2^{\alpha'})\cdot 2^\beta=1\cdot 2^\beta>b$. $\Box$

Note also that the geometry of such straight lines is different from the classical one; i.e., $L_{(-\omega,\omega)}$. E.g., in the generalized plane, i.e., $L^2_{(-\omega_\xi,\omega_\xi)}=L_{(-\omega_\xi,\omega_\xi)}\times L_{(-\omega_\xi,\omega_\xi)}$, $\xi\geq 1$, it is not a case that each of two different points of $R^2_{(-\omega_\xi,\omega_\xi)}$ belongs to a generalized line in $R^2_{(-\omega_\xi,\omega_\xi)}$. For example, there are straight lines in $R^2_{(-\omega_\xi,\omega_\xi)}$ such as $y=x$, or $y=-x$, or $x=const$, or $y=const$, but there is no line given by the following equation: $y=2x$. So for possible straight lines in $R^2_{(-\omega_\xi,\omega_\xi)}$, $\xi\geq 1$, a generalized version of Zeno's paradox arises; however in the direction $y=2x$, Zeno's arrow does not even exist; it is totally destroyed. $\Box$
\bigskip

\parindent=0 cm
{\bf 14. Elements of a generalized calculus}

\parindent=0,5cm
 {\bf Definition 22}. By an $\alpha'$-sequence in $R_\alpha$, $\omega\leq\alpha'\leq\alpha\leq{\bf \Omega}$, $\alpha=\omega\nu$, $\alpha'=\omega\nu'$, $1\leq\nu,\nu'\leq{\bf \Omega}$, we understand a generalized skand $S_{(0,\alpha')}$ whose components $S_{\beta'}$ are elements of $R_\alpha$, $0\leq\beta'<\alpha'$. We denote it temporarily by $S_{\beta'}$, $0\leq\beta'<\alpha'$; in the usual way its denotation is more complicated, i.e., $\{S_{\beta'}\}|_{\beta'<\alpha'}$ and is called a transfinite sequence of type $\alpha'$ (see the corresponding definitions of course for transfinite sequences of ordinal numbers \cite{l6}, p. 287).

 {\bf Remark 19.} So, for $\alpha>\omega$, there are many converging $\alpha'$-sequences in $R_\alpha$, $\omega\leq\alpha'\leq\alpha$, where $\alpha'=\omega\nu$, $\nu\geq 1$. Indeed, since each generalized real number $X_{(-\alpha,\alpha)}\in R_\alpha$ is the intersection of some $\alpha'$-sequence of embedded closed intervals $[X_{(-\alpha,\alpha)},B^{\beta}_{(-\alpha,\alpha)}]\supset[X_{(-\alpha,\alpha)},B^{\beta+1}_{(-\alpha,\alpha)}]$, $0\leq\beta<\alpha'$, (a cofinal system of closed neighborhoods on the right at $X_{(-\alpha,\alpha)}$) and is at the same time the intersection of some $\alpha''$-sequence of embedded closed intervals $[A^\beta_{(-\alpha,\alpha)},X_{(-\alpha,\alpha)}]\supset[A^{\beta+1}_{(-\alpha,\alpha)},X_{(-\alpha,\alpha)}]$, $0\leq\beta<\alpha''$, $\omega\leq\alpha''\leq\alpha$, (a cofinal system of closed neighborhoods on the left at $X_{(-\alpha,\alpha)}$) then $B^{\beta}_{(-\alpha,\alpha)}$ and $A^\beta_{(-\alpha,\alpha)}$ are such $\alpha'$- and $\alpha''$-sequences, respectively. Notice that there are elements $X_{(-\alpha,\alpha)}\in R_\alpha$ for which $\alpha'=\alpha''=\alpha$ (generalized integers, generalized irrationals, some generalized rationals), but there are elements $X_{(-\alpha,\alpha)}\in R_\alpha$ for which $\alpha'<\alpha$ or $\alpha''<\alpha$ though not at the same time; i.e., if $\alpha'<\alpha$, then $\alpha''=\alpha$ or if $\alpha''<\alpha$, then $\alpha'=\alpha$(e.g., $\frac{1}{2^\omega}$, $1-\frac{1}{2^\omega}$, etc.)
 
 {\bf Definition 23}. A generalized real number $X_{(-\alpha,\alpha)}\in R_\alpha$ is a limit of $\alpha'$-sequence $S_\beta$,  $0\leq\beta<\alpha'\leq\alpha$, notation $X_{(-\alpha,\alpha)}=\lim\limits_{\beta\rightarrow\alpha'}S_\beta$, if for each open interval $(A_{(-\alpha,\alpha)},B_{(-\alpha,\alpha)})$, $A_{-(\alpha,\alpha)}<B_{(-\alpha,\alpha)}$, which contains $X_{(-\alpha,\alpha)}$, there exists an ordinal $0\leq\beta_0<\alpha'$ such that $S_\beta\in (A_{(-\alpha,\alpha)},B_{(-\alpha,\alpha)})$, for all $\beta_0<\beta<\alpha'$.
 
In this case a $\alpha'$-sequence $S_\beta$,  $0\leq\beta<\alpha'\leq\alpha$, is called convergent and $X_{(-\alpha,\alpha)}\in R_\alpha$ is its limit. Clearly,  if the $\alpha'$-sequence $S_\beta$,  $0\leq\beta<\alpha'\leq\alpha$, converges to $X_{(-\alpha,\alpha)}\in R_\alpha$, then this limit is unique.

 Using the classical arguments we can easily prove the following theorems.
 
 {\bf Theorem 10.} {\it A mapping $f:R_\alpha\rightarrow R_\alpha$ is continuous in the ordering topology if and only if for each element $X_{(-\alpha,\alpha)}\in R_\alpha$ and every $\alpha'$-sequence $S_\beta$,  $0\leq\beta<\alpha'\leq\alpha$, such that $\lim\limits_{\beta\rightarrow\alpha'} S_\beta=X_{(-\alpha,\alpha)}$, then $\lim\limits_{\beta\rightarrow\alpha'} f(S_\beta)=f(X_{(-\alpha,\alpha)})$}.
 
 {\bf Proof.} Let $f:R_\alpha\rightarrow R_\alpha$ be a continuous mapping and $X_{(-\alpha,\alpha)}$ an arbitrary element in $R_\alpha$. Consider any open interval $(A'_{(-\alpha,\alpha)},B'_{(-\alpha,\alpha)})$ which contains $f(X_{(-\alpha,\alpha)})$ and find an open interval
$(A_{(-\alpha,\alpha)},B_{(-\alpha,\alpha)})$ which contains $X_{(-\alpha,\alpha)}$ such that $f(A_{(-\alpha,\alpha)},B_{(-\alpha,\alpha)})\subseteq (A'_{(-\alpha,\alpha)},B'_{(-\alpha,\alpha)})$. 
Since $\lim\limits_\beta S_\beta=X_{(-\alpha,\alpha)}$ there exists an ordinal number $\beta_0$ such that $S_\beta\in (A_{(-\alpha,\alpha)},B_{(-\alpha,\alpha)})$, for all $\beta_0<\beta<\alpha'\leq\alpha$. Then evidently $f(S_\beta)\in (A'_{(-\alpha,\alpha)},B'_{(-\alpha,\alpha)})$, for all $\beta_0<\beta<\alpha'\leq\alpha$.

Conversely,  suppose the opposite, i.e.,  for every $\alpha'$-sequence $S_\beta$,  $0\leq\beta<\alpha'\leq\alpha$, such that $\lim\limits_{\beta\rightarrow\alpha'} S_\beta=X_{(-\alpha,\alpha)}$ we have $\lim\limits_{\beta\rightarrow\alpha'} f(S_\beta)=f(X_{(-\alpha,\alpha)})$; however $f$ is not continuous at some point $X_{(-\alpha,\alpha)}$ in
$R_{(-\alpha,\alpha)}$. Evidently,  $f$ is not continuous on the right at $X_{(-\alpha,\alpha)}$ or on the left at $X_{(-\alpha,\alpha)}$. We can consider the first case; the second is similar. Thus, if $f$ is not continuous on the right at $X_{(-\alpha,\alpha)}$, then there is an open interval $(A'_{(-\alpha,\alpha)},B'_{(-\alpha,\alpha)})$ which contains $f(X_{(-\alpha,\alpha)})$ such that for every embedded closed interval $[X_{(-\alpha,\alpha)},B^\beta_{(-\alpha,\alpha)}]$,  $0\leq\beta<\alpha'\leq\alpha$, which is cofinal in the system of all neighborhoods on the right at $X_{(-\alpha,\alpha)}$,  we have $f([X_{(-\alpha,\alpha)},B^\beta_{(-\alpha,\alpha)}])$ is not a subset of $(A'_{(-\alpha,\alpha)},B'_{(-\alpha,\alpha)})$. Choosing in each $[X_{(-\alpha,\alpha)},B^\beta_{(-\alpha,\alpha)}]$ an element $S_\beta$, such that $f(S_\beta)\notin (A'_{(-\alpha,\alpha)},B'_{(-\alpha,\alpha)})$, $0\leq\beta<\alpha'\leq\alpha$,  we see that $\lim\limits_{\beta\rightarrow\alpha'}S_\beta=X_{(-\alpha,\alpha)}$ but $\lim\limits_{\beta\rightarrow\alpha'}f(S_\beta)\not=f(X_{(-\alpha,\alpha)})$, which is in contradiction with our assumption; i.e., if for an $\alpha'$-sequence $S_\beta$ we have  $\lim\limits_{\beta \rightarrow\alpha'}S_\beta=X_{(-\alpha,\alpha)}$, then $\lim\limits_\beta f(S_\beta)=f(X_{(-\alpha,\alpha)})$. $\Box$

 {\bf Theorem 11.} {\it Let $X_{(-\alpha,\alpha)}$ and $Y_{(-\alpha,\alpha)}$ be elements of $R_\alpha$. If 
 \begin{equation}
 \label{f998}
 f:[X_{(-\alpha,\alpha)},Y_{(-\alpha,\alpha)}]\rightarrow R_\alpha
 \end{equation}
  be a continuous mapping such that $f(X_{(-\alpha,\alpha)})$ and $f(Y_{(-\alpha,\alpha)})$ have different signs, i.e., $-$ and $+$ or $+$ and $-$, respectively, then there exists an element $Z_{(-\alpha,\alpha)}$ in $R_\alpha$ such that $f(Z_{(-\alpha,\alpha)})=0_{(-\alpha,\alpha)}$}.
 
 {\bf Proof}. If in the interval $[X_{(-\alpha,\alpha)},Y_{(-\alpha,\alpha)}]$ there are no integer ordinal numbers except $X_{(-\alpha,\alpha)}$ and $Y_{(-\alpha,\alpha)}$, then the proof is classical here and is given  by the method of our generalized {\it dichotomy}, or halving the intervals of the unit one. Comp., \cite{l22}, Theorem 25, p. 41. If there are such integers, then choose one of them; let it be $Z_{(-\alpha,\alpha)}$ such that $X_{(-\alpha,\alpha)}<Z_{(-\alpha,\alpha)}<Y_{(-\alpha,\alpha)}$. Clearly, $f(X_{(-\alpha,\alpha)})$ and $f(Z_{(-\alpha,\alpha)})$ or $f(Z_{(-\alpha,\alpha)})$ and $f(Y_{(-\alpha,\alpha)})$ have different signs, and we choose that interval and denote it by $[X^1_{(-\alpha,\alpha)},Y^1_{(-\alpha,\alpha)}]$. We continue this transfinite process up to the first case where there are no integer ordinal numbers except $X^\beta_{(-\alpha,\alpha)}$ and $Y^\beta_{(-\alpha,\alpha)}$, $1\leq\beta<\alpha$; then the proof is also classical  and is given  by the method of our generalized {\it dichotomy} or halving intervals of the unit one. Since we obtain a system of embedded closed intervals  whose lengths converge to $0_{(-\alpha,\alpha)}$, by Theorem 7,
there exists a unique element $Z_{(0,\alpha)}\in R_\alpha$ such that it belongs to all these intervals. By the continuous property of $f$, $f(Z_{(0,\alpha)})$ has at the same time different signs, which is possible only in the  case where $f(Z_{(0,\alpha)})=0_{(-\alpha,\alpha)}$. $\Box$

 One can consider other classical theorems of Mathematical Analysis in the case of $R_\alpha$ since for $\xi=1$ we have the usual real numbers ${\mathbb R}=R_\omega$.
 
{\bf Theorem 12}.  {\it Every closed interval $[A_{(-\alpha,\alpha)},B_{(-\alpha,\alpha)}]$ of $R_\alpha$ is compact.}
 
 {\bf Proof.} Without loss of generality we can assume that $A_{(-\alpha,\alpha)}$ and $B_{(-\alpha,\alpha)}$ are positive or negative ordinals and $A_{(-\alpha,\alpha)}<B_{(-\alpha,\alpha)}$. Moreover, we can restrict the proof to the case where $[A_{(-\alpha,\alpha)},B_{(-\alpha,\alpha)}]\subset R^+_\alpha$. Assume now that $\gamma$ is a covering of $[A_{(-\alpha,\alpha)},B_{(-\alpha,\alpha)}]\subset R^+_\alpha$ consisting of open sets in $R_\alpha$ which has no finite subcovering. Then we divide this interval by the point $A_{(-\alpha,\alpha)}+1$, i.e., $[A_{(-\alpha,\alpha)},B_{(-\alpha,\alpha)}]=[A_{(-\alpha,\alpha)},A_{(-\alpha,\alpha)}+1]\cup[A_{(-\alpha,\alpha)}+1,B_{(-\alpha,\alpha)}]$, the next $[A_{(-\alpha,\alpha)}+1,B_{(-\alpha,\alpha)}]$ by $A_{(-\alpha,\alpha)}+2$, etc. By Theorem 8, restrictions of $\gamma$ on $[A_{(-\alpha,\alpha)},A_{(-\alpha,\alpha)}+1]$, $[A_{(-\alpha,\alpha)}+1,A_{(-\alpha,\alpha)}+2]$,...,$[A_{(-\alpha,\alpha)}+n,A_{(-\alpha,\alpha)}+n+1]$,... have finite subcoverings; we conclude that the restrictions of $\gamma$ on each $[A_{(-\alpha,\alpha)}+n,B_{(-\alpha,\alpha)}]$ have no finite subcoverings, $0\leq n<\omega$. Hence $[A_{(-\alpha,\alpha)}+\omega,B_{(-\alpha,\alpha)}]$ has no finite subcovering of the restriction of $\gamma$, if of course $A_{(-\alpha,\alpha)}+\omega<B_{(-\alpha,\alpha)}$. If
$A_{(-\alpha,\alpha)}+\omega=B_{(-\alpha,\alpha)}$, then $B_{(-\alpha,\alpha)}$ is covered by one element of $\gamma$ and consequently, $[A_{(-\alpha,\alpha)}+n,B_{(-\alpha,\alpha)}]$ covers by this element for some natural $n$. Contradiction. So we continue this transfinite process, passing all possible limit ordinals  and $B_{(-\alpha,\alpha)}$ with the same argument. If $B_{(-\alpha,\alpha)}$ is not a limit ordinal, then the impossibility of choosing a finite subcovering from the restriction of $\gamma$ on $[B_{(-\alpha,\alpha)}-1,B_{(-\alpha,\alpha)}]$ contradicts  Theorem 8. $\Box$

 {\bf Theorem 13}. {\it A subspace X of $R_\alpha$ is a compact Hausdorff space if and only if it is a bounded closed subset of $R_\alpha$}.
 
 The {\bf Proof} is classical, except with reference to Theorem 12.

{\bf Theorem 14}. {\it Every continuous image of a compact Hausdorff space is a compact Hausdorff space.}

The {\bf Proof} is classical.

{\bf Theorem 15}. {\it Every bounded subset $S$ of $R_\alpha$ has a smallest upper bound $M_{(0,\alpha)}=\sup\,S$ and a greatest lower bound $m_{(0,\alpha)}=\inf\,S$ in $R_\alpha$.}

{\bf Proof.} It is a consequence of Theorem 5.

{\bf Theorem 16}. {\it Every continuous function $f$ defined on a compact subset $X$ of $R_\alpha$ is bounded and reaches its maximum and minimum values}.

{\bf Proof.} By Theorem 15, the image $S=f(X)$ is a compact Hausdorff space of $R_\alpha$ and hence, by Theorem 14, bounded. We apply  Theorem 16 and obtain $M_{(-\alpha,\alpha)}=\sup\,X$ and $m_{(-\alpha,\alpha)}=\inf\,X$ which by compactness of $S$ belong to $f(S)$. Since $f:X\rightarrow S$ is a surjection $f$, evidently reaches its maximum and minimum values.

{\bf Remark 20.} If $\alpha=\omega_\xi$, where $\omega_\xi$ is the initial ordinal number, $\xi\geq 1$, then our $R_{\omega_\xi}$ differs from the space denoted by $R_\xi$ in \cite{l23} Hausdorff, [1908], because for the latter $\dim\,R_\xi=0$ and for the former $\dim\,R_{\omega_\xi}=1$.

{\bf Theorem 18}. {\it The class $Q_{\bf \Omega}$ is an $\aleph_\kappa\mbox{-}$universal  linear ordering for all cardinals $\aleph_\kappa$, $\kappa\in{\bf On}$.}

The {\bf Proof} is a consequence of Mendelson's theorem: $Q_\kappa$ is an $\aleph_\kappa\mbox{-}$universal  linear ordering (see \cite{l7} and \cite{l08}, p. 169).

{\bf Theorem 19}. {\it If $X=\bigcup\limits_{1\leq\alpha'<\alpha=\omega_\xi}\frac{1}{2^{\alpha'}}\subset R_\alpha$, $\alpha=\omega_\xi$ is a fixed initial regular  ordinal, $\omega_0\leq\omega_\xi<{\bf \Omega}$, then  there exist proper embeddings $X'\supset X''\supset ... \supset X^{(\beta)}\supset ...$ such that $X^{(\beta)}=\bigcap\limits_{1\leq\gamma<\beta}X^{(\gamma)}$, $\beta=\omega\nu$, $1\leq\nu$, and $X^{(\beta+1)}=(X^{(\beta)})'$, $1\leq\beta<\alpha$, where $X'$ is the Cantor {\it derivative} of $X$, i.e., the set of all limit points of $X$ in the order topology, then $X^{(\alpha)}=0_{(0,\alpha)}$ and hence $X^{(\alpha+1)}=\emptyset$.}

{\bf Proof.} Indeed, $X'=\{\frac{1}{2^\mu}|\mu=\omega\nu, 1\leq\nu\}$, $X''=\{\frac{1}{2^\mu}|\mu=\omega^2\nu, 1\leq\nu\}$,...;$\,\,X^{(\omega)}=\{\frac{1}{2^\mu}|\mu=\omega^{\omega}\nu, 1\leq\nu\}$,...;$\,\,X^{(\varepsilon_2)}=\{\frac{1}{2^\mu}|\mu=\varepsilon_2\nu, 1\leq\nu\}$,...;$\,\,X^{(\alpha)}=\bigcap\limits_{\varepsilon_0\leq\varepsilon_\nu<\varepsilon_\alpha=\alpha}X^{(\varepsilon_\nu)}=0_{(0,\alpha)}$, where $\varepsilon_\nu$, $0\leq\nu\leq\alpha$, are $\varepsilon$-numbers and $X^{(\varepsilon_\alpha+1)}=\emptyset$.

Compare with Cantor's results \cite{l09} of stabilized embeddings, where only the classical real line was considered. $\Box$

\bigskip
\parindent=0 cm
{\bf 15. Conclusion}

\parindent=0,5cm
Skand Theory, as presented in this paper, is actually an axiomatic set theory which extends other axiomatic set theories containing the Foundation Axiom, by representing new mathematical objects, called here \grqq skands", which turn out to be non-well-founded sets, or members of Cantor's paradise, and which were eliminated by von Neuman's regularity axiom.

For formalization we have chosen a von Neumann-Bernays-G\"{o}del-type set theory $NBG$ which includes  the axiom of choice ${\bf C}$  and the axiom of  foundation ${\bf FA}$, where the latter,  within the axiom of choice, is equivalent to the regularity axiom. Skands of length greater than or equal to  $\omega$, which are constituents of the new, extended axiomatic set theory, all sets of which satisfy all axioms of $NBG$, not only enrich von Neumann's universe, but also clarify the binary relation $X\in X$, and reveal the essence of Aczel's anti-foundation axiom. 

Actually, for the formalization of skand theory one  could also use an alternative axiomatic system which contained the foundation axiom, e.g., Zermelo-Frankel axiomatic set theory with the regularity axiom or any other axiomatic system with the foundation axiom. Thus, the skand theory is a definite, economical method for the intake of non-well-founded objects, extending  the pseudo-well-founded sets step by step. On the other hand, skands are generalizations of some (but not all) of Mirimanoff's extraordinary sets. 

The extension of set theory through skands has clarified the relation $X\in X$, which turns out to be essentially many-valued. This ambiguity (skands of length greater than $\omega$) implies an unexpected result: Russell's well-known paradox is not a paradox at all. The same is true of all other set-theoretic paradoxes presented by Cantor, Burali-Forti, Zermelo, Hilbert, Mirimanoff and so on. Moreover within von Neuman-Bernays-G\"{o}del-type set theory $NBG$, one does not need skands to show that Russell's paradox does not obtain, because all sets in this axiomatic theory are well-founded. Having supposed that  Russell's collection $R=\{X|\,\, X\notin X\}$ of sets  is a set, we state that the relation $R\in R$ is always {\it false} in $NBG$, and the relation $R\notin R$ is always {it true} in $NBG$. Hence by the truth-function table values for implication, $R\in R\Longrightarrow R\notin R$ is  {\it true} and $R\notin R\Longrightarrow R\in R$ is  {\it false}. Consequently, Russell's paradox $R\in R\Longleftrightarrow R\notin R$ {\it is false} in $NBG$;  nevertheless, in $NBG$ as well as in $ZFA$ or other axiomatic set theories mathematicians use Russell's argument in proving that $R$ is a proper class, not a set. (Notice that $R$ coincides with the universal class ${\bf V}$ in $NBG$, see \cite{l7}, Chap. 4, \S 5, Proposition 4.40; and \cite{l222}, Chap. II, \ S 6.) So even in $NBG$, the proof of the proposition that $R$ is not a set but a proper class makes no use of Russell's argument. We corrected the proof of this proposition by a very simple principle:

The {\bf Maximality Principle.} {\it If there exists a maximal $($universal$)$ collection $X$ $($sets, classes, hyper-classes$)$, given by some property, predicate, etc., then any assertion which implies the existence of a new element $x$ with the same property and $x\notin X$ is {\bf false}}.

All other set-theoretic paradoxes in Set Theory were diagnosed by the relation $X\in X$ and corrected by this
Maximality Principle. 

Self-similar skands show the existence of the \grqq last cardinal number," called the eschaton, which turns out to be 
a strongly inaccessible (class)-cardinal number that is important for the problem of the existence or non-existence of a strongly inaccessible (set)-cardinal number (an application to the strongly inaccessible cardinal number problem).

Skands as special notations for infinite exponents describe {\it all} $\varepsilon$-numbers in  Cantor's sense (an application to $\varepsilon$-number theory), and as special notations for transfinite sequences skands describe generalized real numbers of arbitrary powers $\aleph_\kappa$, $\kappa\geq 1$, (of course, if we add the Generalized Continuum Hypothesis to the other axiom of set theory; otherwise they describe generalized real numbers of arbitrary  large powers, more precisely, of power $2^{\aleph_\kappa}$, $\kappa\geq 0$). Thus there is another application to  new possible geometries of  the straight line. These investigations go back to Hausdorff \cite{l23}.

We omit another way to construct a hierarchy of generalized real numbers of different powers, which are one-dimensional and continuous but different from these presented above, which we had originally planned  to show in the present paper. Moreover, algebraic operations would be completely defined, as in Conway's approach \cite{l22}, but through a fundamentally different method. There is also a conjecture that Conway's numbers, which come from Game Theory, are subclasses of these new generalized real numbers. Since this paper is already too long, and a construction of generalized real numbers with complete operations is beyond the scope of skand theory, we plan to publish the new material in another paper in the  near future.

\bigskip

\parindent=0 cm
{\bf Acknowledgement}

\parindent=0,5cm
The author wishes to thank Professor Philip T. Grier (USA) for help in editing the manuscript.

\end{document}